\documentclass[12pt]{article}
\usepackage{graphicx}
\usepackage{grffile}

\usepackage{amsmath}
\usepackage{amssymb}

\usepackage{comment}
\usepackage{verbatim}

\usepackage{ifthen}

\newboolean{vmcode}  
\setboolean{vmcode}{false}

\newcommand{\cmdvmcode}[2]{\ifthenelse{\boolean{vmcode}}{\verbatiminput{../Maple/#1}}{#2}}
\def\grpath{.}

\newcommand{\slipparameter}{\beta}

\newcommand{\QQ}{Q}
\newcommand{\usteady}{u}
\newcommand{\EllipticE}{{\rm EllipticE}}
\newcommand{\EllipticK}{{\rm EllipticK}}

\newtheorem{theorem}{Theorem}

\def\grpath{.} 
\begin{document}




\numberwithin{equation}{section}

\title{
Functions with constant Laplacian satisfying Robin boundary conditions on an ellipse}

\author{Grant Keady and Benchawan Wiwatanapataphee\\
Department of Mathematics and Statistics\\
Curtin University, Bentley, 6102\\
Western Australia}
\date{\today}


\maketitle
\vspace{-0.5cm}
\section*{Abstract}
We study
the problem of finding functions, defined within and on an ellipse, whose Laplacian is -1 and which satisfy a homogeneous Robin boundary condition on the ellipse.
The parameter in the Robin condition is denoted by $\beta$.
The integral of the solution over the ellipse, denoted by $Q$,
is a quantity of interest in some physical applications.
The dependence of $Q$ on $\beta$ and the ellipse's geometry is found.
Several methods are used.
\begin{itemize}
\item 
Part I: To find the general solution the boundary value problem
is formulated in elliptic cylindrical coordinates.
Truncated Fourier series solution are then used to achieve good numerical approximations.
The truncation is higher order than earlier papers and also totally systematic.\\
Various asymptotic approximations are found directly from the pde formulations, this being far easier than from our series solution.\\
An ad-hoc approximation that the solutions might have nearly elliptic level curves also leads to estimates.\\
Numerical results from the various approaches are compared.

\item 
Part II. Fourier series are much more directly applicable to a Problem (P($\infty$)),
this problem first arising in connection with asymptotics
as $\beta$ tends to infinity.
This problem is solved by Fourier series and a difference equation which would facilitate rapid computation of the coefficients is studied.

\item
Part III. A variational approach applicable to any cross-section -- 
requiring however information on the solution of Problem (P($\infty$)) --
yields simple and accurate lower bounds.
The method was first
illustrated with the rectangle as an example as for that cross-section Problem (P($\infty$)) 
is easy to solve.
(This material is now published in {\it J. Fluids Engineering}.)

\item 
Part IV. The approximation -- lower bound -- described in Part III is further approximated in the case of
an ellipse.
A variational approach to Problem (P($\infty$))  is used, and connections made to
earlier approximations in Part I and to the Fourier Series in Part II.

\end{itemize}

One notational blemish is that our original usual notation for the perimeter $|\partial\Omega|$
is $L$ but a referee requested, for engineering journals, using $P$.
The latter is used consistently in Part IV.

\noindent
This arXiv preprint will be referenced by journal articles we submit on parts of the work.
The arXiv version contains material,
e.g. codes for calculating $Q$, not in the very much shorter journal version.

\vspace{1cm}

\noindent
Part I is joint work. For the remaining parts GK is the author.
\goodbreak





\newpage

\begin{center}
{\Large{\textbf{\textsc{ Part I:
Approximations by
Fourier series,\\
 asymptotics 
and an ad hoc variational guess
}}}}
\end{center}
\section{Introduction}\label{sec:Intro}

\subsection{
The pde problem}\label{subsec:GovEq}

With $\Omega$ a plane domain and $\beta\ge{0}$ given,
we seek the solution of
$$ - \frac{\partial^2 u}{\partial x^2}- \frac{\partial^2 u}{\partial y^2}
 = 1 \qquad {\rm in\ \ }\Omega ,\qquad
 u + \beta \frac{\partial u}{\partial n}
= 0 \qquad {\rm on\ \ }\partial\Omega .
\eqno({\rm P}(\beta))
$$
Under reasonable conditions on the boundary $\partial\Omega$ (classical) solutions exist and, by Maximum Principle
arguments, are unique and positive in $\Omega$.
\cite{KM93}  surveys some of the easier results.

In the context of flows, this boundary condition is called a  {\it slip boundary condition} or
{\it Navier's boundary condition}; in the wider mathematical literature it is called a {\it Robin boundary condition}.
The $\beta=0$ case is, in the fluid mechanics context, called a `no-slip boundary condition': in elasticity the $\beta=0$ case is called
`the elastic torsion problem'.
A functional of interest, in the context of flows, is the volume flow rate
\begin{equation}
\QQ
:= \int_\Omega u .
\label{eq:Qdef}
\end{equation}
We remark that the same pde problem arises in contexts other than slip flow,
for example, heat-flow with `Newton's law of cooling', e.g.~\cite{KM93,MK94}.
While there are other applications of this pde problem we will, on occasions, use the fluid flow terminology for $u$, $Q$ and $\beta$.

There are many inequalities established for Problem~(P($\beta$)).
There are isoperimetric inequalities, reviewed in~\cite{KWK18} and others, for example, from~\cite{KM93}
\begin{equation}
{\rm max}\left(
\beta\, \frac{|\Omega|^2}{|\partial\Omega|} +\QQ_0,\ 
 \beta\, \frac{|\Omega|^2}{|\partial\Omega|} + \Sigma_\infty +\frac{\Sigma_1}{\beta}\right)
\le \QQ(\beta)
\le \beta\, \frac{|\Omega|^2}{|\partial\Omega|} + \Sigma_\infty 
= U(\beta) , 
\label{in:vBnds}
\end{equation}
where $\Sigma_\infty>0$ and $\Sigma_1\le{0}$  are  defined in equation~(\ref{eq:SigDef1}) at the beginning of~\S\ref{sec:Varl}
and used associated with asymptotics for $\beta$ large.
$U(\beta)$ is defined as the expression immediately to its left.
These old inequalities~(\ref{in:vBnds}) are established using the variational formulation of~\S\ref{sec:Varl}
in~\cite{KM93} : see inequalities~(4.4) and~(4.9).
The inequality (4.4) of~\cite{KM93}  is well-known, (4.9) less so.
The same methods lead to an improvement given in
inequality~(\ref{in:lbbL}).

Except in this subsection and a very small number of introductory subsections, e.g.
\S\ref{subsec:inS0Sinf1}, \S\ref{subsec:VarlGen} and 
\S\ref{subsec:betaSmallGenOm}
where we allow $\Omega$ to be more general,
we will be treating $\Omega$ as the elliptical domains
within ellipses described by
\begin{equation}
\frac{x^2}{a^2}+\frac{y^2}{b^2} \leq 1\qquad{\rm with\ \ } a\geq b .
\label{eq:ellab}
\end{equation}
In Part III we apply the same methods to rectangular domains.

In this Part, we derive, for elliptical domains,
a Fourier series form of the analytical solution of the problem
and provide many checks on it.
Others, e.g.~\cite{DM07}, have used similar separation of variables approaches
but have chosen to approximate the Fourier series of the function
$g$ occuring in the boundary condition, equation~(\ref{bc:4}) by
its first two terms, whereas we treat this more carefully.

Some recent papers treat the problem in the context of slip flows in elliptic microchannels,
e.g.~\cite{DT16,DM07,Duan07}.
(Again in the context of flows in microchannels,
 time-dependent and other aspects are treated in various papers.
\cite{Gupt,Hasl,Suh,Suh18} study pulsatile flows.
Papers~\cite{Day55,Tr73} are relevant to studies of the decay of transients and starting flows.
Other physical effects are studied in~\cite{Ghos,Hsu}.
The steady flow case is, for each of these, a special case when a parameter is set to zero.)
We make considerable effort to check our results.
For example, there are a number of rigorously proved bounds.
An old variational result, a lower bound on $\QQ$ is equation (4.4) of ~\cite{KM93},
\begin{eqnarray}
\QQ(\beta)
&\geq& \QQ(\beta=0) + \frac{\beta\, |\Omega|^2}{|\partial\Omega|}\quad
{\rm \ for\ general\ }\Omega, \nonumber\\
&=& \frac{\pi a^3 b^3}{4 (a^2 + b^2)}  +
 \frac{\pi^2 \beta\, a^2 b^2}{|\partial\Omega|}\quad\quad
{\rm \ for\ ellipse\ }\Omega .
\label{KM93eq4p4}
\end{eqnarray}
Inequality~(\ref{KM93eq4p4}) and the more general form preceding it are very easy to establish: see~\S\ref{sec:Varl}.
Inequality~(\ref{KM93eq4p4}) becomes an equality when the ellipse is a circle, i.e. $b=a$.
The left-hand term in the right-hand side of inequality~(\ref{KM93eq4p4})
gives the dominant term in the asymptotic expansion of $Q$ as $\beta$ tends to 0.
The right-hand term in inequality~(\ref{KM93eq4p4})
gives the dominant term in the asymptotic expansion of $Q$ as $\beta$ tends to infinity.
However (except when $b=a$), in both asymptotic limits, the other term is not the next term in the asymptotic expansion.

We sometimes denote the area $|\Omega|$ by $A$
and the perimeter  $|\partial\Omega|$ by $L$.
For an ellipse, the formula for $L$ is given in equation~(\ref{eq:perimE}).

A much more recent inequality, a generalization of the St Venant Inequality (see~\cite{PoS51}), is
given in~\cite{BG15} (with related work in~\cite{BB13}):
\begin{equation}
\QQ
\leq
Q_\odot(\beta)\,
:= \frac{\pi (a\, b)^{3/2}}{8}\left(\sqrt{a b} + 4\beta\right) .
\label{eq:isoperQ}
\end{equation}
In words: amongst all ellipses with a given area that which has the greatest $\QQ$ is the circular disk.
See equation~(\ref{eq:circleQinf}).
(This, and other isoperimetric inequalities, are  reviewed,
in the context of microchannels in~\cite{KW16,KWK18}.)

\subsection{The structure of this Part}\label{subsec:structure}

The rest of this part is organized as follows.
\begin{itemize}
\item In \S\ref{sec:explicit} we summarize the very well-known elementary solutions, polynomial in Cartesian coordinates, for the simplest special cases against which to  check our general solution.
Each of them forms the lowest order term in asymptotics developed later in this 
part. 
\item
In \S\ref{sec:GovElliptic}, the boundary value problem 
is formulated  in elliptic cylindrical coordinates, and then
the series representation is approximated
(to any accurcy desired).

\item The next sections provide many checks.
To reduce the number of parameters, in these sections, we take $b=1/a$, so the area of the ellipse is $\pi$.
\begin{itemize}
\item In \S\ref{sec:Varl} we present a variational formulation and associated bounds.
\item In \S\ref{sec:nearCircQs} we present asymptotics for the situation where the ellipse is nearly circular.
\item In \S\ref{sec:betaSmall} we give asymptotics for $\beta$ small.
\item In \S\ref{sec:betaLarge} we study asymptotics for $\beta$ large.\\
\end{itemize}
\item Finally, a conclusion is given in \S\ref{sec:Conclusion}.
\item \S\ref{sec:GeometryI}, Appendix A, treats various geometric matters.\\
\S\ref{appx:EllipticInts}, Appendix B, has some notes on elliptic integrals.
These arise in the Fourier series for a functions $g$ and $\hat g$
used at several sections of ths work, notably 
\S\ref{sec:GovElliptic} and 
developments from \S\ref{sec:betaLarge} given in Part II.\\
\S\ref{sec:Numerical}, Appendix C, gives some numerical examples, its first parts comparing our numerics with some previously published examples.
We conclude this appendix with an attempt to suggest numbers which might arise in connection with possible applications
involving blood flow.
\end{itemize}

\section{Some simple explicit solutions:\\ $u$ quadratic polynomial in $x$, $y$}\label{sec:explicit}

\subsection{Circular cross-section with $\beta\ge{0}$}\label{subsec:Circular}

When $\Omega$ is circular, radius $a$, in polar coordinates, $r=\sqrt{x^2+y^2}$, the solution is
\begin{equation}
u_{\beta\odot}=  \usteady 
= \frac{1}{4}\left(a^2-r^2\right) +\frac{\slipparameter a}{2} .
\label{eq:circuSteady}
\end{equation}
In the context of fluid flows,
when $\slipparameter=0$ this is Poiseuille flow.

\begin{equation}
Q_{\beta\odot}=\QQ(\slipparameter) 
=  \frac{\pi a^3}{8} (a+4\slipparameter)  . 
\label{eq:circleQinf}
\end{equation}
See also~\cite{La32}\S{331},~p586.

In later sections
the circle will have radius $a=1$, area $\pi$.

\subsection{Elliptic cross-section with $\beta={0}$}\label{subsec:beta0}

In this subsection (and in sections \S\ref{sec:Varl} to  \S\ref{sec:betaLarge}) we consider ellipses with $b=1/a$ and $a\ge{1}$.
The $\beta=0$ solution dates back at least to St Venant.
In polar coordinates with the origin at the centroid of the ellipse the velocity $u_e$
is, in Cartesian coordinates, a quadratic polynomial in $x$ and $y$ while,
in polar coordinates, it is
$$ u_{0e}
=\frac{1}{4}\left(   \sqrt {1-\epsilon^{2}} -r^2 +\epsilon \, r^2 \, \cos(2\theta)\, \right) .
$$
Here
$$\epsilon = \frac{a^2 - a^{-2}}{a^2 +a^{-2}} . $$
In Cartesian coordinates
\begin{equation}
u_0
=\frac{1- (x/a)^2 - (y/b)^2}{2/a^2 + 2/b^2} .
\label{eq:u0Cart}
\end{equation}
Also
\begin{equation}
Q_{0e}=\QQ({\rm ellipse}, \beta=0)
=\frac{\pi}{4 (a^2 + a^{-2})}
= \frac{\pi}{8} \, \sqrt {1-\epsilon^{2}} .
\label{eq:Qellb0}
\end{equation}

Also relevant is $e$  the eccentricity as defined in equation~(\ref{eq:eDef}).
\begin{figure}[hb]
\centerline{\includegraphics[height=7cm,width=13cm]{\grpath/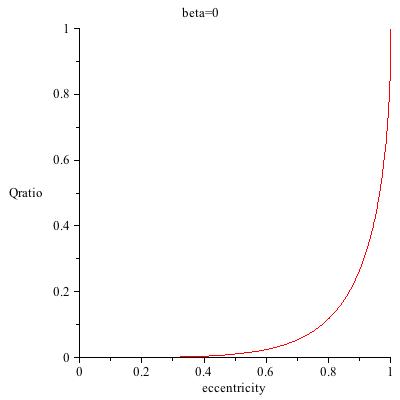}}
\caption{$\beta=0$: plot of Qratio=$(Q_\odot-Q(e))/Q_\odot$  against  eccentricity, $e$.}
\label{fig:QratioBeta0jpg}
\end{figure}

There are various interesting alternative expressions for $\QQ$.
The ellipse's moments of inertia are
$$ I_{xx}=\frac{\pi a^2}{4} , \quad I_{yy}=\frac{\pi}{4 a^2}\
{\rm so\ \ } \QQ =  J:= \frac{I_{xx} I_{yy}}{I_{xx} +I_{yy}} . $$
See~\cite{PoS51} p112.
For any domain (and $\beta=0$) it is known that $\QQ\le{J}$.

\clearpage

\section{Fourier series solution}\label{sec:GovElliptic}

\subsection{Using elliptic coordinates}\label{subsec:EllCoords}

\subsubsection*{Formulae when $a\ge{b}>0$}
\noindent For the geometry of elliptic domains specified as in equation~(\ref{eq:ellab}) with, as always, $a>b$,
it is convenient here  to use elliptic coordinates $(\eta, \psi)$ which are related to the rectangular coordinates
(as in~\cite{DT16} equation~(18)) by
 \begin{equation}\label{eq:5}
 x = c \cosh(\eta) \cos(\psi);  \;\;\;\;\;y = c \sinh(\eta) \sin(\psi); ,
 \end{equation}
\noindent where $0 \leq \eta \leq \infty,\;0 \leq \psi \leq 2 \pi$,
 and $c$ and $-c$ are two common foci of the ellipse.
The origin of coordinates in the Cartesian system 
(the centre of our ellipse) is at $\eta=0$, $\psi=\pi/2$.
The upper side of the positive semi-major axis is $\eta=0$,
$0<\psi<\pi/2$:
the lower side is $\eta=0$, $0>\psi>-\pi/2$.
Let $a$ and $b$ denote lengths of the semi-axes of the ellipse and $a > b > 0$.
Except in this subsection, and elsewhere when so specified, we take $b=1/a$.
Returning to the general situation,
$$a = c \cosh(\eta) {\rm\ \ and\ \ } b = c \sinh(\eta),
{\rm \ \ and\ \ } c=\sqrt{a^2-b^2}.$$

\noindent Defining $\eta_0$ (as in~\cite{DT16} equation~(19)) by
$$\eta_0 = \ln\frac{1+b/a}{\sqrt{1-(b/a)^2}}
= {\rm arctanh}(\frac{b}{a}),$$
the boundary of the ellipse can be represented by
$$ x= c \cosh(\eta_0) \cos(\psi), \qquad y= \sinh(\eta_0) \sin(\psi),$$
so that
$$c=\frac{a}{\cosh \eta_0}=\frac{b}{\sinh \eta_0}.$$
The eccentricity is
$$e = \sqrt{1- \frac{b^2}{a^2}} = \sqrt{1-\tanh^2(\eta_0)}
=\frac{1}{\cosh(\eta_0)}.$$
(Nearly circular ellipses have $\eta_0$ large and $c$ near 0.)

The Jacobian $J$ of the transformation from $(x,y)$-coordinates to $(\eta,\psi)$ coordinates is
$$J(\eta,\psi)= c^2  \left( \cosh^2(\eta) -\cos^2(\psi) \right)
= \frac{c^2}{2}  \left( \cosh(2\eta) -\cos(2\psi) \right) . $$
This is needed in the derivation of equation~(\ref{eq:QJac}).
The Jacobian will be used in other places in this work, e.g.
\S\ref{subsec:Qinf}.

\subsubsection*{Formulae when $b=1/a$}

In all major calculations henceforth  we scale distances so that we have $b=1/a$.
Formulae used later include the following:
\begin{eqnarray}
c^2 =\frac{2}{\sinh(2\eta_0)}
&=& a^2-\frac{1}{a^2}=\frac{e^2}{\sqrt{1-e^2}},
\label{eq:cae}\\
\tanh(\eta_0)
&=& \frac{1}{a^2} ,
\label{eq:tanha}\\
\cosh(\eta_0)
&=& \frac{1}{e} ,
\label{eq:coshe}\\
e^2
&=& 1- \frac{1}{a^4}.
\label{eq:eaDef}
\end{eqnarray}

\subsection{The pde and boundary conditions}

The pde   of Problem (P($\beta$)) 
is recast in these coordinates as
\begin{eqnarray}\label{eq:7}
 \frac{\partial^2 u}{\partial \eta^2} + \frac{\partial^2 u}{\partial \psi^2}
&=&-
c^2 (\cosh^2(\eta)-\cos^2(\psi) ) ,\nonumber \\
&=& -
\frac{c^2}{2}(\cosh(2\eta)-\cos(2\psi)) .
\end{eqnarray}
We will solve the pde (\ref{eq:7}) subject to boundary conditions (\ref{bc:1}), (\ref{bc:2}) and (\ref{bc:3}).
The solution is assumed to be symmetric about the mid-plane for $\psi\;=\;0,\;\frac{\pi}{2}$.
Therefore, symmetric boundary condition is applied, i.e.,
\begin{equation}\label{bc:1}
\frac{\partial u}{\partial \psi}\; = 0\; \mathrm{on}\;\;\psi\;=\;0,\;\frac{\pi}{2}.
\end{equation}
In addition, we set
\begin{equation}\label{bc:2}
 \frac{\partial u}{\partial \eta} = 0\;\;\mathrm{on}\;\;\eta=0.
\end{equation}

The boundary condition on the ellipse,
i.e. at those points when $\eta=\eta_0$ is
\begin{equation}\label{bc:3}
u(\eta_0,\psi)+\frac{\beta}{c\sqrt{\cosh^2\eta_0 - \cos^2 \psi}} \frac{\partial u}{\partial \eta}(\eta_0,\psi)  \;=\;0.
\end{equation}
Define
\begin{equation}\label{eq:gdef}
 g(\psi)= \frac{\cosh(\eta_0)}{ (\cosh^2(\eta_0) - \cos^2(\psi))^{1/2}} .
\end{equation}
Boundary condition (\ref{bc:3}) becomes
\begin{equation}\label{bc:4}
u(\eta_0,\psi)+\frac{\beta}{c \cosh(\eta_0)}g(\psi) \frac{\partial u}{\partial \eta}(\eta_0,\psi)  \;=\;0.
\end{equation}

\noindent In summary,
the problem in elliptic coordinates is given as
the following boundary value problem (BVP).\\

\noindent {\bf BVP}: Find  $u(\eta,\psi)$  such that the above governing equation (\ref{eq:7}) and associated boundary conditions
(\ref{bc:1}), (\ref{bc:2}) and (\ref{bc:4})  are satisfied.

\medskip

\subsection{The Fourier series for $g$}\label{subsec:gFS}

\subsubsection*{Ellipses in general}

The definition of $g$ from~(\ref{eq:gdef}) can be re-cast in several ways:
\begin{eqnarray*}
 g(\psi)
 &=& \left(1-\frac{\cos^2(\psi)}{\cosh^2(\eta_0)}\right)^{-1/2}  ,\\
 &=& C\, g_*(\psi)  \ {\rm with \ \ }
 g_*(\psi) = \left(1-\frac{\cos(2\psi)}{\cosh(2\eta_0)}\right)^{-1/2}\ {\rm and}\ \
 C= \left(\frac{1+\cosh(2\eta_0)}{\cosh(2\eta_0)} \right)^{1/2}.
\end{eqnarray*}

We will, later, use a more economic notation using
 \begin{equation}
 q=\cosh^2(\eta_0) ,
 \label{eq:qEta0}
\end{equation}
and the form for $g_*$
is consistent with the occurence, in Fourier coefficients, of polynomials in $\cosh(2\eta_0)=2 q-1$.

The Fourier series for $g$ is used. Define
\begin{equation}
g_n = \frac{2}{\pi}\int_0^\pi g(\psi)\cos(2 n\psi)\, d\psi ,
\label{eq:bgnDef}
\end{equation}
so
$$ g(\psi) = \frac{g_0}{2} + \sum_{n=1}^\infty g_n\, \cos(2 n\psi) ,$$
In various approximations  the first few terms in the Fourier series are used.
For these define
\begin{equation}
{\rm EllipticE}_0= {\rm EllipticE}(\frac{1}{\cosh(\eta_0)}) \qquad{\rm and\ \  }
{\rm EllipticK}_0= {\rm EllipticK}(\frac{1}{\cosh(\eta_0)}).
\label{eq:Ellnotn}
\end{equation}
We have
\begin{eqnarray}
g_0
 &=& \frac{4}{\pi}\, {\rm EllipticK}_0 ,
 \label{eq:g0Ell}\\
g_1
  &=& -\frac{8\cosh(\eta_0)^2}{\pi}\, {\rm EllipticE}_0  +
   \frac{4(2\cosh(\eta_0)^2-1)}{\pi}\, {\rm EllipticK}_0 \ .
   \label{eq:g1Ell}
\end{eqnarray}

All the terms $g_n$ can be written in terms of Legendre $Q$ functions though we first encountered them in the form presented by maple,
\begin{equation}
g_n =\frac{4}{\pi}\left(  E_n  \, {\rm EllipticE}_0 +
K_n\, {\rm EllipticK}_0 \right), 
\label{eq:gnEK}
\end{equation}
where $E_n$ and $K_n$ are polynomials of degree $n$ in $\cosh(\eta_0)^2$
with rational number coefficients.

Should future work require more terms we remark that a three term recurrence relation exists to
 determine the polynomials $E_n$ and $K_n$.
See Part II \S~\ref{app:Furtherg}.
 In the notation of equation~(\ref{eq:qEta0}), the $K_n$ sequence of polynomials starts with
 $$ K_0=1 , \qquad K_1= 2q-1 = q_2 .
 $$
 The $E_n$ sequence of polynomials starts with
 $$ E_0=0 , \qquad E_1= 2q = q_2-1 .
$$

There are also other representations: see  Part II \S~\ref{app:Furtherg}.
 There are several very elementary facts which have some use.
 These include that
 $$ \cos(2 k\psi)= {\rm ChebyshevT}(k,\cos(2\psi)) .$$
 Also, the Taylor series
 $$ (1-X)^{-1/2}= \sum_{k=0}^\infty
 \frac{(2 k)!}{2^{2k}( k!)^2} X^k ,$$
 can be applied to the preceding expressions for $g$
 to obtain series which are sums of powers of $\cos(\psi)^2$ or of $\cos(2\psi)$.
To date, our only use of these has been in checking the asymptotic approximations, at small eccentricity, to the Fourier coefficients of
$g$ and of $1/g$.

The reciprocal of $g$ also occurs in later calculations
(and also in~\cite{DT16} equation~(61) though they use a Taylor series rather than Fourier series used here).

\cmdvmcode{ellipseFSmpl.txt}{ }

\subsubsection*{Approximating the Fourier series for nearly circular ellipses}

Nearly circular ellipses have $\eta_0$ large.
Then  the lowest approximation is
$$ g(\psi)\sim g_\epsilon(\psi)
= 1+\frac{\cos^2(\psi)}{2\cosh^2(\eta_0)}
=1+\frac{1}{4\cosh^2(\eta_0)}+\frac{\cos(2\psi)}{4\cosh^2(\eta_0)} ,
$$
and
\begin{equation}
\frac{1}{g(\psi)}\sim
1-\frac{1}{4\cosh^2(\eta_0)}-\frac{\cos(2\psi)}{4\cosh^2(\eta_0)} .
\label{eq:gRecipFSnearCirc}
\end{equation}
We will use the next approximation after this, and will use $e=1/\cosh(\eta_0)$:
\begin{equation}
\frac{1}{g(\psi)}\sim
1-\frac{1}{4}e^2 - \frac{3}{64} e^4 -
 \left( \frac{1}{4}e^2+ \frac{1}{16}e^4 \right) \cos(2\psi) -
 \frac{1}{64}e^4 \, \cos(4\psi) .
\label{eq:gRecipFSnearCirc2}
\end{equation}
\ifthenelse{\boolean{vmcode}}{The code deriving this is:}{}
\cmdvmcode{Recurrence/recipgmpl.txt}{ }

We have, for $e$ small the following Fourier coefficients of
${\hat g}=1/g$:
\begin{eqnarray}
{\hat g}_0
&\sim& 2-\frac{1}{2} e^2-\frac{3}{32} e^4
\label{eq:B0e} \\
{\hat g}_1
&\sim& -\frac{1}{4} e^2 -\frac{1}{16} e^4
\label{eq:B1e} \\
{\hat g}_2
&\sim& -\frac{1}{64} e^4
\label{eq:B2e}
\end{eqnarray}
(We will need these in the treatment of $u_\infty$. See equation~(\ref{eq:rgFSB}).)

\medskip

\subsection{Completing the Fourier series solution}\label{subsec:Analytical}

The general solution of equation (\ref{eq:7}) can be expressed in the form of
\begin{equation}\label{us}
  u(\eta, \psi)\;=\;v_h+v_p,
\end{equation}
where $v_h$ and $v_p$ denote the homogeneous and the particular solution of the pde
   (\ref{eq:7}) subjecting to the boundary conditions (\ref{bc:1}),\~(\ref{bc:2}) and(\ref{bc:4}).\\

For the particular solution $v_p$,  from the right hand side of the equation (\ref{eq:7}), we can expressed $v_p$ in the form of
\begin{equation}\label{eq:9b}
v_p\;=\;{\hat g}_1\cosh(2\eta)+{\hat g}_2\sinh(2\eta)+{\hat g}_3\cos(2\psi)+{\hat g}_4\sin(2\psi),
\end{equation}
which gives
\begin{equation}\label{eq:10b}
\frac{\partial^2 v_p}{\partial \eta^2}\;=\;2 {\hat g}_1\sinh(2\eta)+2 {\hat g}_2\cosh(2\eta);
\end{equation}
\begin{equation}\label{eq:11b}
\frac{\partial^2 v_p}{\partial \psi^2}\;=-2 {\hat g}_3\sin(2\psi)+2 {\hat g}_4\cos(2\psi).
\end{equation}
Substituting
$v_p$ 
into the pde (\ref{eq:7}) and comparing the coefficients of $\cosh(2\eta)$ and $\cos(2\psi)$ yields
\begin{equation}\label{eq:12b}
v_p\;=\;-\frac{c^2}{8} (\cosh(2\eta)+\cos(2\psi))
\ \ \left( = -\frac{1}{4}\, (x^2+y^2)\ \right).
\end{equation}

\noindent By separation of variables, the solution $v_h$ of the Laplace equation  subject to above boundary conditions (\ref{bc:1}) and (\ref{bc:2}) is
\begin{equation}\label{eq:8b}
v_h\;=\;\sum_{n=0}^{\infty}A_n \, \cosh(2n\eta)\cos(2n\psi).
\end{equation}
Hence, the general solution can be expressed in the form of
\begin{equation}\label{eq:13}
u(\eta,\psi)\;=\;v_h+v_p\;=\;\sum_{n=0}^\infty A_n \, \cosh(2n\eta)\cos(2n\psi) - \frac{c^2}{8}(\cosh(2\eta)+\cos(2\psi)),
\end{equation}
where term $A_n$ can be determined by using the boundary condition (\ref{bc:4}).
We then obtain
\begin{eqnarray}\label{eq:14}
&\sum_{n=0}^\infty A_n \,  \cosh(2 n\eta_0)\cos(2 n\psi)-
\frac{c^2}{8}\, (\cosh(2\eta_0)+\cos(2\psi)) + \frac{\beta g(\psi)}{c\cosh(\eta_0)}  \nonumber \\
&\sum_{n=0}^\infty n 2 A_n \sinh(2 n \eta_0)\cos(2 n \psi)-\frac{c^2 \beta g(\psi)}{4 c \cosh(\eta_0)}\sinh(2 \eta_0) = 0.
\end{eqnarray}
or
\begin{eqnarray}\label{eq:18}
&\sum_{n=0}^\infty A_n \, \left( \cosh(2 n\eta_0) + \frac{2 n \beta g(\psi)}{c\cosh(\eta_0)} \sinh(2 n \eta_0) \right)\cos(2 n\psi) \nonumber\\
&=\frac{c^2}{8}(\cosh(2\eta_0)+\cos(2\psi)) + \frac{c^2 \beta g(\psi)}{4 c \cosh(\eta_0)}\sinh(2 \eta_0)\;=\;f(\psi).
\end{eqnarray}
As $g$ is an even function so is $f$ and $h_n$ where
$$h_n(\psi)=\cosh(2 n\eta_0) + \frac{2 n \beta g(\psi)}{c\cosh(\eta_0)} \sinh(2 n \eta_0). $$
(The Fourier coefficients of $f$ and of $h_n$
are very simply related to the Fourier coefficients of $g$.)
We then have
%
%
\begin{equation}
\sum_{n=1}^\infty A_n \int_{-\pi}^{\pi} h_n(\psi) \cos(2n\psi)\cos(2m\psi) d\psi = - \int_{-\pi}^{\pi} f(\psi) \cos(2m \psi) d\psi.
\end{equation}
Hence if we truncate the summation at $n_{\rm max}$ we have a system of linear equations for the $A_n$ for $n$ up to $n_{\rm max}$.
Using the $A_n$ so found, we obtain an approximation to the Laplace solution $v_h$ and then steady-state fluid flow, $u(\eta,\psi)=v_h+v_p$ in the elliptical channel.  Finally, the volume flow rate of fluid passing through a elliptical cross sectional area, $Q$, can determined by
\begin{equation}
  \QQ\;=\; {c^2} \int_0^{2\pi} \int_0^{\eta_0} u(\eta,s) \left( \cosh^2(\eta) -\cos^2(s) \right) d \eta ds.
\label{eq:QJac}
\end{equation}

\subsection{$\QQ$ from the Fourier series of $v_h$}

We will substiute $u=v_p+v_h$ with $v_h$ given by its Fourier series into
the formula for $Q$ given in~(\ref{eq:QJac}).
Note that all the Fourier components with $n>1$ contribute nothing to the integral.
First
\begin{eqnarray*} Q_p
&=& c^2 \int_0^{2\pi} \int_0^{\eta_0}
v_p \,  \left( \cosh^2(\eta) -\cos^2(s) \right) d \eta ds , \\
&=& -\frac{1}{4} I_2
= -\frac{1}{4}\, \frac{\pi}{4} \, (a^2+\frac{1}{a^2})
= - \frac{\pi}{16} \, \frac{2-e^2}{\sqrt{1-e^2}} ,
\end{eqnarray*}
$I_2$ being the polar area moment of inertia about the centroid.
On using $c^2=2/\sinh(2\eta_0)$,
\begin{eqnarray}
 Q
 &=& Q_p+ c^2 \int_0^{2\pi} \int_0^{\eta_0}
( A_0 + A_1  \cosh(2\eta)\cos(2 s)) \left( \cosh^2(\eta) -\cos^2(s) \right) d \eta ds ,
 \nonumber \\
&=& - \frac{\pi}{8} \, \frac{\cosh(2\eta_0)}{\sinh(2\eta_0)}+ \pi\, A_0 - \frac{1}{2}\,\pi\, A_1 .
\label{eq:QFS}
\end{eqnarray}

\begin{verbatim}
uA := A0+A1*cosh(2*eta)*cos(2*psi)+A2*cosh(4*eta)*cos(4*psi);
   `cSqr Assign;1/`(sqrt(1/e^2-1)/e):
 hypToe := ->  subs(cosh(eta0) = 1/e, sinh(eta0) = sqrt(1/e), u) :
tmp := subs(cSqr = 1/(cosh(eta0)*sinh(eta0)), int(int(uA*cSqr*(cosh(eta)^2-cos(psi)^2), eta = 0 .. eta0), psi = -Pi .. Pi))
     A0*Pi -A1*Pi/2
\end{verbatim}
The terms involving the $A_j$ are checked in the following code:
\cmdvmcode{QfromFS0mpl.txt}{ }

\ifthenelse{\boolean{vmcode}}{Rather more, with other checks is given in the following:}{}
\cmdvmcode{QfromFSmpl.txt}{ }

\clearpage

\section{Inequalities from a variational formulation}\label{sec:Varl}

Variational formulations are available. For example, in~\cite{KM93} equation (4.1), the functional $\cal J$ is defined as
\begin{equation}
{\cal J}(v) ={\cal A}(v) - \frac{1}{\beta}{\cal B}(v) \quad{\rm where\ \ }
{\cal A}(v)= \int_\Omega(2 v -|\nabla v|^2),  \ \ {\cal B}(v) =\int_{\partial\Omega} v^2 .
\label{eq:Varl}
\end{equation}
For any simply-connected domain with sufficiently smooth boundary,
the maximiser $u$ of this functional 
solves Problem (P($\beta$)).
Furthermore, ${\cal J}(u)$ is $\QQ={\cal Q}(u)$ where
$${\cal Q}(v) = \int_\Omega v . $$

Thus we can find lower bounds for $\QQ$ by evaluating ${\cal J}(v)$ for particular choices of $v$.

In this work, we have two distinct variational approximations.
In the order of importance, though not of our discoveries, these are
as follows.\\
$\bullet$
The more important of these -- applying to general domains,
not just ellipses and now published in~\cite{KW20} --
is included here, in \S\ref{subsec:inS0Sinf1},
 as motivation for further work
on Problem (P($\infty$)) for the ellipse.
We have, in Part II, a series solution for $\Sigma_\infty$ for the ellipse
involving Legendre functions but have yet to find $\Sigma_1$.\\
The lower bound is more easily illustrated with the rectangle
(see~\cite{KW20}) here treated in Part III. 
For the rectangle (and some other domains, see~\cite{Ke20i,Ke20},
$u_\infty$ is available analytically.\\
$\bullet$
The second of these is ad hoc,  good when the solutions of Problem (P($\beta$))
have nearly elliptical level curves.
This is the case for the ellipse, and also, 
{\em though only for large $\beta$}
for a number of other domains, including rectangles.
General domains are treated in~\S\ref{subsec:VarlGen}
and specialised to ellipses in~\S\ref{subsec:VarlEll}.

\subsection{A development from \cite{KM93}}\label{subsec:inS0Sinf1}

Our notation is as in~\cite{KM93}.
For large $\beta$
$$ u \sim \beta\, \frac{|\Omega|}{|\partial\Omega|}  +  u_\infty + o(1) , $$
where $u_\infty$ solves
$$ - \frac{\partial^2 u_\infty}{\partial x^2}- \frac{\partial^2 u_\infty}{\partial y^2}
 = 1\, {\rm in}\,\Omega ,\,
\frac{\partial u_\infty}{\partial n}
= -\frac{|\Omega|}{|\partial\Omega|} \, {\rm on}\,\partial\Omega \,
{\rm and}\,  \int_{\partial\Omega} u_\infty = 0 .
\eqno({\rm P}(\infty))
$$
Define
\begin{equation}
\Sigma_\infty
= \int_\Omega u_\infty , \qquad  {\rm and\ \ } \Sigma_1= - \int_{\partial\Omega} u_\infty^2 ,
\label{eq:SigDef1}
\end{equation}
with $u_\infty$ satisfying Problem (P($\infty$)).
A lower bound on $\QQ(\slipparameter)$ using geometric quantities --
area and perimeter -- and $\QQ_0$ and the two $\Sigma$ functionals is 
presented in Part III~\S\ref{subsec:inS0Sinf}.

\subsection{Quadratic test functions, and a simple first case}\label{subsec:VarlGen}

For domains like ellipses symmetric about both the $x$-axis and $y$-axis
one nice choice for $v$ are the quadratic functions
$$\leqno{(i)}\qquad\qquad\qquad\qquad\qquad   v = c_0+c_2\left( (\frac{x}{a})^2 +a^2 y^2\right) , $$
and, more generally (and with a different $c_0$)
$$\leqno{(ii)}\qquad\qquad\qquad\qquad\qquad v = c_0+c_{xx} x^2 + c_{yy} y^2 . $$
(Case (ii) is motivated by the likelihood that for a range of $\beta$ and a range of $a$
the level curves of $u$ look to be like ellipses. See, for example, Figure 4 of~\cite{SV12}
and our own contour plots shown in Appendix C \S~\ref{sec:Numerical}.)

We briefly return to general domains.
An application of the Divergence Theorem give
\begin{equation}
{\rm if\ \ } -\Delta v=1\ \ {\rm then\ \  }
{\cal J}(v)-{\cal Q}(v) = \int_{\partial\Omega} v (v+\beta \frac{\partial v}{\partial n} )
\label{eq:JmQDgce}
\end{equation}
\medskip

\noindent{\bf Case (i)}
\smallskip

We begin with case (i).
While the beginning of the study with case (ii) will allow $\Omega$ to be more general than merely an ellipse,
for case (i) we restrict to the ellipse.
 The final result, lower bound, will give the same result as we reported in \S\ref{subsec:GovEq}
from equation (4.4) of~\cite{KM93}.
Write $J$ or $J(c_0,c_2)$ for the value of $\cal J$ evaluated at the quadratic function (i)
$$ J = 2 c_0 A+ 2 c_2\left(  \frac{1}{a^2} (1-\frac{2 c_2}{a^2}) I_{XX} +a^2 (1- 2 a^2 c_2)I_{YY}\right) -(c0+c2)^2 \frac{L}{\beta} .
$$
Finding the gradient of $J$ with respect to $[c_0,c_2]$ we find the maximum of $J$ will occur when
the gradient is 0, i.e.
\begin{eqnarray*}
L c_0 + L c_2
&=& \pi \beta ,\\
2L c_0 + c_2\left( 2L +\pi\beta(a^2+\frac{1}{a^2})\right)
&=& \pi\beta .
\end{eqnarray*}
This is readily solved, and on integrating, yields the result given before at
inequality~(\ref{KM93eq4p4}).

\cmdvmcode{VarlQuadr/QKM93derivempl.txt}{ }

\medskip
\cmdvmcode{VarlQuadr/QKM93fnmpl.txt}{ }

It is straightforward to show that the quadratic $u$ which is the variational winner in Case (i)
is such that the result in~(\ref{eq:JmQDgce}) can be applied to establish
${\cal J}(u)={\cal Q}(u)$
(a fact which can be established by other means).
The application of~(\ref{eq:JmQDgce})  is facilitated by the fact that
the quadratic $u=u_0 +\beta|\Omega|/{|\partial\Omega|}$ is constant on $\partial\Omega$.

\medskip

\noindent{\bf Case (ii), an introduction}
\smallskip

Next we begin to treat case (ii).
Write $J$ or $J(c_0,c_{xx},c_{yy})$ for the value of $\cal J$ evaluated at the quadratic function (ii).
\begin{eqnarray*}
J
&=& 2\left( c_0 A + c_{xx}(1-2 c_{xx}) I_{xx} + c_{yy}  (1-2 c_{yy}) I_{yy} \right)-\\
&\ &\frac{1}{\beta}\left(
c_0^2 L + 2 c_0 c_{xx} i_{xx}+ 2 c_0 c_{yy} i_{yy} + 2 c_{xx} c_{yy} i_{xxyy} + c_{xx}^2 i_{xxxx} + c_{yy}^2 i_{yyyy}
\right) \ .
\end{eqnarray*}
where $A$ is the area, $L$ the perimeter and the $I$ are area moments and $i$ boundary moments
$$ I_{xx}=\int_\Omega x^2, \ \  I_{yy}=\int_\Omega y^2 , $$
and
$$
i_{xx}= \int_{\partial\Omega} x^2 , \ \ i_{yy}= \int_{\partial\Omega} y^2 ,$$
and
\begin{equation}
i_{xxyy}= \int_{\partial\Omega} x^2 y^2 , \ \
i_{xxxx}= \int_{\partial\Omega} x^4 , \ \ i_{yyyy}= \int_{\partial\Omega} y^4 .
\label{eq:idef2}
\end{equation}
Define also
\begin{equation}
 Q_q = c_0 A + c_{xx}I_{xx} + c_{yy}I_{yy} ,
\label{eq:Qq}
\end{equation}
the subscript $q$ reminding us of the quadratic approximation to the velocity field.
The gradient of $J$, here denoted $\mathbf g$, is
$$
 {\mathbf g}
 =  -\frac{2}{\beta}  \left(
 \begin{array}{c}
 c_0 L + c_{xx} i_{xx} + c_{yy} i_{yy}
-  \beta A\\
c_0 i_{xx} +c_{xx} (4\beta I_{xx} + i_{xxxx}) + c_{yy} i_{xxyy}
-  \beta I_{xx}\\
c_0 i_{yy} +c_{xx} i_{xxyy} +c_{yy} (4\beta I_{yy} + i_{yyyy})
-  \beta I_{yy}
\end{array}
\right)\ .
$$
We remark that
\begin{equation}
 2(J -Q_q)
= {\mathbf g}[1] c_0 + {\mathbf g}[2] c_{xx} + {\mathbf g}[3] c_{yy} .
\label{eq:JminusQ}
\end{equation}

\ifthenelse{\boolean{vmcode}}{Code for setting up the equations is as follows:}{ }
\cmdvmcode{VarlQuadr/varlGenStartmpl.txt}{ }

This gradient $\mathbf g$ is zero when
\begin{eqnarray*}
c_0 L + c_{xx} i_{xx} + c_{yy} i_{yy}
&=&  \beta A , \\
c_0 i_{xx} +c_{xx} (4\beta I_{xx} + i_{xxxx}) + c_{yy} i_{xxyy}
&=&  \beta I_{xx} , \\
c_0 i_{yy} +c_{xx} i_{xxyy} +c_{yy} (4\beta I_{yy} + i_{yyyy})
&=&  \beta I_{yy} .
\end{eqnarray*}

Up until now the $c_0$, $c_{xx}$, $c_{yy}$ have been general.
Henceforth we use the same letters to denote the $c$ obtained by solving ${\mathbf g}=0$, the
gradient of $J$ equals zero.
We remark that after finding $c_0$, $c_{xx}$, $c_{yy}$ there are two methods of estimating $\QQ$.
On the one hand, we can estimate $\QQ$ from
\begin{equation}
Q_c := c_0 A + c_{xx} I_{xx} + c_{yy} I_{yy} ,
\label{eq:Qc}
\end{equation}
while, on the other hand,  we have the lower bound
$$J_c := J(c_0,c_{xx},c_{yy})\leq{\QQ} . $$
\ifthenelse{\boolean{vmcode}}{We noticed that when maple performed the substitutions it found that $Q_c=J_c$.
This led us to find the identity~(\ref{eq:JminusQ}) confirming this.
We have, as yet, no direct proof via a modest amount of hand calculation,
via appropriate theorems, along the lines of using the Divergence Theorem as in
(\ref{eq:JmQDgce}) used at the end of our treatment of case (i).}{
From identity~(\ref{eq:JminusQ}) we have $Q_c=J_c$.}

For general $\Omega$ inequalities on moments
(such as some in references in~\cite{Ke06})
may lead to lower bounds for $Q$ in terms of  simpler geometric
quantities.
Also for some domains, polygons for example, the various moments
could be calculated.
In general the restriction of the test functions to quadratics
isn't likely to give very tight lower bounds though there 
may be exceptions.
Rectangles at large $\beta$ might be one such as
$u_\infty$ (see~\cite{KM93} or equations~(\ref{eq:uInf}))
is, for rectangles, quadratic.

For an ellipse with $a>1$ (and $b=1/a$) the area moments are elementary,
but the perimeter and other boundary moments involve elliptic integrals:
we present these in the next subsection.
For a circle, all the moments  are elementary.
For a circle of radius 1, the moments are:
\begin{eqnarray*}
A
&=&\pi,\ \ I_{xx}=\frac{\pi}{4}, \ \   I_{yy}=\frac{\pi}{4} , \\
L
&=&2\pi, \ \
i_{xx}=\pi=i_{yy}, \ \ i_{xxxx}=\frac{3\pi}{4}=i_{yyyy}, \ \ i_{xxyy}=\frac{\pi}{4} .
\end{eqnarray*}
For the circle, the variational winner has $c_0$, $c_{xx}=c_{yy}$ giving  the exact solution
previously presented in equation~(\ref{eq:circuSteady}):
both $Q_c$ and $J_c$ evaluate to the formula for $\QQ$ given in
equation~(\ref{eq:circleQinf}).

\medskip

A note of caution concerning the behaviour when $\beta=0$ is appropriate here.
If one naively sets $\beta=0$ in the system of equations ${\rm grad}(J)=0$ the right-hand sides are all zero.
The solution with all the $c$ equal to zero incontravertably gives a lower bound on $\QQ$:
the lower bound being zero.
Of course if $\beta=0$ the term ${\cal B}(v)/\beta$ is the difficulty for functions $v$ not vanishing identically over the whole boundary.
As a quadratic function vanishing on the boundary provides the minimizer for ellipses with $\beta=0$
it seems appropriate, at this stage, to return to the specifics for an ellipse.

\subsection{Variational methods for the ellipse, Case (ii) continued}\label{subsec:VarlEll}

For the ellipse we can write $i_{xx}$ and $i_{yy}$ in terms of the polar second moment about the centre,
$i_2=i_{xx}+i_{yy}$, and the perimeter $L$
Similarly the 4-th moments can all be written in terms of the  polar fourth moment about the centre, $i_4$
and $i_2$ and $L$.
\cmdvmcode{VarlQuadr/checksimpl.txt}{ }
Each of $L$, $i_2$ and $i_4$ can be integrated in terms of elliptic integrals.
\ifthenelse{\boolean{vmcode}}{Maple code for these is as follows:}{ }
\cmdvmcode{VarlQuadr/LiDerivempl.txt}{ }

\ifthenelse{\boolean{vmcode}}{
For evaluating numerical values we use:}{ }

\cmdvmcode{VarlQuadr/Li2i4fnmpl.txt}{ }
\ifthenelse{\boolean{vmcode}}{
These are used in the formulae for the $c$ and hence the $J(v_*)$.
The solution of the system for the gradient being zero is, in code:}{ }

\cmdvmcode{VarlQuadr/varlGenSolvempl.txt}{ }
\ifthenelse{\boolean{vmcode}}{
We cut and paste the solutions for $c_0$, $c_{xx}$ and $c_{yy}$ and use them in the formula~(\ref{eq:Qc})  for $Q_c$.}{ }

\cmdvmcode{VarlQuadr/QcFunmpl.txt}{ }
\medskip


The variational winner over the quadratics gives lower bounds on $Q$
as in the following code.
\newpage
\begin{verbatim}
# FILE varlGenQfnmpl.txt

Qu0fn0:= proc(a,beta)   # simplest as in KM93
 local Q0, ecc, L;
 Q0:= Pi/(4*(1/a^2 + a^2));
 ecc:= sqrt(1 - 1/a^4);
 L:= 4*a*EllipticE(ecc);
 Q0 + beta*Pi^2/L
 end proc:

Qu0fn:= proc(a,beta)    # quadratic test functions
 local ecc, Eecc, Kecc, dQu0s, nQu0s;
 ecc:= sqrt(1 - 1/a^4);
 Eecc:= EllipticE(ecc);
 Kecc:= EllipticK(ecc);
 dQu0s:= 64*a*(4*a^12-15*a^8+4-15*a^4)*Eecc^2+(512*a*(a^8+2*a^4+1)*Kecc+
         720*a^2*Pi*beta*(1-2*a^4+a^8))*Eecc-320*a*(a^4+1)*Kecc^2;
 nQu0s:= 16*a^3*(-19*a^4+4*a^8+4)*Pi*Eecc^2+
    (128*a^3*(a^4+1)*Pi*Kecc+3*Pi^2*beta*(-55*a^4-55*a^8+23+
    23*a^12))*Eecc-80*a^3*Pi*Kecc^2+
    24*Pi^2*beta*(7+7*a^8-6*a^4)*Kecc+180*Pi^3*beta^2*a*(1-2*a^4+a^8);
 nQu0s/dQu0s
 end proc:

Qfn:= proc(a,beta)     # both the above
 local Q0, ecc, Eecc, Kecc, QKM93, nQu0, dQu0beta0, dQu0;
 Q0:= Pi/(4*(1/a^2 + a^2));
 ecc:= sqrt(1 - 1/a^4);
 Eecc:= EllipticE(ecc);
 Kecc:= EllipticK(ecc);
 QKM93:= Q0 + beta*Pi^2/(4*a*Eecc);
 nQu0:= (4*(a^4+1)*Kecc+(a^8-10*a^4+1)*Eecc)^2;
 dQu0beta0:= (4*a^8-19*a^4+4)*Eecc^2+8*(a^4+1)*Eecc*Kecc-5*Kecc^2;
 dQu0:= 45*a^2*Pi*beta*(a^4-1)^2*(a^4+1)*Eecc+ 4*a*(a^4+1)^2*dQu0beta0;
 [QKM93, QKM93+ 5*Pi^2*beta*nQu0/(16*Eecc*dQu0)]
 end proc:
\end{verbatim}

\clearpage

\section{$\QQ$ for nearly circular ellipses}\label{sec:nearCircQs}

\subsection{$\QQ$ at $\beta=0$}\label{subsec:nearCircQbeta0}

Now, at $\beta=0$, the explicit solution as given in \S\ref{sec:explicit}, yields asymptotics for $\epsilon\rightarrow{0}$:
\begin{align}
\QQ
&= 
\frac{\pi}{4(a^2+a^{-2})}  = \frac{\pi}{8}\, \sqrt{1-\epsilon^2} ,
\nonumber \\
&\sim \frac{\pi}{8}\left(1 - 2(a-1)^2+ O( (a-1)^3)\ \right)
\quad
{\rm as\ }\ a\rightarrow 1  ,
\nonumber \\
&\sim \frac{\pi}{8}\left( 1 - \frac{1}{2}\epsilon^2+ O( \epsilon^3)
\ \right)
\quad
{\rm as\ }\ \epsilon\rightarrow 0,
\label{eq:PexactEps}
\end{align}

\subsection{Nearly circular ellipses with $\beta\ge{0}$}\label{subsec:nearCircbetaGt0}

The explicit solution available for when $\beta=0$ is presented above,
and is useful to check against our asymptotics for
the solution when $\beta\ge{0}$.
When $\epsilon\rightarrow{0}$, equation~(\ref{eq:Qellb0}) agrees with the asymptotics given earlier in
equation~(\ref{eq:PexactEps}).

Concerning the ellipse when $\beta>0$,
when $\epsilon$ is small,  $u$ can be approximated in the form
\begin{equation}
 u =  \frac{1}{4}\left(  1-r^2 \right) + \frac{\beta}{2} +\epsilon^2\, t_{02} +
\epsilon \, t_{11} r^2 \cos(2\theta) +
\epsilon^2  \, t_{22} \, r^4\, \cos(4\theta) .
\label{eq:ueps}
\end{equation}
Details of the perturbation analysis are in~\cite{KW16}.
The result of the perturbation analysis is
\begin{align*}
t_{11}
&= \frac{1}{4}\, \frac{1+\beta}{1+2\beta}
 , \\
t_{02}
&=  -\frac{1}{32}\, \frac{4 + 5\beta +6\beta^2}{1+2\beta}
=  -\frac{1}{32}\,\left( 1+3\beta + \frac{3}{1+2\beta}\right)
 , \\
 t_{22}
&=   -\frac{1}{32}\, \frac{\beta (1-2\beta)}{(1+4\beta)(1+2\beta)}
 .
\end{align*}
Integrating $u$ with these parameters over the ellipse gives the
expansion for $\QQ({\rm ellipse})$:
\begin{eqnarray}
\QQ({\rm ellipse})
&\sim& \frac{\pi}{8} (1+4\beta) + q_1\epsilon^2
\label{eq:QsepsSmall}\\ 
 q_1
 &=& - \frac{\pi}{16} \left(1 - 8 t_{11} - 16 t_{02}\right)
 =  - \frac{\pi}{16} \left(1 +{\frac {\beta\, \left( 1+6\,\beta \right) }{2(2\,\beta+1)}}\right)  .
 \label{eq:q1Def}
\end{eqnarray}

Maple code for the near-circular approximation is as follows:
\begin{verbatim}
# FILE nearCircApproxQfnmpl.txt
nearCircApproxQfn:= proc(a,beta)
 local eps,q1;
 eps:= (a^2-1/a^2)/(a^2+1/a^2);
 q1:= -(Pi/16)*(1+beta*(1+6*beta)/(2*(2*beta+1)));
 Pi*(1+4*beta)/8 +q1*eps^2
 end proc;
\end{verbatim}

\subsubsection*{Near circular, $\beta$ small}

When $\beta={0}$, this agrees with the asymptotics given in
equation~(\ref{eq:PexactEps}).
When $\beta$ is small but non-zero, asymptotics given in
\S\ref{subsec:betaSmallEll}
check: see equation~(\ref{eq:Q1bSmallNC}).

\subsubsection*{Near circular, $\beta$ large}

When we first take $\epsilon$ small, then let $\beta$ tend to infinity
we have the approximation
$$ Q\sim \frac{\pi}{2}\beta +\frac{\pi}{8}
-\frac{3\pi}{32}\beta\epsilon^2 -\frac{\pi}{32}\epsilon^2 .
$$
This is used as a check in~\S\ref{subsec:Qinf}.

\cmdvmcode{NearCirc/nearCircCartCheckmpl.txt}{ }

\clearpage

\subsubsection*{Numerical results}

{
\enlargethispage*{1cm}
Table~\ref{tbl:nearCirc} below compares the previous asymptotics (column A) with variational lower bound (column V) and with the Fourier series solution (column F).
\nopagebreak[4]
\vspace*{-0.5cm}

{\small
\noindent
\begin{table}[h]
\begin{center}
\begin{tabular}{cccc}
\hline
$a=5/4$ & & \\
$\beta$& A& V& F\\
1/64  & .3825119944 & .3807427556 &.3807440330\\
1/16  & .4551128646 & .4530404184 &.4530434724\\
 1/4  & .7437766793 & .7406928885 &.7406944966\\
   1  & 1.888863782 & 1.882275712 &1.882277797\\
   4  & 6.450075877 & 6.429038519 &6.429053564\\
  16  & 24.68100694 & 24.60087273 &24.60089638\\
  64  & 97.59955268 & 97.28234644 &97.28237296\\
\hline
$a=17/16$ & & \\
$\beta$& A& V& F\\
      1/64& .4143605368 & .4143493072 &.4143493804 \\
      1/16& .4879061200 & .4878930227 &.4878930894\\
      1/4 & .7819440804 & .7819247200 &.7819248140\\
      1   & 1.957301880 & 1.957260968 &1.957261120\\
      4   & 6.657144997 & 6.657014853 &6.657015098\\
      16  & 25.45536249 & 25.45486728 &25.45486742\\
      64  & 100.6478027 & 100.6458431 &100.6458432\\
\hline
$a=65/64$ & & \\
$\beta$& A& V& F\\
      1/64&  .4170525380 & .4170525625 &.4170524896\\
      1/16& .4906779729 & .4906782880 &.4909655446\\
      1/4 & .7851701838 & .7851700411 &.8350137728\\
      1   & 1.963086617 & 1.963085377 &1.963085489\\
      4   & 6.674647535 & 6.674645198 &6.675381808\\
      16  & 25.52081497 & 25.52081128 &25.51131048\\
      64  & 100.9054563 & 100.9054466 &100.9054478\\
\hline
\end{tabular}
\end{center}
\label{tbl:nearCirc}
\end{table}
} 
}

\begin{figure}[hb]
\centerline{\includegraphics[height=7cm,width=13cm]{\grpath/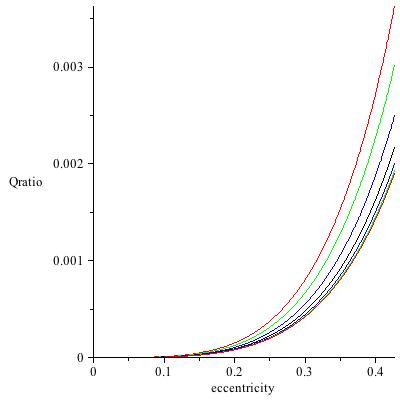}}
\caption{\S\ref{subsec:nearCircbetaGt0}. Plot of Qratio=$(Q_\odot-Q(e,\beta))/Q_\odot$
 at various values of $\beta$ for small eccentricity.}
\label{fig:QratioEpsSmalljpg}
\end{figure}

\clearpage

\noindent{\bf A simpler, but less accurate, variational approach}

\begin{figure}[hb]
\centerline{\includegraphics[height=7cm,width=13cm]{\grpath/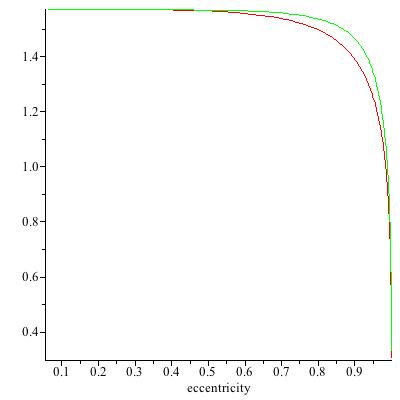}}
\caption{\S\ref{subsec:nearCircbetaGt0}. Plot of small $\beta$ approximation term $Q_1$ 
and a variational lower bound.}
\label{fig:Q1bndBetaSmalljpg}
\end{figure}

The lower bound shown in Figure~\ref{fig:Q1bndBetaSmalljpg} comes from a trial function which is just the $\beta=0$ quadratic in Cartesians plus the (best) constant.
Of course this is exact for eccentricity zero.
Of course we do better with a more general quadratic but the integrals are more ugly.

\clearpage

\medskip

\section{Flows with $\beta$ small}\label{sec:betaSmall}

\subsection{General $\Omega$}\label{subsec:betaSmallGenOm}

The asymptotics are described in~\cite{KM93}. The problem is
\begin{eqnarray*}
- \Delta u
&=& 1 \qquad {\rm in}\ \ \Omega ,\\
\beta \frac{\partial u}{\partial n} + u
&=& 0 \qquad {\rm on}\ \ \partial\Omega .
\end{eqnarray*}
When $\beta$ is small, the solution $u_\beta$ is asymptotically
$$
u_\beta
\sim u_{\beta=0} + \sum_{j=1}^n \beta^j u_j ,
$$
with each $u_j$ ($j\ge{1}$) harmonic (and independent of $\beta$).
Our interest is in $\QQ$, here denoted $Q_\beta$,
$$ Q_\beta = \int_\Omega u_\beta , $$
so
$$ Q_\beta \sim Q_0 +\beta Q_1\qquad{\rm where}\ \ Q_1=\int_\Omega u_1 , $$
and the integral of $u_1$ over $\Omega$ can be found explicitly
(at least up to a single-variable integral)
even when $u_1$ cannot be found explicitly.
Now
\begin{eqnarray*}
- \Delta u_1
&=& 0 \qquad {\rm in}\ \ \Omega ,\\
\beta u_1
&=& -\beta \frac{\partial u_0}{\partial n}  \qquad {\rm on}\ \ \partial\Omega .
\end{eqnarray*}
The Divergence Theorem gives
\begin{eqnarray*}
Q_1=\int_\Omega u_1
&=& \int_\Omega\left( u_0 \Delta u_1 - u_1\Delta u_0\right)
= \int_\Omega {\rm div}\left( u_0\nabla u_1 - u_1\nabla u_0\right) ,\\
&=& \int_{\partial\Omega}\left( u_0 \frac{\partial u_1}{\partial n}. - u_1 \frac{\partial u_0}{\partial n}\right)
= \int_{\partial\Omega}\left( \frac{\partial u_0}{\partial n}\right)^2 \, ds .
\end{eqnarray*}
We remark
$$\int_{\partial\Omega}\, ds=|\partial\Omega|, \qquad
\int_{\partial\Omega}\left( \frac{\partial u_0}{\partial n}\right) \, ds= -|\Omega| .$$
This yields, via Cauchy-Schwarz,
$$ |\Omega|^2 \le |\partial\Omega|\  Q_1 , $$
an inequality which becomes more accurate the closer the ellipse is to a circle.

\subsection{$\Omega$ an ellipse}\label{subsec:betaSmallEll}

In general, and in the notation of the preceding subsection,
$$ Q_1
= \int_{\partial\Omega}\left( \frac{\partial u_0}{\partial n}\right)^2 \, ds
= \int_{\partial\Omega}\, |\nabla u_0|^2 \, ds .
$$
For an ellipse $u_0=u_{0e}$ with $u_{0e}$ given in~\S\ref{sec:explicit}:
thus
$$Q_{1,e}
= \int_{\partial\Omega}\, |\nabla u_{0e}|^2 \, ds .
$$
Although the calculation can be done in polar coordinates, perhaps Cartesians lead to simpler integrals:
$$  |\nabla u_{0e}|^2
= \frac{ {a^{-4}}x^2 + a^4 y^2}{(a^{-2}+a^2)^2} .
$$
The Cartesian equation for that part of the ellipse in the first quadrant is
$$ Y(x)= \frac{1}{a}\, \sqrt{1-(\frac{x}{a})^2} , \qquad{\rm for\ \ } 0<x<a .
$$
The perimeter is
\begin{eqnarray}
 |\partial\Omega|
&=& 4 \int_0^a \sqrt{1+ Y'(x)^2}\, dx ,
\nonumber\\
&=& 4 a {\rm EllipticE}\left(\frac{\sqrt{a^4-1}}{a^2}\right) \quad{\rm for\ \ } a\ge{1} ,
\nonumber\\
&=& 4 a {\rm EllipticE}(e) \quad{\rm for\ \ } a\ge{1} ,
\label{eq:perimE}\\
&=& 4 (1-e^2)^{-1/4}\, {\rm EllipticE}(e) ,
\nonumber
\end{eqnarray}
where $e$ is the eccentricity as defined in equation~(\ref{eq:eaDef}).\\

We now calculate $Q_1$, beginning with
$$  |\nabla u_{0e}|^2
= \frac{ a^2  - a^{-2} (a^2 -a^{-2}) x^2}{(a^{-2}+a^2)^2}  \qquad{\rm on\ \ } y=Y(x) .
$$
On integrating this $Q_1$ is given by
\begin{equation}
 Q_1
= \frac{4}{3}\, \frac{a^3}{(1+a^4)^2}\left( 2(1+a^4){\rm EllipticE}(e)- {\rm EllipticK}(e)\right)
\quad {\rm for\ \ } a\ge{1} .
\label{eq:Q1bSmall}
\end{equation}
There are many equivalent ways to write this.
With $i_2=i_{xx}+i_{yy}$ we have
$$ Q_1 = \frac{L}{a^2+\frac{1}{a^2}} -
\frac{i_2}{(a^2+\frac{1}{a^2})} . $$

\begin{verbatim}
# FILE QsmallBetampl.txt
betaSmallQfn:= proc(a,beta)
 local ecc, Q0, Q1;
 ecc:=sqrt(1-1/a^4);
 Q0:= Pi/(4*(a^2 +1/a^2));
 Q1:= (4/3)*(a^3/((1+a^4)^2))*(2*(1+a^4)*EllipticE(ecc)-EllipticK(ecc));
 Q0+beta*Q1
 end proc;

 
\end{verbatim}{ }

\subsubsection*{$\beta$ small, near circular}

The asymptotic expansion for nearly circular ellipse, $a$ near 1, has
\begin{equation}
Q_1
\sim \frac{\pi}{2} - \frac{\pi}{8} (a-1)^2\qquad {\rm for\ \ } a\rightarrow{1} .
\label{eq:Q1bSmallNC}
\end{equation}
This agrees with the earlier expansion given in equations~(\ref{eq:QsepsSmall},\ref{eq:q1Def}) where,
for $\beta$ small
$$ \QQ \sim \frac{\pi}{8}(1+ 4\beta)- \frac{\pi}{4}(1+ \frac{\beta}{2}) (a-1)^2
\ \quad{\rm for }\ \ a\rightarrow{1}, \ \beta\rightarrow{0} . $$
\medskip

\subsubsection*{Numerical results}

In Table~\ref{tbl:betaSmall} below we compare with the numerical values given by the lower-bound variational winner amongst quadratics.

\begin{table}[h]
\begin{center}
\begin{tabular}{cccc}
\hline
$\beta=1/4$& & & \\
$a$& A& V& F\\
  65/64 &   .7851858133   &   .7851702110   &   .8350137728  \\
  17/16 &   .7821607512   &   .7819248090   &   .7819248140  \\
  5/4   &   .7432200201   &   .7406928855   &   .7406944966  \\
    2   &   .4953609194   &   .4909982744   &   .4910897996  \\
    4   &   .2147581278   &   .2139371008   &   .2142092926  \\
   16   &   .4473378672e-1 &  .4463047672e-1 &  .4473097574e-1 \\
\hline
$\beta=1/16$& & & \\
$a$& A& V& F\\
65/64  &  .4906792259  &  .4906780838  &  .4909655446 \\
17/16  &  .4879127252  &  .4878930761  &  .4878930894 \\
5/4    &  .4532504638  &  .4530404192  &  .4530434724 \\
  2    &  .2624399058  &  .2620398525  &  .2621170972 \\
  4  & .9036181975e-1  &  .9023065471e-1 & .9032360830e-1 \\
 16  & .1348438275e-1  &  .1345885965e-1 & .1348419675e-1 \\
\hline
$\beta=1/64$& & & \\
$a$& A& V& F\\
  65/64  &  .4170525792  &  .4170526698  &  .4170524896 \\
  17/16  &  .4143507188  &  .4143493767  &  .4143493804 \\
5/4  &  .3807580746  &  .3807427560  &  .3807440330 \\
  2  &  .2042096524  &  .2041611381  &  .2041881644 \\
  4  &   .5926274272e-1  &   .5923472697e-1  &   .5926025316e-1 \\
 16  &   .5672031762e-2  &   .5665670014e-2  &   .5671955380e-2 \\
\hline
\end{tabular}
\end{center}
\label{tbl:betaSmall}
\end{table}

\clearpage

\subsection{$u_1$ for $\beta$ small: Fourier series?}

We begin by recalling (ref{eq:12b})
$$ v_p\;=\; -\frac{c^2}{8} (\cosh(2\eta)+\cos(2\psi))
=  -\frac{1}{4} (x^2+y^2)  .$$
Several items relating to $v_p$ are used, in particular its normal derivative at the boundary,
the $\eta$-derivative being closely related.
\cmdvmcode{vpEta0mpl.txt}{ }

The diagonal matrix $M_0$ has as its $(n,n)$ entry
$$m_{0,0}= 2\pi , \quad{\rm and\ for\ }n\geq{1}\quad
m_{0,n}=\pi\,\cosh(2n\eta_0).$$
The first subscript, here $0$, reminds us that the entries come from $M_0$.
The vector $F_0$ has nonzero entries only at $n=0$ and $n=1$ when
$$ f_{0,0}= \frac{c^2\pi}{4}\, \cosh(2\eta_0) \quad{\rm and\ } f_{0,1}= \frac{c^2\pi}{8} . $$
Again, the first subscript, here $0$, reminds us that the entries come from $f_0$.

The symmetric matrix $M_1$ is full with entries involving EllipticK and EllipticE.
The vector $f_1$, related to the Fourier series of $g$, similarly involves  EllipticK and EllipticE.

For $\beta$ positive but small, we seek a representation of the form
$$u \sim u_0 + \beta \, u_1\qquad{\rm for}\ \ \beta\rightarrow{0} . $$
Our main interest is checking with an asymptotic result obtained using Cartesian coordinates:
$$Q \sim Q_0 + \beta \, Q_1\qquad{\rm for}\ \ \beta\rightarrow{0} . $$
Our calculation in Cartesian coordinates was completed a year earlier than the elliptic coordinate work here.

\subsubsection*{$\beta=0$ re-visited}

As one small check on the equations determining the $A_n$, we note that when $\beta=0$ the equations are very easy to solve.
The matrix $M=M_0$ is diagonal.
The only nonzero $A_n$ are  those corresponding to $n=0$ and $n=1$..

When $\beta=0$, and subscripting with the first 0 in a pair merely reminding us that the values are for $\beta=0$,
$u_0= v_p +A_{0,0} +A_{0,1}\cosh(2\eta)\cos(2\psi)$ and the equation $M_0 A= f_0$ gives
\begin{eqnarray*}
A_{0,0}
&=& \frac{f_0}{m_{0,0}}= \frac{c^2\, \cosh(2\eta_0) }{8}=\frac{\cosh(2\eta_0)}{4\sinh(2\eta_0)} \\
 A_{0,1}
 &=& \frac{f_1}{m_{0,1}}=   \frac{c^2}{8\cosh(2\eta_0)  }= \frac{1}{4\cosh(2\eta_0)\sinh(2\eta_0)} .
\end{eqnarray*}

As we have already calculated $Q_0$ from other methods, the  calculation here merely serves as a
check on $A_{0,0}$ and $A_{0,1}$.

\subsubsection*{$\beta>0$ small}

Consider now approximations when $\beta$ is small.
In \S{\ref{sec:betaSmall}} we find
$$ Q\sim Q_0 +\beta\, Q_1\qquad{\rm as}\ \ \beta\rightarrow{0} , $$
with $Q_0$ as in equation~(\ref{eq:Qellb0}).
The approach began with assuming
$$ u\sim u_0 +\beta\, u_1\qquad{\rm as}\ \ \beta\rightarrow{0} , $$
but did not require finding $u_1$.
If one wished to find $u_1$, a Fourier series approach as in this section might be appropriate.
Approximate the vector $A$ of Fourier coefficients by
$$ A \sim A_0 + \beta \, A_1 , $$
where the vector $A_0$ is that already calculated in the preceding subsection concerning $\beta=0$.
Then we require, to order $\beta$,
$$ (M_0+\beta M_1)(A_0 + \beta A_1) = (f_0+ \beta f_1) , $$
i.e., on approximating to order $\beta$,
$$ M_0 A_1 = f_1 -M_1 A_0 . $$
Equivalently, using, as noted previously, that $M_0$ is diagonal,
\begin{eqnarray*}
m_{0,0} A_{1,0}
&=& f_{1,0} - m_{1,0,0} A_{0,0} - m_{1,0,1} A_{0,1},\\
 m_{0,1} A_{1,1}
 &=&f_{1,1} -m_{1,1,0} A_{0,0} - m_{1,1,1} A_{0,1} \\
 \vdots
 &=& \vdots \\
 m_{0,n} A_{1,n}
 &=&f_{1,n} -m_{1,n,0} A_{0,0} - m_{1,n,1} A_{0,1}
\end{eqnarray*}
Again the first subscript denotes the appropriate order in $\beta$.

While the calculations are routine,  it may be that the solution $A_1$ is ugly, and,
in any event, the human effort at keeping all the details correct is  considerable.
We suspect that, at this stage in the study and use of microchannels, the greater accuracy for $u$ is not needed.
A calculation which would provide an additional check is
to finding $A_{1,0}$ and $A_{1,1}$ which is all that is needed to check against the
result for $Q_1$ given in equation~(\ref{eq:Q1bSmall}).\\

\clearpage

\section{$\beta$ large}\label{sec:betaLarge}

\subsection{The dominant term}

For general $\Omega$
$$Q(\beta)\sim \beta\frac{|\Omega|^2}{|\partial\Omega|}
\qquad{\rm for\ \ } \beta\rightarrow\infty .
$$
For our ellipses of area $\pi$ the perimeter is given by equation~(\ref{eq:perimE}),
Hence
$$ Q(\beta)
\sim \beta\frac{\pi^2}{4 a {\rm EllipticE}(e)}
\qquad{\rm for\ \ } \beta\rightarrow\infty .
$$

\medskip
Maple code for this is:
\begin{verbatim}
betaLargeQfn:= (a,beta) -> beta*Pi^2/(4*a*EllipticE(sqrt(1-1/a^4)));
\end{verbatim}

\subsubsection*{$\beta$ large, nearly circular ellipses}

Further asymptotic approximation of that above for nearly circular ellipses is:
$$
Q(\beta)
\sim \beta\left(\frac{\pi}{2} - \frac{3\pi}{8} (a-1)^2 + O((a-1)^3)\quad \right)
\qquad{\rm for\ \ } \beta\rightarrow\infty, \ a\rightarrow 1 .
$$

This may be compared with the asymptotics found in
equations~(\ref{eq:QsepsSmall},\ref{eq:q1Def}) further approximated for $\beta$ large.
Then $q_1\sim -\pi\beta/32$ so
$$\QQ
\sim \frac{\pi\beta}{2}+q_1\epsilon^2
\sim \frac{\pi\beta}{2} - \frac{3\pi\beta}{32}\epsilon^2
\sim \frac{\pi\beta}{2} - \frac{3\pi\beta}{8}(a-1)^2 ,$$
which agrees with the result of the preceding paragraph.

We have yet to find higher order approximations.
\medskip

The variational approximation using quadratic test functions agrees with the lowest order term.
\ifthenelse{\boolean{vmcode}}{
\newline
assume(ap $>$ 1);\\
limit(Qvals(ap, beta)/beta, beta = infinity);
}{ }

The first term $O(\beta)$ is, at fixed $\beta$  a decreasing function of $a$ on $a>1$.
Note that the quadratic lower bound in column V in the following table
is a better approximation than
$\beta |\Omega|^2/|\partial\Omega|$ which is less than the lower bound.\\

\newpage
\subsubsection*{Numerical results $\beta$ large} 

\begin{table}[h]
\begin{center}
\begin{tabular}{cccc}
\hline
$\beta=4$& & & \\
$a$& A& V& F\\
 65/64 & 6.282052744 & 6.674647024 & 6.675381808\\
 17/16 & 6.265911941 & 6.657014849 & 6.657015098\\
   5/4 & 6.056752387 & 6.429038511 & 6.429053564\\
2 & 4.602060688 & 4.869625399 & 4.870154656\\
4 & 2.449872080 & 2.627235185 & 2.627747294\\
    16 & .6168200091 & .6676825774 & .6691074764\\
\hline
$\beta=16$& & & \\
$a$& A& V& F\\
65/64 & 25.12821097 & 25.52081287 & 25.51131048\\
 17/16 & 25.06364776 & 25.45486727 & 25.45486742\\
   5/4 & 24.22700955 & 24.60087275 & 24.60089638\\
2 & 18.40824275 & 18.69310262 & 18.69427704\\
4 & 9.799488320 & 10.09964112 & 10.10542728\\
    16 & 2.467280036 & 2.656743920 & 2.661123514\\
\hline
$\beta=64$& & & \\
$a$& A& V& F\\
 65/64 & 100.5128439 & 100.9054480 & 100.9054478\\
 17/16 & 100.2545911 & 100.6458431 & 100.6458432\\
   5/4 & 96.90803820 & 97.28234641 & 97.28237296\\
2 & 73.63297100 & 73.92334640 & 73.92489450\\
4 & 39.19795328 & 39.57699513 & 39.59760516\\
    16 & 9.869120146 & 10.54995058 & 10.55945235\\
\hline
\end{tabular}
\end{center}
\label{tbl:betaLarge}
\end{table}

Further work concerning large $\beta$ behaviour is given,
via Fourier series, in Part II, and, 
by variational techniques in Part IV.

\clearpage

\section{Conclusion and Open Problems}\label{sec:Conclusion}

The solution to the Robin boundary problem has been approached from several directions:
Fourier series, variational bounds, and asymptotic approximations.

There are several questions that remain.
\begin{enumerate}
\item What can be said about asymptotics when the eccentricity tends to 1?
The small $\beta$ and the large $\beta$ approximations yield some information but more should be possible.
(We remark that at large $\beta$, small eccentricity will
cause the Fourier coefficients $V_n$ of Part II to decay more slowly.)
Matched asymptotic approximations may be appropriate.

\item Can $u_\infty$ be found explicitly, or at least the integrals  $\Sigma_\infty$ and $\Sigma_1$?\\
Given that there are 3-term recurrence relations for the Fourier coefficients of both $g$ and $1/g$, might there be a recurrence relation for the coefficients $V_n$?
One might combine equation~(\ref{eq:Vninf})
with items from Part II \S~\ref{app:Furtherg}.
See also Part IV.

\item In the asymptotics for $\beta$ tending to zero, we found how $Q$ changed, but can one find tidy formulae giving
the departure of $u$ from $u_0$?

\item Might it be possible to make better use of the Fourier series for $g$ and $1/g$?

\end{enumerate}

In connection with items 2 and 3 above, we remark that there is
a `Poisson Integral Formula' for solving the Dirichlet problem for
Laplace's equation in an ellipse.
This follows from the conformal map between ellipse and disk:
see~\cite{Kob}, p177 and~\cite{Sz50}.
See also~\cite{Ro64,Mi90,LK14,SD16}.

\bigskip


\section{Appendix A: Geometry of ellipses}\label{sec:GeometryI}

\subsection{General ellipses}\label{subsec:GeomGen}

\subsubsection{Polar coordinates}
Consider the ellipse $\frac{x^2}{a^2}+a^2 y^2 \le{1}$.
In polar coordinates relative to the centre the boundary curve is
\begin{align}
 r
&= \frac{1}{\sqrt{\frac{\cos(\theta)^2}{a^2} + a^2 \sin(\theta)^2}}
= \sqrt{\frac{2}{{a}^{2}+{a}^{-2}-\left( {a}^{2}-{a}^{-2} \right) {\cos(2\theta)}}}
\nonumber \\
&= r\left(\frac{\pi}{4}\right) 
\left( 1 -
\frac{a^2-a^{-2}}{a^2+a^{-2}}
\cos(2\theta)
\right)^{-1/2}\ \
{\rm where}\ \  r\left(\frac{\pi}{4}\right) = \sqrt{\frac{2}{{a}^{2}+{a}^{-2}}}
\label{eq:rEllipsec2} \\
&=  \sqrt{\frac{1 +\tan^2(\theta)}{a^{-2} + a^2 \tan^2(\theta)} } \ .
\label{eq:rEllipset}
\end{align}
In our computations we take $a\ge{1}$.\\
Without any assumptions on $a$, the coefficients in the Fourier series for $r(\theta)$ involve elliptic integrals.
Our interest in some later sections will be  in $a$ near 1.
However $(a-1)$ might not be the best perturbation parameter.
We have also used, in computations,
$$\epsilon
=\frac{a^2-a^{-2}}{a^2+a^{-2}}
\sim 2(a-1) \qquad {\rm\ as\  }\ a\rightarrow{1}\  ,
$$
The ellipse is
$$ r^2= \frac{\sqrt{1-\epsilon^2}}{1- \epsilon\cos(2\theta)} .
$$
We have also used  (as have others, e.g.~\cite{Day55}),  the eccentricity,
as in equation~(\ref{eq:eaDef}),
\begin{equation}
 e=\sqrt{1-\frac{1}{a^4}} \sim 2(a-1)^{1/2}  \qquad {\rm\ as\  }\ a\rightarrow{1} .
\label{eq:eDef}
\end{equation}
As both $e^2$ and $\epsilon$ are relatively similar low-order rational functions of $a^2$,
one readily finds that
$$ \frac{\epsilon}{e^2}= \frac{1}{2}\, (1+\epsilon) .$$

\cmdvmcode{epsToempl.txt}{ }

We remark that the polar equation in terms of the eccentricity is
$$ r= \frac{1}{a\sqrt{1- (e\cos(\theta))^2}} = \frac{g(\theta)}{a}. $$
\cmdvmcode{polarCheckmpl.txt}{ }

Including an $a$ dependence in $r(a,\theta)$ we remark that $r(1/a,\theta)=r(a,\theta+\pi/2)$ and
$\epsilon(1/a)=-\epsilon(a)$.
The binomial expansions
\begin{align}
\frac{r(0)}{r(\pi/4)}
&= (1-\epsilon)^{-1/2}
= 1+ \sum_{k=1}^\infty \frac{ (2k)! }{2^{2k} (k!)^2}\ \epsilon^k ,
\nonumber \\
\frac{r(\theta)}{r(\pi/4)}
&= (1-\epsilon\cos(2\theta))^{-1/2}
= 1+ \sum_{k=1}^\infty \frac{ (2k)! }{2^{2k} (k!)^2}\ \cos(2\theta))^k \epsilon^k ,
\label{eq:rCosk}
\end{align}
may be useful in finding higher terms in the perturbation expansions (for $a$ near 1) of some domain functionals.
The symmetries of the ellipse explain the form of the expansion in~(\ref{eq:ueps}):
\begin{itemize}
\item It is symmetric about $\theta=0$, hence only cosine terms.
\item It is symmetric about $\theta=\pi/2$ and hence only the even order cosine terms  $\cos(2 m\theta)$.
\item When $\epsilon$ is replaced by $-\epsilon$ and $\theta$ by $\theta+\pi/2$ the expression is unchanged and hence the form of the polynomial coefficients in $\epsilon$ forming the Fourier coefficients.
For $m$ is odd, only odd powers of $\epsilon$ appear: for $m$ is even, only even powers of $\epsilon$ appear.
\end{itemize}
We see these symmetries in connection  the solutions of our pde problem.
Returning to the study of the boundary curve,
we also need the expansion for $r(\pi/4)$:
$$ a^2=\sqrt{\frac{1+\epsilon}{1-\epsilon}} , \qquad
\frac{2}{r(\pi/4)}= a^2 + a^{-2} =\frac{2}{\sqrt{1-\epsilon^2}},\qquad
r\left(\frac{\pi}{4}\right) = \left(1 -\epsilon^2 \right)^{1/4} .
$$

\medskip
There are a few items of undergraduate calculus that are used.
Let $s$ denote arclength measured around the curve.
Then
\begin{eqnarray*}
\frac{d s}{d\theta}
&=& \sqrt{ r(\theta)^2 +\left(\frac{d  r(\theta)}{d\theta}\right)^2}\qquad
{\rm in \ general}, \\
&=& 2 a\sqrt{\frac{ (1+a^8)-(a^8-1)\cos(2\theta)}
{((1+a^4)-(a^4 -1)\cos(2\theta))^3}}\qquad
{\rm for \ an \ ellipse} , \\
&=& 2 a\sqrt{\frac{ (1+\epsilon)(1+\epsilon^2-2\epsilon\cos(2\theta))}
{(1-\epsilon\cos(2\theta))^3}}\qquad
\ \ \  {\rm for \ an \ ellipse} .
\end{eqnarray*}

We remark that, for the ellipse the normal can 
be found from the gradient of the ellastic torsion function,
$\nabla u_{0e}  = \nabla u({\rm ellipse},\beta=0)$.

\subsubsection{The usual parametric description of an ellipse}

The boundary of the ellipse can be described by
\begin{equation}
x=a\cos(\psi),\qquad y=\frac{1}{a}\sin(\psi) .
\label{eq:param}
\end{equation}
(The relation between the parameter $\psi$ and the polar angle
$\theta$ is $\tan(\theta)=\tan(\psi)/a^2$.)
\begin{eqnarray*}
(\frac{ds}{d\psi})^2
&=& (\frac{dx}{d\psi})^2 + (\frac{dy}{d\psi})^2 ,\\
&=& a^2 \sin(\psi)^2 + \frac{1}{a^2} \cos(\psi)^2 ,\\
&=& a^2\left( 1 - e^2 \cos(\psi)^2 \right)
\end{eqnarray*}
From this the perimeter is calculated:
\begin{eqnarray*}
|\partial\Omega|
&=& 4 a \int_0^{\pi/2}
\sqrt{1- (e\cos(\psi))^2}\, d\psi ,\\
&=& 4 a \int_0^{\pi/2}
\sqrt{1- (e\sin({\hat\psi}))^2}\, d{\hat\psi} ,\\
&=& 4 a\, {\rm EllipticE}(e) .
\end{eqnarray*}

\bigskip

\pagebreak[3]

\par\noindent{\bf The support function, $h$}

There are several equivalent definitions of the
support function $h$ which we consider as a function
of points on the boundary.
Thus at the point associated with the value $\psi$
in definition~(\ref{eq:param})
$$ h =\frac{1}{a\sqrt{1+e^2\cos^2(\psi)}} . $$
As noted in~\cite{Kim07}, the curvature is $h^3$.
Once again, the functions $g$ and $\hat g$
make an appearance.

\subsubsection{Elliptic coordinates}

The elliptic coordinates are related to Cartesians by
$$  x= c\cosh(\eta)\cos(\psi), \qquad
y= c\sinh(\eta)\sin(\psi), $$
and we will set the parameter $c$ by
$$ c= \sqrt{a^2-a^{-2}} . $$

The boundary of the ellipse can be represented, with fixed $\eta_0$,
$$ x= c\cosh(\eta_0)\cos(\psi), \qquad
y= c\sinh(\eta_0)\sin(\psi). $$
In this
$$ a= c\cosh(\eta_0), \qquad  a^{-1}= c\sinh(\eta_0) , $$
so that
$$ \tanh(\eta_0) = a^{-2},\quad (c^2=a^2-a^{-2}),\quad
\epsilon= \frac{1}{\cosh(2\eta_0)} .
$$
Nearly circular ellipses will have $\eta_0$ large.

The perimeter is calculated:
\begin{eqnarray*}
|\partial\Omega|
&=& 4 \int_0^{\pi/2} \sqrt{ (\frac{d x}{d\psi})^2 +(\frac{d y}{d\psi})^2}\, d\psi ,\\
&=& 4 c \int_0^{\pi/2} \sqrt{ (\cosh(\eta_0))^2 - (\cos(\psi))^2} \, d\psi , \\
&=& 4 c \cosh(\eta_0) \, {\rm EllipticE}\left(\frac{1}{\cosh(\eta_0)}\right) , \\
&=& 4 a\, {\rm EllipticE}(e) .
\end{eqnarray*}
We remark that toroidal functions occur elsewhere in this study
and that~\cite{Ab09} gives a formula, symmetric in $a$ and $b$,
\begin{equation}
|\partial\Omega|
= 2\pi \sqrt{a\,b}\, 
{\rm LegendreP}_{1/2}\left(\frac{a^2+b^2}{2 a\, b}\right) .
\label{eq:ab09P}
\end{equation}

\subsection{Nearly circular ellipses, polar coordinates}\label{subsec:GeomNearCirc}

The first few terms in the expansion of $r(\theta)$ are
$$ r(\theta)
= 1+\frac{1}{2}\,\cos \left( 2\,\theta \right) \epsilon+ \left( \frac{3}{16}\,\cos \left( 4\,\theta \right) -\frac{1}{16} \right) {\epsilon}^{2} +
 O(\epsilon^3) .
$$
Alternatively, we can consider asymptotics as $a\rightarrow{1}$.
Then $\rho=r_{\rm ellipse}-1$ satisfies
{\small
\begin{align}
\rho(\theta)
&\sim (a-1) \cos(2\theta) - (\frac{1}{4} +\frac{1}{2} \cos(2\theta)-\frac{3}{4} \cos(4\theta))  (a-1)^2 + O(  (a-1)^3)
\nonumber \\
&= -\frac{1}{4} (a-1)^2 +  \left((a-1)-\frac{1}{2}  (a-1)^2\right) \cos(2\theta)+\frac{3}{4}  (a-1)^2\cos(4\theta)  + O(  (a-1)^3)   .
\label{eq:ellipseRho}
\end{align}
}

\medskip
For calculations extending the use of higher order terms in the Fourier series~(\ref{eq:ellipseRho}).
one might need to replace $x^{2n}$ terms using formulae like
 $$   x^{2n}=2^{1-2n}\left(\frac{1}{2} \binom{2n}{n}+\sum_{j=1}^n \binom{2n}{n-j}T_{2j}(x) \right),
 $$
and, at some future date, we may  for elliptical $\Omega$ use this but have yet to do so.

\medskip
\subsection{Moments of inertia}\label{subsec:Ic}

\subsubsection{General $\Omega$}

\noindent{\bf Area moments}

The polar moment of inertia, taking the origin at the centroid, is
$$ I_c = \int_\Omega \left( (x-x_c)^2 + (y-y_c)^2 \right) ,
$$
where $z_c=(x_c,y_c)$ is the centroid of $\Omega$.
When the boundaries are given in polar coordinates, this is
$$I_c =\frac{1}{4}  \int_0^{2\pi} r(\theta)^4\ d\theta .
$$

Higher order moments arise in connection with calculations
based on polynomial test functions in the variational approach 
of~\S\ref{sec:Varl}.

\medskip\par\noindent{\bf Boundary moments}

See~\S\ref{sec:Varl}.

\subsubsection{Ellipse}

\noindent{\bf Area moments}

For our disk and ellipse these are
\begin{equation}
I_c({\rm disk})= \frac{\pi}{2} a^4, \qquad
I_c({\rm ellipse})= \frac{\pi}{4} (a^2 + a^{-2})
=\frac{\pi}{2\sqrt{1-\epsilon^2}} .
\label{eq:icEllipseEps}
\end{equation}
The asymptotics below check with the entries in the 
table of~\cite{PoS51} treating domain functionals for
nearly circular domains:
\begin{align*}
\frac{2 I_c({\rm ellipse})}{\pi}
&\sim 1 + 2(a-1)^2 \qquad{\rm as\ } a\rightarrow 1 ,\\
\left(\frac{2 I_c({\rm ellipse})}{\pi}\right)^{1/4}
& \sim 1 + \frac{1}{2}(a-1)^2  + o(a-1)^2 ) \qquad{\rm as\ } a\rightarrow 1  ,\\
& \sim 1 + a_0 + \frac{3}{4}a_2^2    + o(a-1)^2 )\qquad{\rm as\ } a\rightarrow 1 .
\end{align*}

\medskip\par\noindent{\bf Boundary moments}

For doubly-symmetric domains like our ellipse, our notation 
in~\S\ref{sec:Varl} is
$$ i_{2n} = \int_{\partial\Omega} r^{2n}\, ds . $$
For the ellipse, these can be evaluated in terms of
elliptic integrals.

\subsection{Ellipse geometry: miscellaneous}\label{subsec:ellipseMisc}

The modulus of asymmetry for an ellipse is calculated in~\cite{Fi14}.
The calculation involves calculating ${\hat g}_0$.

In the case $\beta=0$ there are improvements 
to the St Venant inequality
in terms of the modulus of asymmetry.
(For the St Venant inequality, see~\cite{PoS51}, and, 
for ellispes, trivially the $\beta=0$
case in inequality~(\ref{eq:isoperQ}).)
\bigskip

\section{Appendix B: Elliptic integrals}\label{appx:EllipticInts}

The functions denoted by EllipticE and EllipticK in this work
are, in this appendix, written as $E(k)$ and $K(k)$.
Our notation follows that of many (but not all!) authors,
in particular~\cite{Law89}\S3.8 and~\cite{GR}:
\begin{eqnarray*}
E(k)
&=& \int_0^{\pi/2} \sqrt{1-k^2\,\sin(\theta)^2}\, d\theta,\\
K(k)
&=& \int_0^{\pi/2} \frac{1}{\sqrt{1-k^2\,\sin(\theta)^2}}\, d\theta .
\end{eqnarray*}

\medskip

For $0\le{z}<1$
$$ 
K(k)=\int_0^1 \frac{1}{\sqrt{(1-t^2)(1-k^2 t^2)}}\, dt , $$
$$ 
E(k)_=\int_0^1\sqrt{ \frac{1-k^2 t^2}{1-t^2}}\, dt . $$
The Cauchy-Schwarz inequality yields
$$\frac{\pi^2}{4} \le K(k)\, E(k) .$$

The functions can be written as hypergeometric functions:
\begin{eqnarray*}
K(k)
&=&\frac{\pi}{2} {\ }_2{F}_1(\frac{1}{2},\frac{1}{2};1; k^2) ,
\\
E(k)
&=& \frac{\pi}{2} {\ }_2F_1(-\frac{1}{2},\frac{1}{2};1; k^2) .
\end{eqnarray*}

In showing, in Part II \S\ref{app:Furtherg}
equation~(\ref{eq:gHathy}), that two different representations of
${\hat g}_n$ are equivalent, we needed to tell maple that,
as in~\cite{GR} 8.126, with $k'=\sqrt{1-k^2}$,
(and, in a different notation in~\cite{AS}17.3.29)
\begin{eqnarray*}
K\left(\frac{1-k'}{1+k'}\right)
&=& \frac{1+k'}{2}\, K(k) ,
\\
E\left(\frac{1-k'}{1+k'}\right)
&=& \frac{1}{1+k'}\, \left(E(k)+k'\, K(k)\right) .
\end{eqnarray*}
(The first of these is given, in a different notation, on
{\tt functions.wolfram.com}.)
In our application $k=1/\sqrt{q}=e$.

Another identity (which we will see in connection with
toroidal functions equation~(\ref{eq:niceW0})) is
`Legendre's relation' (\cite{AS}17.3.13):
\begin{equation}
K(k) E(k') - K(k) K(k') +K(k') E(k) = \frac{\pi}{2} .
\label{eq:LegRel}
\end{equation}

\subsection{Complex k}

The EllipticK function has a mirror symmetry:
$$ K({\bar k}) = {\bar K(k)} .$$
\medskip

The real and imaginary parts of $K(\exp(i\phi))$ are given,
for example, at the `Complex Characteristics' section of\\
{\tt functions.wolfram.com}.
From this one finds
\begin{equation}
|K(\exp(i\phi)|^2 
=\frac{1}{4}\left(K(\cos(\phi))^2 + K(\sin(\phi))^2\right) . 
\label{eq:Kdisk}
\end{equation}
See equation~(\ref{eq:weightK}).

\section{Appendix C: Result Validation and Numerical Examples}\label{sec:Numerical}
\subsection{Slip flow in Elliptic Channel with $a=2,\;b=\frac{1}{2}$}\label{subsec:Spiga}
\begin{figure}[hb]

\begin{tabular}{cc}
  
\includegraphics[height=2cm,width=7cm]{\grpath/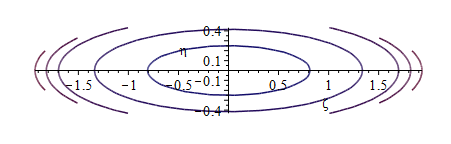} & \includegraphics[height=2cm,width=7cm]{\grpath/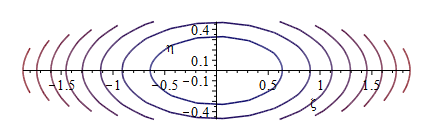} \\

  $\beta=\frac{1}{4}$         &  $\beta=4$  \\

  \includegraphics[height=2cm,width=7cm]{\grpath/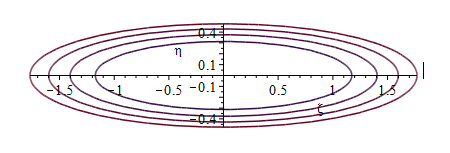}  & \includegraphics[height=2cm,width=7cm]{\grpath/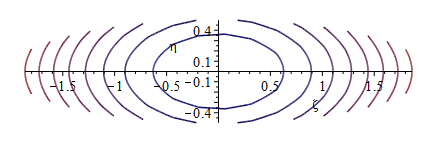} \\

  $\beta=\frac{1}{16}$         &  $\beta=16$  \\

  \includegraphics[height=2cm,width=7cm]{\grpath/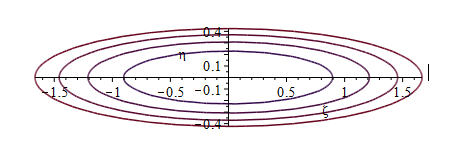} &\includegraphics[height=2cm,width=7cm]{\grpath/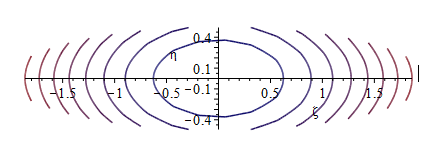} \\

  $\beta=\frac{1}{64}$         &  $\beta=64$  \\

\end{tabular}

\caption{$a=2$. Contour plot of axial velocity on the ellipse for six different $\beta$ values.}
\label{fig:ContourPlotBetajpg}
\end{figure}

Figure~\ref{fig:ContourPlotBetajpg} gives, at left,
contour plots of $u$ for $\beta=\frac{1}{4},\;\frac{1}{16}$ and $\frac{1}{64}$, and, at right,
those for $\beta=4,~16,~64$.  
The contours are very similar to those in Figure 4 of~\cite{SV12}.
It is noted that when $\beta>0$ on the boundary
(as well as in the interior) $u>0$. 
At any given fixed position $z$ on the boundary, $u(z)$ increases when $\beta$ value increases.

\clearpage
\subsection{Comparing the results with Ritz Method for Slip Flow}\label{subsec:Wang}

Here we compare with results given in~\cite{Wa14}.
\medskip

\noindent
\begin{tabular}{|c|r|r|r|r|r|r|r|r|}
  \hline
 $\lambda$&     $c$     &    $a=\frac{1}{\sqrt{c}}$  & $\beta=\frac{\lambda}{\sqrt{c}}$ & A & V & F &  Ritz method  \\
  \hline
          &    .25      &   2.000000    &      .200000      &   .027080     &   .026893     &   .0268994  & .0268994 \\
    0.1   &     .5      &   1.414213    &      .141421      &   .131628     &   .131224     &   .1312286  & .1312286  \\
          &    .75      &   1.154700    &      .115470      &   .313497     &   .313320     &   .3133204  & .3133204   \\
  \hline
          &   .25      &   2.000000     &      .400000      &   .042611     &   .041996     &   .0419990  & .0419990  \\
   0.2    &    .5      &   1.414213     &      .282843      &   .184716     &   .183372     &   .1833741  & .1833741  \\
          &   .75      &   1.154700     &      .230940      &   .414937     &   .414338     &   .4143388  & .4143388   \\
  \hline
          &   .25      &   2.000000     &      1.000000     &   .071907     &   .086515     &   .0865169  &  .0865167  \\
  0.5     &   .5       &   1.414214     &      .707107      &   .343981     &   .338284     &   .3382847  &  .3382847  \\
          &   .75      &   1.154700     &      .577350      &   .719257     &   .716699     &   .7166986  &  .7166986  \\
  \hline
          &   .25      &   2.000000     &       2.000000    &    .143814    &   .159568     &   .1595802  &  .1595797  \\
   1.0    &    .5      &   1.414214     &      1.414214     &    .509349    &   .594535     &   .5945397  &  .5945397  \\
          &   .75      &   1.154701     &      1.154701     &   1.004665    &  1.219750     &  1.2197505  & 1.2197505  \\
  \hline
          &  .25      &   2.000000      &      4.000000     &     .287629   &   .304352     &   .3043838  &   .3043838  \\
   2.0    &   .5      &   1.414214      &      2.828427     &    1.018698   &  1.105106     &  1.1051201  &  1.1051201  \\
          &  .75      &   1.154700      &      2.309401     &    2.009330   &  2.224974     &  2.2249748  &  2.2249748  \\
  \hline
          &  .25      &   2.000000      &      10.000000    &    .719072    &   .736625     &   .7366869  &   .7366851 \\
  5.0     &  .5       &   1.414214      &       7.071068    &   2.546745    &  2.634130     &  2.6341541  &  2.6341541 \\
          &  .75      &   1.154701      &       5.773503    &   5.023326    &  5.239416     &  5.2394175  &  5.2394175 \\
  \hline
          & .25       &   2.000000      &      20.000000    &   1.438144    &  1.456036     &  1.4561139  &  1.4561106 \\
   10     & .5        &   1.414214      &      14.142136    &    5.093491   &  5.181257     &  5.1812861  &  5.1812861 \\
          & .75       &   1.154701      &      11.547005    &  10.046652    & 10.262916     & 10.2629186  & 10.2629186 \\
  \hline
\end{tabular}

\subsection{Blood flow problem}\label{subsec:Numerical}
\noindent Most of the veins and arteries in our bodies can be taken as having circular cross-section.
However in places where the vein or artery has to go through a region which is squeezed in by muscle or bone
one might expect the cross-section of the vein or artery to depart somewhat from circular.
In this appendix, we take  fluid density,$\rho$, of 1.05 $g/mL$ and fluid viscosity, $\mu$, of 0.04 $Poise$. See Table~2 in~\cite{Wblood}.
The geometry of elliptical cross-section of the channel is described by coordinates $(x, y)$ with $$x=\overline{r}(\psi)
\cos(\psi);\;\;\;  y=\overline{r}(\psi) \sin(\psi),$$
where $$\overline{r}(\psi)=\lambda k(1+\kappa(k, \psi)),\;\;\;0 \leq \psi \leq 2\pi$$  and
$$\kappa = -\frac{1}{4}(k-1)^2+(k-1-\frac{1}{2}(k-1)^2)\cos(2\psi)+\frac{3}{4}(k-1)^2\cos(k \psi).$$
The tube ellipticity $\varepsilon=\sqrt{1-b^2/a^2}$ with the lengths of the half-axes $a=\bar{r}(0)\cos(0)$ and
$b=\bar{r}(\pi/2)\sin(\pi/2)$ is determined by setting $\lambda=0.005$ and $k=1.15$. Thus, the tube has the ellipticity of 0.6720
with $a=6.6\; \mu m$, $b=4.9\; \mu m$, $\eta_0\;=\;0.9518$ and $c=4.4 \;\mu m$. The boundary of the elliptical cross-section is described by
$$\partial \Omega: (x,y)=(c \cosh \eta_0 \cos \psi, c \sinh \eta_0 \sin \psi).$$
\begin{figure}
\centering
\includegraphics[scale=0.65]{\grpath/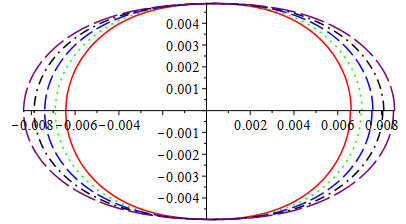} \\
\caption{The cross sections of elliptical shapes with the same length of minor axis $b = 4.9 \;\mu m.$ and
various lengths of major axis $a=6.6,\; 7.1,\; 7.6,\; 8.1$ and $8.5 \;\mu m.$ }
\label{fig2}
\end{figure}
\begin{figure}
\centering
\begin{tabular}{c}
\includegraphics[scale=0.55]{\grpath/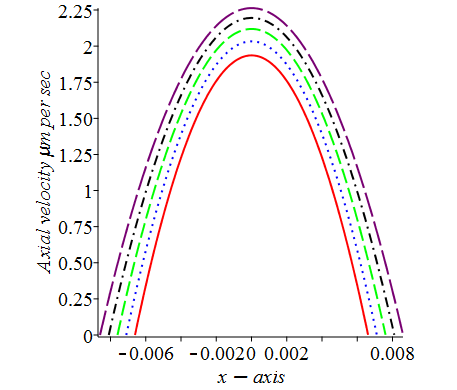} \includegraphics[scale=0.55]{\grpath/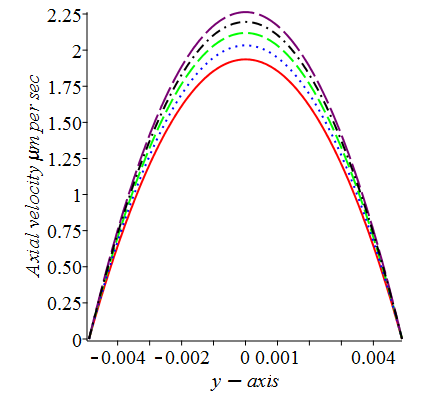} \\
\end{tabular}
\caption{Axial velocity profiles in $x$ and $y$ directions obtained from the model with
no slip lengths: $\ell =0.0$, and various shapes of cross section.}
\label{fig3}
\end{figure}
%
\begin{figure}
\centering
\includegraphics[scale=0.65]{\grpath/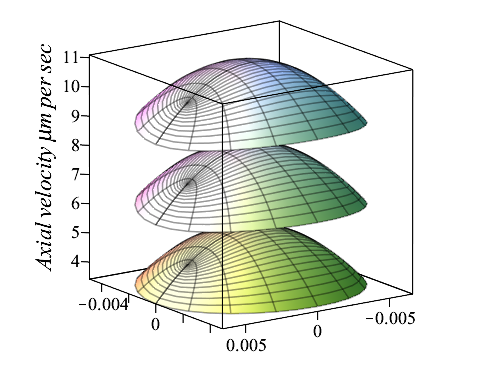} \\
\caption{Axial velocity profiles obtained from the model with $a=6.6 \;\mu m.$ and  $b := 4.9 \;\mu m.$}
\label{fig4}
\end{figure}
\begin{figure}
\centering
\begin{tabular}{c}
\includegraphics[scale=0.55]{\grpath/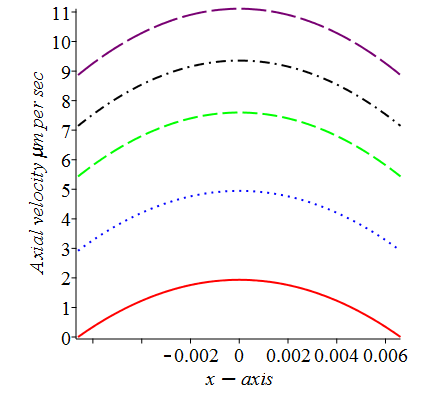} \includegraphics[scale=0.55]{\grpath/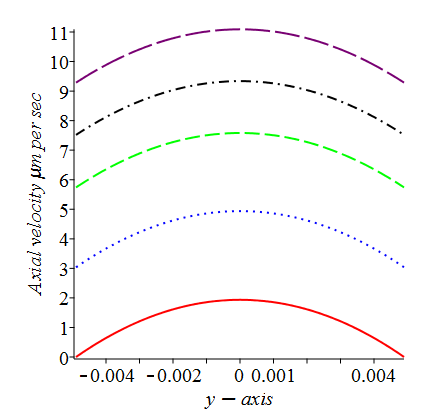} \\
\end{tabular}
\caption{Axial velocity profiles in $x$ and $y$ directions obtained from the model having $a=6.5 \;\mu m.$
 and  $b = 4.9\;\mu m$ and various slip lengths $\ell (\mu m.)$ of 0.0 (solid line), 6.25(dotted line),
 25 (dashed line), 37.5  (dash-dotted line)  and 50 (long dashed line).}
\label{fig5}
\end{figure}
\noindent The investigate the influence of elliptical shapes on constant pressure-driven flow of fluid, various sizes of the major axis $a$ are chosen to be vary between  $6.6 \;\mu m.$ to $8.6 \;\mu m.$  while the length of minor axis is fixed as
 $b = 4.9 \;\mu m$.\\

\noindent To demonstrate the impact of the slip length $\ell$ on the constant pressure-driven flow of fluid, values of slip length
 $\ell$  are chosen to be vary from $0.0$ to $50 \;\mu m.$.
The results show that the slip length has a direct influence on the axial velocity.
Larger slip length provides higher velocity as shown in Figures \ref{fig3}, \ref{fig4} and \ref{fig5}.

\clearpage

\clearpage

\newpage

\begin{center}
{\Large{\textbf{\textsc{ Part II:
$u_\infty$ by Fourier series\\
  \qquad\\
By G. Keady
}}}}
\end{center}
\section{More about $u_\infty$}

In Part I we only noted the dominant approximation to $u$ when
$\slipparameter$ tends to infinity.

\subsection{The next term, and $u_\infty$, for ellipses in general}\label{subsec:Qinf}

We seek the term $\Sigma_\infty$ (independent of $\beta$) in the asymptotic expansion
$$ Q \sim \frac{\beta\, |\Omega|^2}{|\partial\Omega|} + \Sigma_\infty
\qquad {\rm for\ \ } \beta\rightarrow\infty . $$
The notation here is, as in~\cite{KM93},
$$ \Sigma_\infty
= \int_\Omega u_\infty ,
$$
with $u_\infty$ satisfying equations~(\ref{eq:uInf}).


The approach here is to seek a Fourier series solution for $u_\infty$,
investigate its contour plot, and after this
to integrate $u_\infty$ thereby finding
the next term in the expansion for $Q$.
The integration task only requires the constant ($n=0$) and the next ($n=1$) Fourier coefficients of $u_\infty$.
The tasks have been completed  for nearly circular ellipses, but the general case,
specifically determining the $n=0$ constant coefficient of $u_\infty$, remains `work in progress'.

We recast the elliptic coordinates Fourier series approach so that one uses the Fourier series of $1/g$ rather than $g$.

\medskip
{\noindent\bf Fourier series of ${\hat g}=1/g$}

Write the Fourier series of $1/g$ as
\begin{equation}
 \frac{1}{g(s)} = \frac{{\hat g}_0}{2} + \sum_{n=1}^\infty {\hat g}_n \cos(2 n s) .
\label{eq:rgFSB}
\end{equation}
As with the earlier calculation of the Fourier series of $g$, the coefficients ${\hat g}_n$ can be found in terms of
EllipticE and EllipticK functions.
In particular, with
the notation of~(\ref{eq:Ellnotn}),
we have
\begin{eqnarray}
|\partial\Omega|
&=& 4 c\, \cosh(\eta_0)\, {\rm EllipticE}_0 ,
\nonumber\\
{\hat g}_0
&=&  \frac{4}{\pi} \, {\rm EllipticE}_0
=  \frac{1}{\pi\,  c\, \cosh(\eta_0)} \, |\partial\Omega|
\label{eq:gBar0Ell}\\
{\hat g}_1
&=&  \frac{4}{3\pi} \, \left( -\cosh(2\eta_0) {\rm EllipticE}_0 +
 (\cosh(2\eta_0) -1) {\rm EllipticK}_0 \right) ,
\label{eq:gBar1Ell}
\end{eqnarray}
and all the higher ${\hat g}_n$ can also be expressed similarly in terms of ${\rm EllipticE}_0$ and ${\rm EllipticK}_0$:
\begin{equation}
{\hat g}_n =\frac{4}{\pi}\left(  {\hat E}_n  \, {\rm EllipticE}_0 +
{\hat K}_n\, {\rm EllipticK}_0 \right), 
\label{eq:gBarnEK}
\end{equation}
where ${\hat E}_n$ and ${\hat K}_n$ are polynomials of degree $n$ in 
$q=\cosh(\eta_0)^2=1/e^2$
with rational number coefficients.

Again there is  a three term recurrence relation exists to
 determine the polynomials ${\hat E}_n$ and ${\hat K}_n$.
 In the notation of equation~(\ref{eq:qEta0}), the ${\hat K}_n$ sequence of polynomials starts with
 $$ {\hat K}_0=0 , \qquad {\hat K}_1= \frac{2}{3} q =\frac{1}{3} ( q_2 -1).
 $$
 The ${\hat E}_n$ sequence of polynomials starts with
 $$ {\hat E}_0=1 , \qquad {\hat E}_1= -\frac{1}{3} (2 q-1)=-\frac{1}{3}  q_2 .
$$
See \S~\ref{app:Furtherg}, and the recurrence~(\ref{eq:gBarnrec}).

\medskip
{\noindent\bf Calculating $u_\infty$ (at least up to an additive constant}
\smallskip

As in~\cite{KM93}, $u\sim~\beta|\Omega|/|\partial\Omega|+u_\infty$ where
\begin{equation}
-\Delta u_\infty = 1\ {\rm and\ } \frac{\partial u_\infty}{\partial n}= -\frac{|\Omega|}{|\partial\Omega|}, \ \
\int_{\partial\Omega} u_\infty = 0 . 
\label{eq:uInf}
\end{equation}
Set 
$$ u_\infty = v_\infty +  v_p , $$
with $v_p$ as in (\ref{eq:12b})  before.
Once again the harmonic function $v_\infty$ can be represented as a Fourier series
\begin{equation}
 v_\infty = \sum_{n=0}^\infty V_n \cosh(2n\eta)\cos(2n\psi) .
 \label{eq:vInfDef}
\end{equation}
To determine $V_n$ for $n\ge{1}$ we will need
$$ \frac{\partial v_\infty}{\partial \eta}(\eta_0,\psi)
=  \sum_{n=1}^\infty 2n\,V_n \sinh(2n\eta_0)\cos(2n\psi) \qquad{\rm on\ } \eta=\eta_0.
$$
The Neumann boundary condition is
$$ \frac{g(\psi)}{c\cosh(\eta_0)} \frac{\partial u_\infty}{\partial \eta}(\eta_0,\psi)
=  \frac{1}{c\sqrt{\cosh^2\eta_0 - \cos^2 \psi}} \frac{\partial u_\infty}{\partial \eta}(\eta_0,\psi)
\;=\;  -\frac{|\Omega|}{|\partial\Omega|} .
$$
In terms of the harmonic function $v_\infty$ this is
\begin{equation}
\frac{\partial v_\infty}{\partial \eta}
= -\frac{\partial v_p}{\partial \eta} - \frac{|\Omega|}{|\partial\Omega|}\, \frac{c\cosh(\eta_0)}{g(\psi)}
\qquad{\rm on\ } \eta=\eta_0.
\label{eq:vInfNormal}
\end{equation}
Now one finds
\begin{equation}
\frac{\partial v_p}{\partial \eta}=  -\frac{1}{2} \qquad{\rm on\ } \eta=\eta_0.
\label{eq:vpeta}
\end{equation}
(In Cartesian coordinates we have
$\nabla{v_p}\cdot{\nabla{u_0}}$ is constant around $\partial\Omega$,
but we seek $u_\infty$ such that
$\nabla{u_\infty}\cdot{\nabla{u_0}}$ is a constant
multiple of $|{\nabla{u_0}}|$ there.)
The integral around the boundary of the normal derivative of a harmonic function must be 0,
and this accords with the value of ${\hat g}_0$,  the constant term in the Fourier series of $1/g$ to ensure this as
combining~(\ref{eq:vpeta}) and
$$  \frac{|\Omega|}{|\partial\Omega|}\, {c\cosh(\eta_0)}{{\hat g}_0}
= \frac{\pi}{|\partial\Omega|} \,   \frac{1}{2\pi} \, |\partial\Omega| =\frac{1}{2}\\
$$
in~(\ref{eq:vInfNormal}) we have the required result.

From the Fourier series of $1/g(\psi)$ we can find all  the Fourier coefficients~(\ref{eq:vInfDef}),
except the constant term, in the series  for $v_\infty$.
In particular the term $V_1$ is found from~(\ref{eq:vInfNormal}) (on using~(\ref{eq:vpeta}), contributing 0)
$$ 2 V_1 \sinh(2\eta_0)
 = 0 - \frac{|\Omega|}{|\partial\Omega|}\, {c\cosh(\eta_0)}{{\hat g}_1}
 = -\frac{\pi\, {\hat g}_1}{4{\rm EllipticE}_0} .
$$
Hence
\begin{equation}
 V_1= - \frac{\pi\, {\hat g}_1}{8{\rm EllipticE}_0\, \sinh(2\eta_0)} .
 \label{eq:V1inf}
\end{equation}
Similarly
\begin{equation}
 V_2= - \frac{\pi\, {\hat g}_2}{16\,{\rm EllipticE}_0\, \sinh(4\eta_0)} ,
 \label{eq:V2inf}
\end{equation}
and, more generally, for $n>0$,
\begin{equation}
 V_n= - \frac{\pi\, {\hat g}_n}
          {8\, n\,{\rm EllipticE}_0\, \sinh(2n\eta_0)} .
 \label{eq:Vninf}
\end{equation}

\medskip
\par\noindent{\bf Contour plots for $u_\infty$}

As we have all the Fourier coefficients $V_n$ except $V_0$ we 
can produce contour plots.
We show one of these, that for $a=2$.

\begin{figure}[!h]
\centerline{\includegraphics[height=6cm,width=10cm]{\grpath/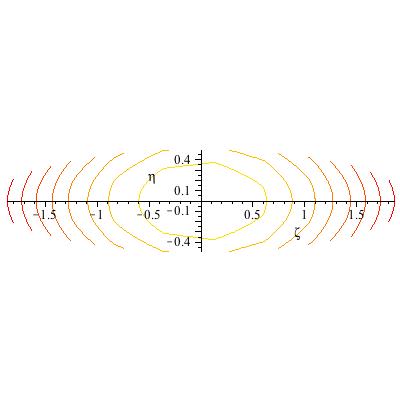}}
\caption{$a=2$. Contours of $u_\infty$}
\label{fig:contoursuInfaEq2jpg}
\end{figure}

The contourplot (at $a=2$) is visually indistinguishable from that produced
from the quadratic found using the $c_0$, $c_XX$, $c_YY$ values found
in~\S\ref{subsec:VarlEll}.
Except when $\eta_0$ is small, the $V_n$ will decrease rapidly with
$n$, and the plots are essentially those from truncating the
Fourier series so that it only includes the $V_1$ term.
The numerical solutions given in~Figure~\ref{fig:ContourPlotBetajpg} 
of Part I Appendix~C 
are, for $\beta\ge{4}$, very much like that
of Figure~\ref{fig:contoursuInfaEq2jpg}.

\clearpage
\medskip
\par\noindent{\bf Finding $V_0$ and $Q$?}

The calculation of $V_0$ remains.
We will outline how it can be found as an infinite series,
but shall only complete the details for nearly circular ellipses.

The condition that the boundary integral of $u_\infty$ is zero 
remains to be used.
 Now, with, as usual, $s$ for arc length around the boundary, $\eta=\eta_0$,
 $$(\frac{d s}{d\psi})^2 ={(\frac{d x}{d\psi})^2+(\frac{d y}{d\psi})^2}
 = c^2\left(\cosh(\eta_0)^2-\cos(\psi)^2\right)
 = J(\eta_0,\psi) .
$$
The square root of the Jacobian is related to the function $g$:
 $$ \sqrt{J(\eta_0,\psi)}= \frac{c\cosh(\eta_0)}{g(\psi)} = \frac{a}{g(\psi)}. $$
We need to satisfy
$$ 0 = \int_{\partial\Omega} u_\infty
 = 4  \int_{0}^{\pi/2} (v_p+v_\infty) \, \sqrt{J(\eta_0,\psi)}\, d\psi . $$
The integral $I_{p,\partial}$ around the boundary of $v_p$ is  
$-i_2/4$
where $i_2$ is the boundary moment of inertia and can also be evaluated directly, 
\begin{eqnarray*}
I_{p,\partial}
  &=&  4  \int_{0}^{\pi/2} v_p \, \sqrt{J(\eta_0,\psi)}\, d\psi 
= 4a  \int_{0}^{\pi/2} v_p \,{\hat g}(\psi)\, d\psi ,\\
&=& - \frac{a^3}{3}\left(  (2-e^2) {\rm EllipticE}_0 + 
 (1-e^2) {\rm EllipticK}_0 \right).
\end{eqnarray*}
\cmdvmcode{IpBndryeSmallmpl.txt}{ }

When evaluating the boundary integral of $v_\infty$, the square~root brings in all the $V_n$.
\begin{itemize}
\item
It may be that there are useful Mean Value results applying to ellipses, and harmonic functions in them, of
 $$\int_{\partial\Omega} v_\infty - v_\infty(0)\, |\partial\Omega|\ ? $$
 It is, of course, zero for a circle. In general it wou't be zero, as, for example,
$$ \int_{\partial\Omega}  (x^2-y^2) = i_X - i_Y \ne{0} {\rm \ \ when\ \ } a>1 . $$
Some boundary integrals can be related to integrals along the 
major axis of the ellipse, and~\cite{Ro64} is given as 
a reference in papers by Symeonidis.
For any function $h$ harmonic in the ellipse, with $c$ as before
$c=\sqrt{a^2-1/a^2}$,
$$ \int_{-c}^c \frac{h(x,0)}{\sqrt{c^2-x^2}}\, dx
= \frac{1}{2}\int_{-\pi}^\pi h(a\cos(\psi),\frac{\sin(\psi)}{a}) \, d\psi .
$$
(The right-hand side is not the integral of $h$ with respect to arc length.)
\item The integral
$$
I_{\infty,\partial}
=  4  \int_{0}^{\pi/2} v_\infty \, \sqrt{J(\eta_0,\psi)}\, d\psi 
= 4a  \int_{0}^{\pi/2} v_\infty \,{\hat g}(\psi)\, d\psi ,
$$
can be evaluated, with $V_0$ left as an unknown using
Parseval's Theorem. 
We could find $V_0$ if we could evaluate
$$ \sum_{n=1}^\infty \frac{{\hat g}_n^2}{n\tanh(2n\eta_0)} ,
$$
which can be recast as a series involving Legendre functions
$Q_{n-1/2}^{-1}$: see~\S\ref{app:Furtherg}.
At each fixed $a>1$ (or equivalently $\eta_0$ or 
eccentricity $e$) this would be a simple float numerical task.
However, we have not yet found
a simple formula for $V_0$,
and hence $\Sigma_\infty$, as simple as that
we have for $Q_1$, in equation~(\ref{eq:Q1bSmall}).

\end{itemize}

\cmdvmcode{checksEllCoordsmpl.txt }{ }

Finally to determine $\Sigma_\infty$ using $u_\infty= v_p + v_\infty$ and recalling that we will only need, as in~(\ref{eq:QFS}), the
$n=0$ and $n=1$ Fourier coefficients of $v_\infty$:
\begin{equation}
\Sigma_\infty
= - \frac{\pi}{8} \, \frac{\cosh(2\eta_0)}{\sinh(2\eta_0)}+ \pi\, V_0 - \frac{1}{2}\,\pi\, V_1 .
\label{eq:QinfFS}
\end{equation}
Combining this with the dominant term we have
$$ Q\sim \frac{\beta |\Omega|^2}{|\partial\Omega|} + \Sigma_\infty \qquad{\rm as\ } \beta\rightarrow\infty .$$

We remark that for $a$ large $\Sigma_1$ becomes large.
The asymptotics for $Q(a,\beta)$ for $\beta$ large
should only include the $\Sigma_1$ term when 
$\Sigma_1/\beta$ is smaller than $\Sigma_\infty$.
In practice this means that, for $a$ large,
the large $\beta$ asymptotics are good only for
$\beta$ very large.
(When the domain is a rectangle, rather than an ellipse,
there is the same behaviour and this is very easy to
demonstrate as $u_\infty$ is, in Cartesian coordinates,
a simple quadratic.)

\subsubsection*{$\beta$ large  for nearly circular ellipses}

Recall that the eccentricity $e=1/\cosh(\eta_0)$.
From equation~(\ref{eq:gRecipFSnearCirc})
$$ {\hat g}_0 = 1-\frac{e^2}{4},\qquad {\hat g}_1 =  -\frac{e^2}{4} . $$
Also
$$ {\rm EllipticE}_0 \sim
\frac{\pi}{2}  - \frac{\pi}{8}\, e^2 - \frac{3\pi}{128}\, e^4 ,
\qquad
{\rm EllipticK}_0\sim \frac{\pi}{2} + \frac{\pi}{8}\, e^2 + \frac{9\pi}{128}\, e^4 ,
 $$
and
$$ a=\left(\frac{1}{1-e^2}\right)^{1/4}
\sim 1+ \frac{1}{4}\, e^2+ \frac{5}{32}\, e^4  ,\quad
c=a\,e ,$$
and
$$|\partial\Omega|=4a\, {\rm EllipticE}_0 \sim
2\pi + \frac{3\pi}{32} e^4 .
$$
(The terms in $e^4$ in the series for ${\rm EllipticE}_0$ and $a$ are only there to check against the series for $|\partial\Omega|$.)
We also note
$$ \sinh(2\eta_0)=2 \frac{1}{e}\,\sqrt{\frac{1}{e^2}-1}
\sim \frac{2}{e^2} -1 - \frac{1}{4}\, e^2 .
$$

The Fourier coefficients of $v_\infty$ are needed.
From equation~(\ref{eq:V1inf})
$$V_1=  -\frac{\pi\, {\hat g}_1}{8{\rm EllipticE}_0\, \sinh(2\eta_0)}
\sim \frac{1}{32} e^4 .
$$
\cmdvmcode{V1eSmallmpl.txt}{ }

\medskip
\noindent{\bf Comparison with earlier near-circular results}
\smallskip

Because, for nearly circular ellipses,
the asymptotics for all the $V_n$, $n\geq{1}$
are uncomplicated, and the terms get rapidly smaller in $n$,
we can find $V_0$.

\cmdvmcode{NearCirc/nearCircApproxQfnEpsmpl.txt}{  }

Maple code gives, for $\beta\rightarrow\infty$ and $e\rightarrow{0}$,
$$ Q\sim
 \left( \frac{1}{2}\,\pi -{\frac {3}{128}}\,\pi \,{e}^{4} \right) \beta+\frac{1}{8}\,
\pi -{\frac {1}{128}}\,\pi \,{e}^{4}-{\frac {1}{256}}\,{\frac {\pi \,{ e}^{4}}{\beta}}
. $$

\ifthenelse{\boolean{vmcode}}{
Of course the dominant term in $\beta$ is readily checked.}{ }
\cmdvmcode{NearCirc/nearCircBetaLarge0mpl.txt}{ }

Another quick look at the $\epsilon$ small solution
of \S\ref{subsec:nearCircbetaGt0}, in $\eta$, $\psi$ coordinates,
produces the same result.

\cmdvmcode{NearCirc/nearCircEtaPsimpl.txt}{ }

\cmdvmcode{uInfFS/findV0eSmallmpl.txt}{ }

\subsection{$u_\infty$ for other domains}

It seems rarly to happen that both $u_0$ and $u_\infty$
have simple explicit formulae representing them.
An exception is the equilateral triangle: see~\cite{KM93,MK94}.
And, of course, the circular disk has $u_0=u_\infty$
and is the only domain for which there is equality.
(There may well be other less trivial domains, for example,
the semicircle.)

For our ellipse, $u_0$ is simply, in Cartesians, 
a quadratic polynomial, or equivalently has a Fourier
series (in elliptic coordinates) with just 2 terms, while
$u_\infty$ has a Fourier series with all terms nonzero.

The situation with a rectangle is as follows.
There is an elaborate Fourier series in Cartesian coordinates
for $u_0$.
However $u_\infty$ is simply a quadratic in the Cartesian coordinates.
Indeed our function $u_0$, given in equation~(\ref{eq:u0Cart}),
is $u_\infty$ for the rectangle
$(-a^2,a^2)\times(-b^2/b^2)$.
The methods in~\S\ref{subsec:VarlGen} would give good
lower bounds for $Q({\rm rectangle})$ when $\beta$ is large.
See Part III~\S~\ref{subsec:VarlRectq}.
(For large $\beta$ asymptotics for a rectangle, see~\cite{KWK18}.)

\section{A difference equation -- and toroidal functions}
\label{app:Furtherg}

The Fourier coefficients $g_n$ and ${\hat g}_n$ are central
to our separation of variables solution of
Problems (P($\beta$)) and (P($\infty$)) respectively.
On using maple to find these, maple returned (instances of) them
as the sum of two terms: the first term is the product of
a polynomial with an EllipticE function, and
the second term is the product of
a polynomial with an EllipticK function.
The dependence on $n$ is only in the polynomials.
See equations~(\ref{eq:gnEK}), (\ref{eq:gBarnEK}),
repeated at equations~(\ref{eq:gnEll}), (\ref{eq:gBarnEll}).
The polynomials satisfy the same recurrence relations as
$g_n$ and ${\hat g}_n$: see (RE($\alpha$)) below
where $\alpha=1/2$ for $g_n$ and $\alpha=3/2$ for ${\hat g}_n$.

It happens that $g_n$ and ${\hat g}_n$ can also be expressed 
in terms of LegendreQ functions -- specifically toroidal functions.
Again sequences of polynomials satisfying (RE(1/2)) and (RE(3/2))
arise: see equations~(\ref{eq:Q1pol}), (\ref{eq:Q3pol}).
More generally, for any $\alpha$, these polynomial sequences
are present in the entries on Legendre functions at functions.wolfram.com.
There is also some literature concerning the polynomials
in families of recurrence relations which include
our case $\alpha=1/2$:
see~\cite{BD67,Gr85} 
and subsection~\ref{subapp:orthog} below.
It isn't obvious what, for general $\alpha$, 
the polynomials should be called,
though for $\alpha=0$ and $\alpha=1$ they are
Chebyshev polynomials.
In the case $\alpha=1/2$, the sequence of polynomials
(particularly those labelled $p_{01}$)
might be called `associated toroidal Legendre polynomials',
and we now refer to all
the sequences of polynomials $p_{01}(n)$, $p_{10}(n)$ of equation~(\ref{eq:Q1pol})) and
${\hat p}_{01}(n)$, ${\hat p}_{10}(n)$ of equation~(\ref{eq:Q3pol}))
(or linear combinations thereof, e.g. $p_{10}+x\,p_{01}$
and, more generally, $\alpha$ half an odd integer)
by this name.
Though others have used the adjective `associated' in this context,
we remark that the adjective is used in several different ways.
While the previously cited~\cite{BD67,Gr85} on `associated Legendre functions' includes our case of $\alpha=1/2$ it doesn't include
$\alpha=3/2$.
`Associated Gegenbauer polynomials' as in~\cite{BI82}\S3
or, more generally, the family of Pollaczek polynomials
introduced in~\cite{Poll50} (see~\cite{HTF2}p220) cover general $\alpha$ but have
more parameters than our simpler cases.
(See also our \S\ref{subapp:orthog}.)

Questions which arise include the following.
\begin{itemize}
\item Are there simple formulae for these sequences of polynomials?
E.g. (i) Can each sequence be expressed as a single hypergeometric function (perhaps a ${}_{3}F_2$)?
(ii) What are the formulae for coefficients?

\item Finally, how might the information on the polynomials be best
applied to progress the tasks associated with the pde questions
which we have listed below in a subsection concerning `Hopes'.

\end{itemize}

A related problem, also involving sequences of polynomials
satisfying (RE(1/2)) or (RE(3/2)) is to discover formulae,
for general $n$, agreeing with maple's output for the
$g_n$ and ${\hat g}_n$ involving EllipticE and EllipticK functions.
This also seems to be an old problem:
finding formulae for the toroidal functions $Q_{n-1/2}$
involving EllipticE and EllipticK functions.
We will do relatively little on this last problem.


\subsection{Notation}\label{subapp:Notation}

We denote the reciprocal of $g$ by $\hat g$:
$${\hat g}(\psi)=\frac{1}{g(\psi)}=\sqrt{1-\frac{\cos^2(\psi)}{q}} .$$
As described earlier, the Fourier cosine coefficients of $g$ are denoted $g_n$ (see equation~(\ref{eq:bgnDef})),
and those of ${\hat g}$ by ${\hat g}_n$:
\begin{equation}
g_n = \frac{2}{\pi}\int_0^\pi g(\psi)\cos(2 n\psi)\, d\psi ,
\label{eq:agnDef}
\end{equation}
\begin{equation}
{\hat g}_n = \frac{2}{\pi}\int_0^\pi {\hat g}(\psi)\cos(2 n\psi)\, d\psi .
\label{eq:bgBarnDef}
\end{equation}
Both $g_n$ and ${\hat g}_n$ can be written concisely in terms of
toroidal functions, Legendre functions:
see equations~(\ref{eq:gnQ}),~(\ref{eq:gBarnQ}).
The sequences $g_n$ and ${\hat g}_n$ are used, so far, 
just in numerics.
In both sequences our interest is in $q>1$ but, the integral for ${\hat g}_n$ evaluates at
other values of $q$, for example at $q=1$ to
\begin{equation}
{\hat g}_n(q=1) 
=\frac{2}{\pi} \int_0^\pi \sin(\psi) \cos(2n\psi) \, d\psi
= \frac{1}{\pi}\, \frac{1}{n^2-1/4} . 
\label{eq:gBarnq1}
\end{equation}
Though the integrals are complex, one finds
$$ \frac{g_{n+1}(q=1/2)}{g_{n-1}(q=1/2)}
= -\, \frac{n-1/2}{n+1/2} ,$$
and a similar formula for the ${\hat g}_n$ sequence.

\subsection{Hopes for application to the pde problems}\label{subapp:Hopes}

We have various hopes for future results.
\begin{itemize}
\item
One hope is that further results on ${\hat g}_n$ may ultimately lead to further analytical results
concerning $u_\infty$ and, in particular, its values on the boundary, and thence
to analytical expressions for $\Sigma_\infty$ and $\Sigma_1$.

\item
A more modest goal, also as yet not implemented, is that it may be possible to improve the efficiency
of the numerical calculations of $g_n$ and of ${\hat g}_n$
as the present codes call on maple to evaluate each of the integrals separately for each $n$
and uses maple's results given in terms of EllipticE and EllipticK functions.

\item
There also remains the prospect, should there be any demand for codes in languages
used by engineers, e.g. matlab, to provide such codes.
(We remark that as at 2019a matlab's legendre function
does not compute toroidal functions.
We have yet to check the matlab File Exchange site.)
We remark that there are numerical stability issues with recurrence
relations when performed purely within fixed precision float computation
but these can be circumvented using the Symbolic Toolbox.
It isn't essential for the code to be entirely within matlab 
as it would also be possible to use existing Fortran or C codes
(such as described in~\cite{GST00}) mex-ed into matlab.
(Matlab codes and related algorithms for toroidal functions
are mentioned by several authors, Majic is one,
Gil and Segura others, e.g~\cite{GST00}.
There is a library, POLPAK, implementing algorithms in various languages,
including matlab.)\\
Noting the present limitations of matlab's legendre function, 
our present suggestion is to compute $Q_{-1/2}$, $Q_{1/2}$
and $Q_{-1/2}^{-1}$, $Q_{1/2}^{-1}$ from their expressions in
terms of EllipticE and EllipticK functions,
and then the toroidal $Q$ functions using the polynomials
$p_{01}$, $p_{10}$ and ${\hat p}_{01}$, ${\hat p}_{10}$
defined below. 
(Matlab's Symbolic Toolbox might have a role in connection
with the polynomials.)

\end{itemize}

\subsection{Introducing the difference equation (RE($\alpha$))}\label{subapp:REalpha}

In the body of this paper we noted, at equations~(\ref{eq:gnEK}) and
(\ref{eq:gBarnEK}), that maple,
when asked to perform the integrals defining $g_n$ and ${\hat g}_n$
returned a result in which the terms multiplying elliptic integrals
are polynomials in $q$ or rather $q_2=2q-1$.
We denoted the polynomials by $E_n(q_2)$, $K_n(q_2)$
(these satisfying the recurrence (RE(1/2)) defined below)
and by ${\hat E}_n(q_2)$, ${\hat K}_n(q_2)$
(satisfying the recurrence (RE(3/2))).

While we have yet to make any significant use of the items here for
our Robin boundary condition problems,  we have discovered that the Fourier series 
of $g$ and of $\hat g$,
and related studies of toroidal functions, have been studied in connection with various different problems since at least the 1880s.
In spite of the bewilderingly large literature both
on Legendre functions and on orthogonal polynomials, 
we have not seen a systematic account of the few results
we present
concerning sequences of polynomials satisfying the simple recurrences satisfied by
each of $g_n$ and of ${\hat g}_n$.

Our sequences of Fourier coefficients satisfy difference equations of the form
$$ (n+\alpha) u_{n+1} =2 q_2 n u_n - (n-\alpha) u_{n-1} ,
\eqno{({\rm RE}(\alpha))}$$
where $q_2=2q-1$, and where conformity to conventions in the 
orthogonal polynomial literature seems appropriate, we
write $x$ rather than $q_2$.
The $g_n$ satisfy this with $\alpha=1/2$ while the ${\hat g}_n$ satisfy it with $\alpha=3/2$.

The recurrence relation (RE($\alpha$)) satisfied by Legendre functions is given in standard references such as
\cite{AS}~8.5.3 and  \cite{HTF1}~\S3.8 (2).
The order $\mu$ of the Legendre functions,
$Q^\mu_\nu(q_2)$ is determined by
$\alpha=\frac{1}{2}-\mu$:
the degrees in the sequence are of the form $\nu=n-1/2$.

Our interest is in $\alpha=1/2$ 
(for which we give later polynomial sequences
$p_{01}(n)$ and $p_{10}(n)$)
and $\alpha=3/2$ 
(for which we give later polynomial sequences
${\hat p}_{01}(n)$ and ${\hat p}_{10}(n)$)
but results at other $\alpha$ might be useful in contexts
other than the pde problem of this paper.

\subsection{An aside: (RE(0)) and (RE(1)): Chebyshev}\label{subapp:RE0}

At $\alpha=0$ the $n$ cancels and the recurrence is
constant coefficient, so solvable by $u_n=r^n$ for  appropriate $r$.
The general solution is a linear combination of the Chebyshev polynomial solutions
$$u_n(\alpha=0)= c_T(q_2) T_n(q_2) + c_U(q_2) U_n(q_2), $$
where $c_T$ and $c_U$ are independent of $n$.
As $n$ no longer occurs in the coefficients in the recurrence
(RE(0)), equally good as a polynomial solution is $u_n=T_{n-m}$ for integer $m$.

Chebyshev polynomials also arise in solving the recurrence
when $\alpha=1$:
solutions are 
$$u_n(\alpha=1)=\frac{u_n(\alpha=0)}{n} . $$

As remarked before, solutions can be expressed in terms of Legendre functions.
For the $\alpha=0$, $\mu=1/2$ case, see~\cite{HTF1}3.6.1(12).
Also we have
\begin{eqnarray*}
P^{1/2}_{n-1/2}(\cos(\theta))
&=& \sqrt{\frac{2}{\pi\sin(\theta)}} \cos(n\theta)
= \sqrt{\frac{2}{\pi\sin(\theta)}} T_n(\cos(\theta)) ,
 \\
Q^{1/2}_{n-1/2}(\cos(\theta))
&=& -\sqrt{\frac{\pi}{2\sin(\theta)}} \sin(n\theta) ,
\\
P^{1/2}_{n-1/2}(\cosh(\eta))
&=& \sqrt{\frac{2}{\pi\sinh(\eta)}} \cosh(n\eta) 
= \sqrt{\frac{2}{\pi\sinh(\eta)}} T_n(\cosh(\eta)),
 \\
Q^{1/2}_{n-1/2}(\cosh(\eta))
&=& -\sqrt{\frac{\pi}{2\sinh(\eta)}} \exp(-n\eta) .
\end{eqnarray*}
We remark that
\begin{eqnarray*}
U_n(x)
&=& \frac{x\, T_{n+1}(x) - T_{n+2}(x)}{1-x^2} ,\\
U_n(\cosh(\eta))
&=& \frac{\sinh((n+1)\eta)}{\sinh(\eta)} ,\\ 
\exp(-n\eta)
&=&T_n(\cosh(\eta))-\sinh(\eta)U_{n-1}(\cosh(\eta)),
\end{eqnarray*}
and $u_n=\exp(-n\cosh^{-1}(q_2))$ solves (RE(0)).
Similar formulae are available for the case $\alpha=1$,
order $\mu=-1/2$.

Maple's {\tt rsolve} solves (RE($m$)) for other integer $m$.

None of the items in this elementary subsubsection are relevant to
the $\alpha=1/2$ and $\alpha=3/2$ cases which arise
in our partial differential equation problem.
However, relations between families of polynomials and
various connections between Legendre functions will arise in
the context of the $\alpha=1/2$ and $\alpha=3/2$ cases
and the easier situation sometimes suggests possible results.

\subsection{Deriving (RE($\alpha$)) from integrals defining Fourier coefficients}\label{subapp:fromFS}

While the expressions for $g_n$ and ${\hat g}_n$ in terms of toroidal functions have value, much can be derived, simply from their definitions, including the recurrences (RE($\alpha$))
and inter-relationships.
The simple relationships between $g$ and ${\hat g}$ lead to
relationships between their Fourier coefficients.
\begin{itemize}
\item Since
$$ {\hat g}(\psi)=\left(1- \frac{\cos^2(\psi)}{q}\right) \, g(\psi) 
= \frac{1}{2 q} \left( (2q-1)-\cos(2\psi)\right)\, g(\psi) ,$$
we have by multiplying both sides by $\cos(2n\psi)$, using cosine formulae, and integrating
\begin{equation}
4 q {\hat g}_n = 2 (2q-1) g_n - g_{n+1}- g_{n-1}  .
\label{eq:cross0}
\end{equation}

\item Since
$$ 2 q \frac{ d {\hat g}(\psi)}{d q} +  {\hat g}(\psi) - g(\psi) = 0, $$
the Fourier coefficients $g_n$ can be found from those of $ {\hat g}$:
\begin{equation}
g_n =  {\hat g}_n  +  2 q \frac{ d {\hat g}_n}{d q} = 2\sqrt{q} \frac{d}{dq} (\sqrt{q}{\hat g}_n ) .
\label{eq:gDq}
\end{equation}

\item Since
$$
\frac{ d {\hat g}(\psi)}{d q}
=\frac{\cos(\psi)^2 \, g}{2 q^2} 
= \frac{(\cos(2\psi)+1 \, g}{4 q^2} ,
$$
\begin{equation}
 \frac{ d {\hat g}_n}{d q} 
= \frac{1}{8q^2} (2g_n +g_{n-1}+g_{n+1}) .
\label{eq:ghDq}
\end{equation}
(Equations~(\ref{eq:gDq}) and~(\ref{eq:ghDq}) check against~(\ref{eq:cross0}).)

\end{itemize}

Define also
$$ q_2 = 2 q -1 . $$
There are three-term recurrence relations between the coefficients:
\begin{eqnarray}
(n-\frac{1}{2}) g_n
&=& 2 q_2 (n-1) \, g_{n-1} - (n-\frac{3}{2}) \, g_{n-2} ,
\label{eq:gnrec}\\
(n+\frac{1}{2}) {\hat g}_n
&=& 2 q_2 (n-1) \, {\hat g}_{n-1} - (n-\frac{5}{2} )\, {\hat g}_{n-2} .
\label{eq:gBarnrec}
\end{eqnarray}
The derivation is treated below, where we derive the following (RE($1/2$))
equivalent to~(\ref{eq:gnrec}):
$$
(n+\frac{1}{2}) g_{n+1}
= 2 q_2 n \, g_{n} - (n-\frac{1}{2}) \, g_{n-1} .
\eqno{({\rm RE}(1/2))}
$$ 
Corresponding to~(\ref{eq:gBarnrec}) we have (RE($3/2$))
$$
(n+\frac{3}{2}) {\hat g}_{n+1}
= 2 q_2 n \, {\hat g}_{n} - (n-\frac{3}{2}) \, {\hat g}_{n-1} .
\eqno{({\rm RE}(3/2))}
$$ 
We now repeat equations given in the main part of the paper.
In these  the argument  of the elliptic integrals is $1/\sqrt{q}$ and we wrote,
\begin{eqnarray}
g_n
&=& \frac{4}{\pi}\, \left( E_n {\rm EllipticE}_0  + K_n {\rm EllipticK}_0 \right),
\label{eq:gnEll}\\
 {\hat g}_n
&=&  \frac{4}{\pi}\, \left( {\hat E}_n {\rm EllipticE}_0  + {\hat  K}_n {\rm EllipticK}_0 \right) ,
\label{eq:gBarnEll}
\end{eqnarray}
with the $E_n$ ,  $K_n$,  ${\hat E}_n$, $ {\hat  K}_n$ polynomials in $q_2$ of degree at most $n$.
The polynomials  $E_n$ ,  $K_n$ satisfy the recurrence~(\ref{eq:gnrec}), (RE(1/2)):
the polynomials  ${\hat E}_n$ ,  ${\hat K}_n$ satisfy the recurrence~(\ref{eq:gBarnrec}), (RE(3/2)).
The starting values for the iterations to determine the polynomials are given,
for the $g_n$ iteration by (\ref{eq:g0Ell}),~(\ref{eq:g1Ell}) and
for the ${\hat g}_n$ iteration by (\ref{eq:gBar0Ell}),~(\ref{eq:gBar1Ell}).
Maple's rsolve will solve (RE($\alpha$)) in the special cases when $q_2=0$ and $q_2=1$.
When $q_2=1$ there are constant solutions ($g_n=1$ for all $n$), and reduction of order gives the other.



The method to establish (\ref{eq:gnrec}),~(\ref{eq:gBarnrec})
involves a further
relation between the coefficients $g_n$, ${\hat g}_n$ of the two series,
and is suggested in a posting by Jack D'Aurizio on stackexchange.
See\\
{\small
\verb$math.stackexchange.com/questions/930003/fourier-series-of-sqrt1-k2-sin2t$\\
}
We have
\begin{equation}
g_{m+1}= (2 q-1) g_m + (4 m-2) q \, {\hat g}_m ,
\label{eq:cross1}
\end{equation}
and also, exactly as in D'Aurizio's post,
\begin{equation}
8 m q {\hat g}_m = g_{m+1}- g_{m-1} .
\label{eq:cross2}
\end{equation}
Eliminating ${\hat g}_m$ between
equations~(\ref{eq:cross1}),~(\ref{eq:cross2}) we find
the recurrence~(RE(1/2)) 
with $n$ there replaced by $m$.
Also, equations~(\ref{eq:cross1}) and~(\ref{eq:cross2})
together yield
\begin{equation}
g_{m-1}= (2 q-1) g_m - (4 m+2) q \, {\hat g}_m ,
\label{eq:cross1a2}
\end{equation}

It may be that, for $q_2>1$,  all the $g_m$ are positive, and the sequence is decreasing.\\
We have ${\hat g}_0>0$ and it may be that, for $q_2>1$ and $m\geq{1}$ that the terms ${\hat g}_m$ are
all negative and form an increasing sequence.

\subsubsection{Generating functions?}

Correcting a misprint in~\cite{Ch78} p202, equation~(12.5),
(which reports~\cite{BD67}) for $\alpha=1/2$,
a generating function is
\begin{eqnarray}
G(x,w)
&=& \frac{1}{2}\,\frac{1}{\sqrt{w(1-2xw+w^2)}}\,
\int_0^w  \frac{1}{\sqrt{t(1-2xt+t^2)}} \, dt ,
\label{eq:GenDef}\\
G(x,w)
&=& \sum_{n=0}^\infty p_{01}(n+1,x) w^n .
\label{eq:Genp01}
\end{eqnarray}
The $p_{01}(n)$ are polynomials of degree $n-1$
satisfying (RE(1/2)) starting from $p_{01}(0)=0$, $p_{01}(1)=1$.
These are studied extensively in later parts of this  section. 
(Equation~(\ref{eq:p01Leg}) follows from this: see~\cite{Ch78} equation~(12.6).)

This subsection is `incomplete', 
with everything after this paragraph
merely an early personal attempt predating finding
the result reported in equation~(\ref{eq:GenDef}) above.
It may be that, for general $\alpha$, elaborate special functions, 
e.g. Apell hypergeometrics, may be needed as they were in the
more general situation in~\cite{BI82}.
At least for our special case of $\alpha=1/2$
incomplete elliptic integrals suffice:
see also equations (1.10), (1.11) of ~\cite{VZ07}.

The stackexchange posting suggests there might be reasonably neat expressions for
generating functions,
$$ G(X) = g_0 + g_1\, X +\sum_{n=2}^\infty g_n X^n, \qquad
 {\hat G}(X) = {\hat g}_0 + {\hat g}_1\, X +\sum_{n=2}^\infty {\hat g_n} X^n ,$$
 and that these functions satisfy first order linear differential equations.\\
 The differential equation satisfied by $G$ is
 $$ {\cal L} G =(X-2 q_2 X^2 +X^3)\frac{d G}{d X} -\frac{1}{2} (1-X^2) G = \frac{1}{2} (-g_0 +g_1 X) .$$
 $X=0$ is a singular point.
 The homogeneous de ${\cal L}G_h=0$ is solved by a function $G_h$ with $G_h(0)=0$,
 $$ G_h(X) = \sqrt{\frac{X}{1-2q_2 X + X^2}} . $$
 We would be content with representations of $G$ valid for $q_2>1$ and $0\leq{X}<1/(2q_2)$.
 Applying the `variation of parameters' formula to the nonhomogenous de above,
 using $G_h$, leads to a messy particular solution in terms of 
incomplete elliptic integrals.
(Though possibly unrelated, we will see elliptic integrals in other contexts:
see~\S\ref{subapp:orthog}.)
\\
The corresponding de for ${\hat G}$ is:
$$ {\hat{\cal L}}{\hat G} =(X-2 q_2 X^2 +X^3)\frac{d {\hat G}}{d X} +
\frac{1}{2} (1-X^2) {\hat G}
= \frac{1}{2} ({\hat g}_0 +3{\hat g}_1 X) .$$
 $X=0$ is a singular point.
 The homogeneous de ${\cal L}{\hat G}_h=0$ is solved by a function ${\hat G}_h$ with ${\hat G}_h(0)=0$,
 $$ {\hat G}_h(X) = \sqrt{\frac{1-2q_2 X + X^2}{X}} , $$
(which is the reciprocal of $G_h$).
The relations~(\ref{eq:cross0},\ref{eq:gDq}) lead to
relations between $G$ and $\hat G$, for example
$$  G= 2\sqrt{q}\frac{d}{dq} (\sqrt{q} {\hat G}) . $$
One can derive a differential equation for ${\hat G}(q)$ involving derivatives with respect to $q$.
It may be that if one starts from a closed form for ${\hat G}(1)$, i.e. at $q=1$, solving the initial value problem
might be useful.
\medskip

\subsubsection{Hypergeometric representations}

Also presented on stackexchange is a 
hypergeometric formula for ${\hat g}_n$.
See also~\cite{Cv09} at Lemma~1 which gives a
slightly different representation.
Define
$$ \lambda= \frac{1}{\sqrt{q}+\sqrt{q-1}}
= \sqrt{q}-\sqrt{q-1} , $$
from which
$$ \lambda^2 = \frac{1}{q_2 +\sqrt{q_2^2-1}} . $$
We have
\begin{equation}
{\hat g}_n = - \frac{2}{\lambda\sqrt{q}}\,
\frac{(2n)!}{(2n-1)\,2^{2n+1}\, (n!)^2}\, (\lambda^2)^n\,
{\ }_2{F}_1(-\frac{1}{2},n-\frac{1}{2};n+1;\lambda^4) .
\label{eq:gHathy}
\end{equation}
(Quick check, ${\hat g}_0>0$ and, for $n>0$, ${\hat g}_n<0$,
though we only have numeric evidence for this.)

\subsubsection{Negative $n$}

We remark that from the $\cos(2n\psi)$ in the defining integrals
we have
$$ g_{-n}= g_n, \qquad {\hat g}_{-n}={\hat g}_n . $$
See also equations~(\ref{eq:AS822}),~(\ref{eq:Pminus}).
It is also more generally true that if
$u_n$ is a solution of (RE($\alpha$)) and $\alpha\ne{0}$ then
$u_{-n}=u_n$.

\subsection{Toroidal functions, LegendreQ functions}\label{subapp:Toroidal}

\subsubsection{Main facts concerning $Q$ functions}

We defer to a later subsection items on LegendreP functions.
LegendreP functions are equally appropriate as solutions of
(RE(1/2)) and (RE(3/2)) but our initial concern is with
the Fourier coefficients $g_n$ and ${\hat g}_n$ and these
are LegendreQ functions.
Also in this subsection we are interested in the functions
defined on the interval $(1,\infty)$ though in some 
later sections we are also concerned with the interval $(-1,1)$.
(Branch points will then need to be considered even though our application
involves real-valued functions of a real argument.)

Toroidal functions are special cases of Legendre functions,
specifically $Q^\mu_{n-1/2}$ and $P^\mu_{n-1/2}$ with $n$ integer,
and, sufficient for our purposes also $\mu$ integer.
We have, as in~\cite{HTF1}{\S}3.10(3) and other references such as~\cite{Co07}
\begin{eqnarray}
g_n
&=& \frac{4\sqrt{q}}{\pi}\, Q_{n-1/2}(q_2) ,
\label{eq:gnQ}\\
{\hat g}_n
&=& \frac{2\sqrt{q-1}}{\pi}\, Q^{-1}_{n-1/2}(q_2) .
\label{eq:gBarnQ}
\end{eqnarray}
(See equations~(43) and (44) of~\cite{Co07}.)
 
Toroidal functions can be related to Elliptic Integrals. 
In particular
\begin{eqnarray}
Q_{-1/2}(q_2)
&=&\frac{1}{\sqrt{q}}\,  {\rm EllipticK}(\frac{1}{\sqrt{q}})
 ,
 \label{eq:Qmh}
\\
Q_{1/2}(q_2)
&=& q_2 \frac{1}{\sqrt{q}}\,  {\rm EllipticK}(\frac{1}{\sqrt{q}} )-
2\sqrt{q} \, {\rm EllipticE}(\frac{1}{\sqrt{q}} )  .
 \label{eq:Qph}
 \end{eqnarray}
That for $Q_{-1/2}$ is \cite{AS}~8.13.3: that for $Q_{1/2}$ is \cite{AS}~8.13.7.
The corresponding formulae for the Legendre $P$ functions,
$P_{-1/2}(q_2)$ and $P_{1/2}(q_2)$,
though not immediately relevant to the $g_n$ and ${\hat g}_n$ are relevant to
solving (RE($\alpha$)) and are treated in a later section.
We also remark, as in~\cite{CD10} p339,
{\small
\begin{eqnarray}
Q^{1}_{-1/2}(q_2)
&=&  - \frac{1}{2\sqrt{q-1}}\,  {\rm EllipticE}(\frac{1}{\sqrt{q}})
 ,
 \label{eq:Qp1mh}
\\
Q^{1}_{1/2}(q_2)
&=&   - \frac{q_2}{2\sqrt{q-1}}\,  {\rm EllipticE}(\frac{1}{\sqrt{q}}) +\sqrt{q-1} {\rm EllipticK}(\frac{1}{\sqrt{q}} )
 .
 \label{eq:Qp1ph}
\end{eqnarray}
}
and, looking ahead, from~(\ref{eq:Qp1m1})
\begin{eqnarray}
Q^{-1}_{-1/2}(q_2)
&=&   -4\, Q^{1}_{-1/2}(q_2)
 ,
 \label{eq:Qm1mh}
\\
Q^{-1}_{1/2}(q_2)
&=&   \frac{4}{3} Q^{1}_{1/2}(q_2)
 .
 \label{eq:Qm1ph}
\end{eqnarray}
 
\cite{AS}~8.2.2 gives, for integer $\mu$,
\begin{equation}
 Q_{-n-1/2}^\mu(z) = Q_{n-1/2}^\mu(z) . 
 \label{eq:AS822}
\end{equation} 
 
\cite{HTF1}{\S}3.8(9) gives
\begin{equation}
 \frac{d Q_{n-1/2}(z)}{dz}
 = \frac{n^2-1/4}{\sqrt{z^2-1}} \, Q^{-1}_{n-1/2} .
 \label{eq:HTF1Qp3p8p9}
\end{equation}

As noted above, the $g_n$ are given by a function of $q_2$ times $Q_{n-1/2}(q_2)$,
the latter functions solving (RE($1/2$)).
The ${\hat g}_n$ are given by a function of $q_2$ times $Q^{-1}_{n-1/2}(q_2)$,
the latter functions solving (RE($3/2$)).
Where it seems appropriate we prefer to arrange our formulae so that the $Q^{-1}_{n-1/2}$
are written in terms of $Q_{m-1/2}$, for example from 
\cite{HTF1}\S3.8(5), 
\begin{equation}
Q^{-1}_{n-1/2}(z)
= \frac{z\,Q_{n-1/2}(z)-Q_{n-3/2}(z)}{(n+1/2)\sqrt{z^2-1}} .
\label{eq:HTF385}
\end{equation}
(This agrees with (\ref{eq:cross1a2}).)
Also, as in~\cite{HTF1}{\S}3.8(3)
\begin{equation}
Q^{-1}_{n-1/2}(z)
= \frac{Q_{n+1/2}(z)-Q_{n-3/2}(z)}{2n\sqrt{z^2-1}} .
\label{eq:HTFQ3p8p3}
\end{equation}
(See also Theorem~\ref{thm:thm1} below.)
Particular cases of equation~(\ref{eq:HTF385}), 
at $n=0$ and $n=1$,
are:
\begin{eqnarray}
Q^{-1}_{-1/2}(z)
&=& \frac{z\,Q_{-1/2}(z)-Q_{1/2}(z)}{(1/2)\sqrt{z^2-1}} ,
\label{eq:HTF0Q}\\
Q^{-1}_{1/2}(z)
&=& \frac{z\,Q_{1/2}(z)-Q_{-1/2}(z)}{(3/2)\sqrt{z^2-1}}  .
\label{eq:HTF1Q}
\end{eqnarray}  
We note, as in \cite{AS}~8.2.6,
\begin{equation}
Q^{1}_{n-1/2} (z) =  (n^2-1/4) Q^{-1}_{n-1/2} (z) 
\label{eq:Qp1m1}
\end{equation}
which can be deduced from \cite{HTF1}~\S3.8 equations (5) and (7).
(This accords with Theorem~\ref{thm:thm3a} given in the next section.)

As noted above Conway~\cite{Co07} and various other authors give
formulae for $g_n$ (see~\cite{Co07} equation (43))
and for ${\hat g}_n$
in terms of Legendre $Q$ functions.
Concerning the latter, see~\cite{Co07} equation (44),
\begin{eqnarray}
{\hat g}_n
&=& \frac{2(\sqrt{q-1})}{\pi} Q^{-1}_{n-1/2}(q_2)
=  \frac{2(\sqrt{q-1})}{\pi} \frac{1}{n^2-1/4} Q^1_{n-1/2}(q_2) ,  \label{eq:ConQh1}\\
&=& \frac{2(\sqrt{q-1})}{\pi}\left( \frac{q_2\,Q_{n-1/2}(q_2)-Q_{n-3/2}(q_2)}{(n+1/2)\sqrt{q_2^2-1}} \right)
 .
 \label{eq:ConQh2}
 \end{eqnarray}
(See also~\cite{CD10} p339.
It also accords with the corollary to Theorem~\ref{thm:thm1}.)
 
As a check we recall~(\ref{eq:gBarnq1})
 $${\hat g}_n(q_2=1)= - \frac{1}{\pi(n^2-1/4)} .$$
 \cite{HTF1}~\S3.9(6) gives asymptotics for $q_2$ tending down to 1:
 $$ Q^{-1}_{n-1/2}(q_2) \sim -\frac{1}{\sqrt{2}}\frac{\Gamma(n-1/2)}{\Gamma(n+3/2)}\frac{1}{\sqrt{q_2-1}}
 =   -\frac{1}{\sqrt{2}}\frac{1}{n^2-1/4} \frac{1}{\sqrt{2q-2}}. $$ 
 Using this last expression in~(\ref{eq:ConQh1}) checks with the elementary result that
 ${\hat g}_n(q_2=1)= - {1}/({\pi(n^2-1/4)})$.
 
\subsubsection{Integrals of sums of squares of our Fourier coefficients}

Using Parseval's Theorem it is easy to find 
$\sum_0^\infty{\hat g}_n^2$, etc. 

See also~\cite{Co07}\S{5}.

\subsection{(RE($\alpha$)) and polynomial solutions}\label{subapp:REalphaPoly}

The earlier subsections have had a narrow focus on the $g_n$ and ${\hat g}_n$ 
Fourier coefficients.
In this subsection we will focus on polynomial solutions
of (RE($\alpha$)).
Since the general solution of the (RE($\alpha$)) difference equation can be formed from a
linear combination of $Q^\mu_{n-1/2}$ and $P^\mu_{n-1/2}$
we will, later in this section,  be led to consider the Legendre $P$ functions.

\subsubsection{Easy results}

While too elementary to suggest immediate generalization, we note that when we set $x=q_2=1$ in (RE($\alpha$)), the recurrence is easy to solve.
For all $\alpha$ one solution is $u_n=1$ for all $n$.
Reduction of order can be used to find another solution.

In particular (RE(1/2)) is then solved in terms of the digamma function $\psi$, specifically by $u_n=\psi(n+1/2)$.
We have
$$ \psi(n+\frac{3}{2}) = \psi(n+\frac{1}{2}) +\frac{1}{n+1/2} .$$
Forming the corresponding equation with $n$ replaced by $n-1$, and subtracting the two equations we find a solution.
With
\begin{eqnarray*} 
u_n
&=&  (\psi(n+1/2)-\psi(1/2))/2,
\qquad u_0=0,\ u_1=1 ,\\
u_n
&=&  1-(\psi(n+1/2)-\psi(1/2))/2,
\qquad u_0=1,\ u_1=0 ,
\end{eqnarray*}
we have solutions to the case $x=q_2=1$ of (RE(1/2))

When $\alpha=3/2$ the formulae are simpler.
Then we have
that the sequence ${\hat u}_n=(n^2-1)/(4n^2-1)$
satisfies (RE(3/2)) and has ${\hat u}_0=1$ and ${\hat u}_1=0$.
A simple relationship between the ${\hat u}_n$ of this paragraph
and the $u_n$ of the preceding is
$$ \frac{2}{4n^2-1} 
= \frac{\psi(n+3/2)-\psi(n-1/2)}{n} .$$
This can be regarded as a simple lead-in to Theorem~\ref{thm:thm1}.

Next some results follow from very routine algebraic manipulation.
Define
$${\cal L}(\alpha,u,n):=
2 n q_2 u_n - (n+\alpha) u_{n+1} - (n-\alpha) u_{n-1} .
$$

\begin{theorem}\label{thm:thm1}
If $u_n$ solves {\rm{(RE(1/2))}} and
$$ {\hat u}_n = 
\frac{1}{n} \left( u_{n+1}- u_{n-1}\right ) , $$
then $ {\hat u}_n $  solves {\rm (RE(3/2))}.
\end{theorem}

\noindent{\it Proof.}
\begin{eqnarray*}
{\cal L}(\frac{3}{2},{\hat u},n)
&=& 2 n q_2 \frac{u_{n+1}-u_{n-1}}{n} - 
(n+\frac{3}{2})  \frac{u_{n+2}-u_{n}}{n+1} -
 (n-\frac{3}{2})  \frac{u_{n}-u_{n-2}}{n-1} , \\
 &=& \frac{1}{n^2-1}\left(
(n-1) {\cal L}(\frac{1}{2},{ u},n+1) +
(n+1) {\cal L}(\frac{1}{2},{ u},n-1)\right) ,\\
&=& 0 .
\end{eqnarray*}
This establishes the theorem.
\medskip

\noindent{\bf Corollary.}
{\it If $u_n$ solves {\rm{(RE(1/2))}} and
$$ {\hat u}_n
=\frac{x u_{n}- u_{n-1}}{2n+1} , $$
then $ {\hat u}_n $  solves {\rm (RE(3/2))}, and
so too does
$$ {\hat u}_n 
=\frac{x u_{n}- u_{n+1}}{2n-1} . $$}
(The latter is suggested by~\cite{HTF2}3.8(6).)

\noindent{Remarks.}
\begin{itemize}
\item
We have already seen in~(\ref{eq:cross2}) that
$$ {\hat g}_n = \frac{1}{4 n (q_2+1)}
\left( g_{n+1}- g_{n-1}\right ) .$$
See also \cite{HTF1}~3.8.3

\item With the definitions in the later subsection on
`Polynomial sequences' we have the further examples
\begin{eqnarray}
-4{\hat p}_{01}(n)+12 x{\hat p}_{10}(n)
&=& \frac{3}{n}\left(p_{10}(n+1)-p_{10}(n-1)\right) ,
\label{eq:L1ppBarn}
 \\
4 x{\hat p}_{01}(n)-12{\hat p}_{10}(n)
&=& \frac{3}{n}\left(p_{01}(n+1)-p_{01}(n-1)\right) .
\label{eq:L2ppBarn}
\end{eqnarray}
(Substituting $x=1$ and adding the equations has
each side summing to 0 which provides one simple check.)

\end{itemize}

\medskip
\begin{theorem}\label{thm:thm2} 
If $u_n$ solves \rm{(RE(1/2))} and
$$ {\hat u}_n = 
\left( 2q_2 u_n - u_{n+1}- u_{n-1}\right ) , $$
then $ {\hat u}_n $  solves {\rm (RE(3/2))}.
\end{theorem}

\noindent{Remark.}
We have already seen in~(\ref{eq:cross0}) that
$$ {\hat g}_n = \frac{1}{2(q_2+1)}
\left( 2q_2 g_n - g_{n+1}- g_{n-1}\right ) .$$

\medskip

There are many miscellaneous results.
For example a routine calculation shows that 
\begin{theorem}\label{thm:thm3a}
If ${\hat u}_n$ satisfies {\rm (RE(3/2))}, then
$(n^2-1/4){\hat u}_n$ satisfies {\rm (RE(-1/2))}.
\end{theorem}
(See also~\cite{HTF2}3.8(1) with $\mu=-1$.)

\subsubsection{Equivalent recurrence relations}

There are many equivalent forms of recurrence relations.

\medskip
\begin{itemize}
\item
Defining, from a sequence satisfying (RE(3/2)),
the $a_n$ sequence via
$$ {\hat u}_{n+1} = a_{n+1} \, (2 q_2)^n \,
\frac{\Gamma(n+1)\Gamma(5/2)}{\Gamma(n+5/2)}
= a_{n+1} \, (2 q_2)^n \,
\frac{3\, 2^n\, n!}{(2n+3)!!} ,
$$
we have, for $n>1$
$$ a_{n+1}= a_n +
\frac{(n-3/2)(n+1/2)}{4\, q_2^2\, n(n-1)} a_{n-1} .
\eqno{\rm{(RE11(3/2))}}
$$
(See also~\cite{VZ07} equations~(5.12) and~(5.17)
for (RE(1/2)).)
While ${\hat u}_n$ might, for some initial conditions be polynomial in $q_2$ degree $n-1$, the $a_n$ are no longer polynomial,
merely rational.
Also, we cannot set $n=1$ in the last equation.

\item
Some studies of recurrence equations for sequences of  polynomials
involve recurrences for monic polynomials, of the form
$$ v_{n+1}= x\,v_n - f(n)\, v_{n-1} . \eqno{{\rm (RE1xf)}}
$$
Starting with a sequence of polynomials $u_n$ satisfying 
(RE($\alpha$)),
one sets $v_n=u_n/g(n)$ and proceeds as follows:
$$ g(n+1) (n+\alpha) v_{n+1}
= 2 n x g(n) v_n - g(n-1) (n-\alpha) v_{n-1} .$$
For this to be of the required form we must have
$$ g(n+1) (n+\alpha) = 2 n g(n) ,$$
which solves to
$$ g(n) = c\, \frac{n!}{\Gamma(n+\alpha)} . $$
From this
$$ f(n)=\frac{g(n-1) (n-\alpha)}{g(n+1) (n+\alpha)}
= \frac{(n-\alpha)(n-1+\alpha)}{4(n-1)n} .$$
Thus, for $n>1$
$$ v_{n+1}= x\, v_n +
\frac{(n-\alpha)(n-1+\alpha)}{4(n-1)n}\, v_{n-1} .
\eqno{(\rm{RE1x}(\alpha))}
$$
Again we cannot set $n=0$ or $n=1$ in the preceding
(except, taking appropriate limits, in the cases $\alpha=0$
and $\alpha=1$).
For the $\alpha=1/2$ case see~\cite{Ch78}~Chapter~6, equation~(12.9).
See also~\cite{VZ07} equations~(5.12) and (5.17),
but note, as occurs with some other references, the indexing
is different from ours, in this case their $S_n$ is our
$v_{n+1}$ (monic form of $p_{01}$).

There is another well-known approach to the calculation of $f$ for
(RE1xf).
Let $(v_n)$ be a sequence of monic polynomials, mutually orthogonal.
We do not, at this stage, need to know the weight function
but we will assume it is positive, even in $x$ and the interval is
symmetric about $x=0$.
Taking inner products of each side of (RE1xf) with $v_{n-1}$ gives
$$ 0 = \langle x v_n, v_{n-1}\rangle - f_n \langle v_{n-1}, v_{n-1}\rangle .$$
Shifting the subscripts of the terms in (RE1xf) down by 1,
i.e. $n$ replaced by $n-1$, and then taking the inner product with $v_n$
gives
$$\langle v_n, v_{n}\rangle 
= \langle x v_{n-1}, v_{n}\rangle . $$
The preceding two equations yield
$$ f_n = \frac{\langle v_n, v_{n}\rangle}
{\langle v_{n-1}, v_{n-1}\rangle} .$$
We can apply this at $\alpha=1/2$ to our polynomial sequence
$p_{01}(n)$ defined later, see \S\ref{subapp:PolyQ} and (\ref{eq:p01PQ}).
The weight function is known explicitly and calculations 
(from~\cite{BD67}, reported in various places including ~\cite{Ch78})
give
$$ \frac{\langle p_{01}(n), p_{01}(n)\rangle}
{\langle p_{01}(n-1), p_{01}(n-1)\rangle}
=\frac{n-1}{n} .$$
We need monic polynomials. So
$$ f_n = \frac{n-1}{n}\,
\left(\frac{{\rm lcoeff}(p_{01}(n-1)}
{{\rm lcoeff}(p_{01}(n)}\right)^2 . $$
Using the result on the leading coefficient given in
equation~(\ref{eq:lc01}) we find the formula gives
$$ f_n =\frac{n-1}{n}\,
\left( \frac{2n-1}{4(n-1} \right)^2 ,$$
which agrees with the $\alpha=1/2$ case given in
(RE1x($\alpha$) above.

\item A different transformation $v_n=\psi(n,\alpha)\sqrt{n}\, u_n$,
with $\psi$ specified later, but with $\psi(n,1/2)=1$, leads to another
family of recurrence relations, of the form
$$ a_n v_{n+1} =  x v_n - a_{n-1} v_{n-1} . 
\eqno{({\rm RE}a_n{\rm x}(\alpha))}$$
The recurrence $({\rm RE}a_n{\rm x}(\alpha))$ can be written
\begin{equation}
\begin{pmatrix}
v_{n+1}\\
a_n v_n 
\end{pmatrix}
=A(n) \,
\begin{pmatrix}
v_{n}\\
a_{n-1} v_{n-1} 
\end{pmatrix}
\ {\rm where\ }
A(n)=\begin{bmatrix}
\frac{x}{a_n}& -\frac{1}{a_n}\\
a_n& 0
\end{bmatrix} .
\label{eq:Amat}
\end{equation}
(We remark that $({\rm RE}a_n{\rm x}(\alpha))$ can be transformed to
(RE1xf), the latter with a different $v$ of course,
with $f(n)=a_{n-1}^2$.)
The fact that ${\rm det}(A(n))=1$ leads to generalizations of
the result we have at equation~(\ref{eq:niceph}).
Specialise now to our (RE($\alpha$)).
For (RE(1/2)) we have
\begin{equation}
a_n = \frac{n+1/2}{2\sqrt{n(n+1)}} .
\label{eq:anDef}
\end{equation}
More generally, for (RE($\alpha$)) we have
$$\psi(n,\alpha) =\sqrt{\frac{\Gamma(n+\alpha)}{\Gamma(n+1-\alpha}},\quad 
a_n = \frac{\sqrt{\Gamma(n+1+\alpha)\Gamma(n+2-\alpha)}}
{\sqrt{n(n+1)\Gamma(n+1-\alpha)\Gamma(n+\alpha)}} .
$$
Returning to the case $\alpha=1/2$,
on using $a_n$ of equation~(\ref{eq:anDef}) in
$A(n)$ of equation~(\ref{eq:Amat}), we find
that the product of matrices
$$ A(n)\, A(n-1)\, \ldots A(2)\, A(1)
= \begin{bmatrix}
\sqrt{n+1}\,p_{01}(n+1)& 4\, \sqrt{n+1}\,{ p}_{10}(n+1)\\
\sqrt{n}\, a(n)\,p_{01}(n)& 4\, \sqrt{n}\,a(n)\, {p}_{10}(n) 
\end{bmatrix} ,
$$
where $p_{01}$ and $p_{10}$ are the polynomials satisfying
(RE(1/2)) starting 0,1 and 1,0 respectively.

\item
It is possible that further work related to the items in
this paragraph might be productive.
We record them in the hope that some reader might suggest
how it might be developed.

\begin{itemize}
\item
Next a definition. See~\cite{Ch78}\S5.5.
A sequence of positive numbers $(f_n)$ is called a {\it chain sequence}
if there exists another sequence $(g_n)$ such that
$$ f_n = g_n\, (1-g_{n-1}) ,\ {\rm with\ }
0\le g_0<1,\ 0<g_n<1 {\ \rm for\ } n=1,2,\ldots \ \ .
$$
The $f(n,\alpha)=a_{n-1}^2$ above is a chain sequence.
Use
$$ g_n=\frac{1}{2}-\frac{1}{4n} \ {\rm for\ \ }
\alpha=\frac{1}{2} . $$

\item In some papers the 
$a_n$ is written as
a ratio $a_n=\gamma(n)/\gamma(n-1)$.
For $\alpha=1/2$ the $a_n$ of equation~(\ref{eq:anDef}) has
$$ a_n=\frac{\gamma(n)}{\gamma(n-1)}\ \ {\rm with\ \ }
\gamma(n) = \frac{2^{n-1}\, n!}{\sqrt{n}}\,
\frac{\Gamma(1/2)}{\Gamma(n+1/2)} .
$$

\end{itemize}

\end{itemize}

\subsubsection{Legendre $P_{n-1/2}$ functions}

Formulae. involving Legendre $P_{n-1/2}$ functions,
 similar to those in
(\ref{eq:Qmh},\ref{eq:Qph},\ref{eq:Qp1mh},\ref{eq:Qp1ph}) can be found.
{\small
\begin{eqnarray}
P_{-1/2}(\cos(\theta))
&=& \frac{2}{\pi} {\rm EllipticK}(\sin(\frac{\theta}{2})) ,
\label{eq:PmhDL}\\
P_{1/2}(\cos(\theta))
&=& \frac{2}{\pi}\left(
2 {\rm EllipticE}(\sin(\frac{\theta}{2})) -
 {\rm EllipticK}(\sin(\frac{\theta}{2}))\right) ,
\label{eq:PphDL}\\
Q_{-1/2}(\cos(\theta))
&=& {\rm EllipticK}(\cos(\frac{\theta}{2})) ,
\label{eq:QmhDL}\\
Q_{1/2}(\cos(\theta))
&=& 
 {\rm EllipticK}(\cos(\frac{\theta}{2})) -
2 {\rm EllipticE}(\cos(\frac{\theta}{2})) ,
\label{eq:QphDL}\\
P_{-1/2}(\cosh(\xi))
&=& \frac{2}{\pi\cosh(\frac{\xi}{2})}
{\rm EllipticK}(\tanh(\frac{\xi}{2})) ,
\label{eq:PmhDLx}\\
P_{1/2}(\cosh(\xi))
&=& \frac{2\exp(\xi/2)}{\pi}
{\rm EllipticE}(\sqrt{1-\exp(-2\xi)}) ,
\label{eq:PphDLx}\\
Q_{-1/2}(\cosh(\xi))
&=& \frac{2}{\sqrt{\pi}}\exp(-\frac{\xi}{2})
{\rm EllipticK}(\exp(-{\xi})) ,
\label{eq:QmhDLx}\\
Q_{1/2}(\cosh(\xi))
&=& \frac{2}{\sqrt{\pi}}
\frac{\cosh(\xi)}{\cosh(\frac{\xi}{2})}
{\rm EllipticK}(\frac{1}{\cosh(\frac{\xi}{2})}) -\\
&\ &\qquad\qquad
\frac{4}{\sqrt{\pi}}\cosh(\frac{\xi}{2})
{\rm EllipticE}(\frac{1}{\cosh(\frac{\xi}{2})})
 .\label{eq:QphDLx}
\end{eqnarray} 
}
See~\cite{AS}8.13. There are many equivalent ways of
expressing the formulae, and we will use, in a later section,
{\small
$$ P_{-1/2}(x)=\frac{2}{\pi} {\rm EllipticK}\left(
\sqrt{\frac{1-x}{2}}\right), \quad
Q_{-1/2}(x)= {\rm EllipticK}\left(
\sqrt{\frac{1+x}{2}}\right) \quad
{\rm for\ } |x|<1 .
$$
}
Though we don't explicitly use the following it does relate to the weight $w$ in a later section being even,
so we note
$$ \frac{\pi}{2}\,P_{-1/2}(-x)=Q_{-1/2}(x),\qquad
{\rm and\ similarly\ \  } \frac{\pi}{2}\,P_{1/2}(-x)=-Q_{1/2}(x) . $$
The functions with order $\mu=-1$ arise in the case $\alpha=3/2$.
The left-hand parts of the following formulae, corresponding to $\mu=1$ in~\cite{AS}8.2.5,
parallel those given earlier, (\ref{eq:Qm1mh},\ref{eq:Qm1ph}),
\begin{eqnarray}
P^{-1}_{-1/2}(q_2)
&=&   -4\, P^{1}_{-1/2}(q_2)
 ,
 \label{eq:Pm1mh}
\\
P^{-1}_{1/2}(q_2)
&=&   \frac{4}{3} P^{1}_{1/2}(q_2)
 .
 \label{eq:Pm1ph}
\\ 
Q^{-1}_{-1/2}(q_2)
&=&   -4\, Q^{1}_{-1/2}(q_2)
= \frac{4}{\sqrt{2(q_2-1)}}\, 
{\rm EllipticE}\left(\sqrt{\frac{2}{q_2+1}}\right)
 ,
 \label{eq:Qm1mhC}
\\
Q^{-1}_{1/2}(q_2)
&=&   \frac{4}{3} \, 
\frac{-q_2}{\sqrt{2(q_2-1)}}\,
{\rm EllipticE}\left(\sqrt{\frac{2}{q_2+1}}\right)+\\
&\ & \qquad\qquad
\frac{4}{3}\,\sqrt{\frac{q_2-1}{2}}
{\rm EllipticK}\left(\sqrt{\frac{2}{q_2+1}}\right) 
 .
\label{eq:Qm1phC}  
 \end{eqnarray}
 
There are many remarkable formulae
involving both $P$ and $Q$ functions, including Whipple's, which for integer $n$ and $\mu$ is,
 for $z>1$,
  \begin{eqnarray*}
 P_{n-{\frac {1}{2}}}^{\mu}(z )
 &=&{\frac {(-1)^{n}}{\Gamma (n-\mu+{\frac {1}{2}})}}
\sqrt{\frac {2}{\pi}}\, \left(z^2-1\right)^{-1/4}
 Q_{\mu-{\frac {1}{2}}}^{n}(\frac{z}{\sqrt{z^2-1}} )\\
 Q_{n-{\frac {1}{2}}}^{\mu}(z )
 &=& {\frac {(-1)^{n}\pi }{\Gamma (n-\mu+{\frac {1}{2}})}}
 \sqrt{\frac {\pi}{2}}\, \left(z^2-1\right)^{-1/4}
 P_{\mu-{\frac {1}{2}}}^{n}(\frac{z}{\sqrt{z^2-1}} )
  .
\end{eqnarray*}
When both $n$ and $\mu$ are zero, we have
\begin{eqnarray*}
 P_{-{\frac{1}{2}}}(z )
&=&\frac{\sqrt{2}}{\pi}
\left(z^2-1\right)^{-1/4}
Q_{-{\frac{1}{2}}}(\frac{z}{\sqrt{z^2-1}} ),\\
P_{-{\frac{1}{2}}}(\cosh(\eta))
&=& \frac{\sqrt{2}}{\pi}
\left(\sinh(\eta)\right)^{-1/2}  
Q_{-{\frac{1}{2}}}(\coth(\eta)) ,\\
Q_{-{\frac{1}{2}}}(z )
 &=& 
 {\frac{\pi}{\sqrt{2}}}\, \left(z^2-1\right)^{-1/4}
 P_{-{\frac{1}{2}}}(\frac{z}{\sqrt{z^2-1}} ) ,\\
Q_{-{\frac{1}{2}}}(\cosh(\eta))
&=& \frac{\pi}{\sqrt{2}}
\left(\sinh(\eta)\right)^{-1/2}  
P_{-{\frac{1}{2}}}(\coth(\eta))  
  .
\end{eqnarray*}

Another notable formula is from~\cite{HTF2}3.3.1(1)
\begin{equation}
P_{-n-1/2}^\mu(z) = P_{n-1/2}^\mu(z) . 
\label{eq:Pminus}
\end{equation}
The corresponding formula for $Q$ is~(\ref{eq:AS822})
given in~\cite{HTF2}3.3.1(3).

 There are formulae involving Wronskians and similar determinantal forms, e.g. for $z>0$
 \begin{equation}
 P_{1/2}(z) Q_{-1/2}(z)- Q_{1/2}(z) P_{-1/2}(z) = 2 .
 \label{eq:niceW0}
 \end{equation}
On substituting from equations (\ref{eq:PmhDL})--(\ref{eq:QphDL})
one finds that equation~(\ref{eq:niceW0}) is equivalent to
Legendre's relation~(\ref{eq:LegRel}) for the EllipticE and EllipticK functions.
 
\cite{AS}8.1.8 gives, when $\mu=0$,
 $$ {\rm Wronskian}(P_\nu(z),Q_\nu(z))
 =\frac{-1}{z^2-1} .
$$
(See also ~\cite{HTF1}3.4(25).)
 
We also have from
 \cite{HTF1}~3.8 (9)  and~\cite{AS}8.5.2
\begin{equation}
 \frac{d P_{n-1/2}(z)}{dz}
 = \frac{n^2-1/4}{\sqrt{z^2-1}} \, P^{-1}_{n-1/2} .
 \label{eq:HTF1Pp3p8p9}
\end{equation}
Using this and the corresponding formula for the $Q$, namely (\ref{eq:HTF1Qp3p8p9})
$$ P_{n-1/2}(z) Q^{-1}_{n-1/2}(z) - P^{-1}_{n-1/2}(z) Q_{n-1/2}(z)
 = -\frac{1}{(n^2-1/4)\, \sqrt{z^2-1}} . 
$$
(See also~\cite{GST00} equation (7), as in~(\ref{eq:GST007}) below.)

The Wronskian expression and \cite{AS}~8.5.4
\begin{equation}
P_{n+1/2}(z) Q_{n-1/2}(z)- Q_{n+1/2}(z) P_{n-1/2}(z) = \frac{1}{n+1/2} .
\label{eq:niceWn}
\end{equation}
(See also equation~(\ref{eq:niceph}).)
Equation~(\ref{eq:niceWn}) is the $m=0$ case of
equation~(6) of~\cite{GST00}:
\begin{equation}
 P_{n+1/2}^m(z) Q_{n-1/2}^m(z)- Q_{n+1/2}^m(z) P_{n-1/2}^m(z) 
 = \frac{\Gamma(n+\frac{1}{2}+m)}{\Gamma(n+\frac{3}{2}-m)}
 \, (-1)^m .
\label{eq:GST006}
\end{equation}
Equation~(7) of~\cite{GST00} allows for orders to be changed by 1
\begin{equation}
 P_{n-1/2}^m(z) Q_{n-1/2}^{m+1}(z)- Q_{n-1/2}^m(z) P_{n-1/2}^{m+1}(z) 
 = \frac{\Gamma(n+\frac{1}{2}+m)}{\Gamma(n+\frac{1}{2}-m)}
 \, \frac{(-1)^m}{\sqrt{z^2-1}} .
\label{eq:GST007}
\end{equation}

Writing $z$ where formerly we wrote $q_2$, the 
general solution of (RE(1/2)) is
$$ u_n(z) = f_P(z) P_{n-1/2}(z)+f_Q(z) Q_{n-1/2}(z) ,
$$
where $f_P$ and $f_Q$ are independent of $n$.
We can prescribe $u_0$ and $u_1$ and determine $f_P$ and $f_Q$.
Using equation~(\ref{eq:niceW0}) we have
$$ f_P= -\frac{1}{2} (u_0 Q_{1/2}- u_1 Q_{-1/2}),\ 
f_Q= \frac{1}{2} (u_0 P_{1/2}- u_1 P_{-1/2}).
$$
In particular, the sequence starting with $u_0=0$, $u_1=1$,
and denoted by $p_{01}(n)$ in the next subsection is
\begin{equation}
p_{01}(n)=\frac{1}{2} (Q_{-1/2} P_{n-1/2} -P_{-1/2} Q_{n-1/2} ).
\label{eq:p01PQ}
\end{equation}
Similarly
\begin{equation}
p_{10}(n)= -\frac{1}{2} (Q_{1/2} P_{n-1/2} -P_{1/2} Q_{n-1/2} ).
\label{eq:p10PQ}
\end{equation}
While these must be the polynomial sequences tabulated in the 
next subsection, it isn't immediately obvious how to extract 
information about the polynomials from the formula above.
(For general $\alpha$, see~\cite{BI82} equation (3.5).)
One item that follows readily  is that the above 
with equations~(\ref{eq:AS822}) and~(\ref{eq:Pminus})
imply $p_{01}(-n)=p_{01}(n)$, etc.

\subsubsection{Reduction of order}

We have already noted that reduction of order can be used
when $x=q_2=1$ producing from the obvious constant solution to
(RE($\alpha$)) a second linearly independent solution.
This subsection fails to help so far in our quest to find more information on the polynomial solutions.
However reduction of order gives equation~(\ref{eq:wSol})and 
this can be used to give a relation between toroidal $P$ and $Q$ functions.

Assume $u_n$ is one solution of (RE($\alpha$)) and seek
a second linearly independent solution $v_n$ with
$$ v_n = w_n \, u_n\qquad {\rm with\ \ } w_1\ne w_0 
\quad {\rm and\ with\ all\ } u_n\ne{0-}{\rm function}.$$
Writing out the equation stating that $w_n u_n$ satisfies
(RE($\alpha$)) and then eliminating $u_n$ from the fact that
$u_n$ also satisfies (RE($\alpha$)) gives
$$(n+\alpha) u_{n+1} (w_{n+1}-w_n)+
(n-\alpha) u_{n-1} (w_{n-1}-w_n)= 0. $$
Let $\theta_n=w_{n+1}-w_n$.
Then
$$ (n+\alpha) u_{n+1}\theta_n -(n-\alpha) u_{n-1}\theta_{n-1} = 0. $$
Hence
$$ \frac{\theta_n}{\theta_{n-1}}
=\frac{(n-\alpha)u_{n-1}}{(n+\alpha)u_{n+1}} , $$
and hence
$$ \frac{\theta_n}{\theta_{n-2}}
=\frac{(n-\alpha)(n-1-\alpha) u_{n-1} u_{n-2}}
{(n+\alpha)(n-1+\alpha) u_{n+1} u_n} . $$
Now set $\alpha=1/2$. The last equation becomes
$$ \frac{\theta_n}{\theta_{n-2}}
=\frac{(n-2+\frac{1}{2}) u_{n-1} u_{n-2}}
{(n+\frac{1}{2}) u_{n+1} u_n} . $$
Continuing to multiply the ratios of $\theta$ we get
$$ \frac{\theta_n}{\theta_{0}}
=\frac{(\frac{1}{2}) u_{1} u_{0}}
{(n+\frac{1}{2}) u_{n+1} u_n} . $$
Our concern is with $w$ rather than $\theta$ and
$$ w_{n+1}-w_0 = \sum_{k=0}^n \theta_k .$$
Now $\theta_0$ is known and nonzero as
$\theta_0=w_1-w_0$.
Hence
\begin{equation}
w_{n+1}
= w_0 + \theta_0 u_0 u_1
\sum_{k=0}^n \frac{1}{(2k+1) u_{k+1} u_k} .
\label{eq:wSol}
\end{equation}

The preceding equation concludes the general reduction of order.
However, a quick check is to return to the case $x=q_2=1$
and the known solution $u_k=1$ for all $k$.
The reduction of order formula for $w_{n+1}$ involves
the sum of the reciprocals of odd integers, and this
is consistent with our previous formula involving
the digamma function $\psi$.

\medskip
We remark that one could view one instance as $w_n=P_{n-1/2}/Q_{n-1/2}$
and then one has an identity involving Legendre functions.
It presumably is a much disguised consequence of the 
recurrence relations (RE(1/2)) satisfied by the $P$ and $Q$.

\medskip
We remark also that if the input $u_n$ to the reduction of order is
polynomial in $x=q_2$, the output $v_n$ is rational but not necessarily polynomial.\\

\medskip
For a formula similar to equation~(\ref{eq:wSol})
see~\cite{Gr85} equation (14).

\subsection{Polynomial sequences relating to LegendreQ functions}\label{subapp:PolyQ}
Clearly we can write
\begin{equation}
Q_{n-1/2}(z) = p_{10}(\alpha=1/2,n)\, Q_{-1/2}(z) + p_{01}(\alpha=1/2,n)\, Q_{1/2}(z) ,
\label{eq:Q1pol}
\end{equation}
where the polynomials $p$ satisfy (RE(1/2)) and
{\small
$$  p_{10}(\alpha=1/2,0)=1,\  p_{10}(\alpha=1/2,1)=0,\quad
p_{01}(\alpha=1/2,0)=0,\  p_{01}(\alpha=1/2,1)=1 . $$
}
When space is too limited, we drop from the argument list
the $\alpha=1/2$, leaving $n$ as the sole argument.

\begin{center}
\begin{tabular}{|c|c|c|}
\hline
$n$& $p_{01}(n)$ & $p_{10}(n)$\\
\hline
$0$&$0$&$1$\\
$1$&$1$&$0$\\
$2$&$\frac{4}{3}\,x$&$-\frac{1}{3}$\\
$3$&${\frac {32}{15}}\,{x}^{2}-\frac{3}{5}$&$-{\frac {8}{15}}\,x$\\
$4$&${\frac {128}{35}}\,{x}^{3}-{\frac {208}{105}}\,x$&
$-{\frac {32}{35}}\,{x}^{2}+{\frac {5}{21}}$\\
$5$&${\frac {2048}{315}}\,{x}^{4}-{\frac {544}{105}}\,{x}^{2}
+{\frac {7}{15}}$&
$-{\frac {512}{315}}\,{x}^{3}+{\frac {88}{105}}\,x$\\
\hline
\end{tabular}
\end{center}

\begin{figure}[hb]
\begin{center}
\includegraphics[height=6cm,width=8cm]{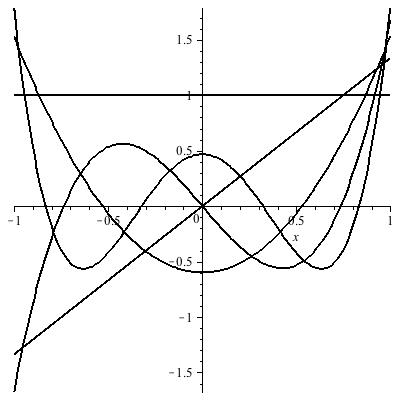} \\
\end{center}
\caption{Plots of $p_{01}(n)$ for $n=1$ to 5.}
\label{fig:p01plot}
\end{figure}

Plots of the $p_{01}$ polynomials are given in Figure~\ref{fig:p01plot}.

\goodbreak

{\small
\cmdvmcode{Recurrence/iter1906mpl.txt}{ }
}

Similarly, we write
\begin{equation}
Q_{n-1/2}^{-1}(z) = {\hat p}_{10}(\alpha=1/2,n)\, Q_{-1/2}(z) + 
{\hat p}_{01}(\alpha=1/2,n)\, Q_{1/2}(z)  ,
\label{eq:Q3pol}
\end{equation}
where the polynomials ${\hat p}$ satisfy (RE(3/2)) and
$$  {\hat p}_{10}(\alpha=3/2,0)=1,\  
{\hat p}_{10}(\alpha=3/2,1)=0,\qquad
{\hat p}_{01}(\alpha=3/2,0)=0,\  
{\hat p}_{01}(\alpha=3/2,1)=1 . $$
Thus, $p_{01}$ and $p_{10}$ are the polynomials corresponding to $\alpha=1/2$ and
${\hat p}_{01}$ and ${\hat p}_{10}$ are those corresponding to 
$\alpha=3/2$.

\medskip
\begin{center}
\begin{tabular}{|c|c|c|}
\hline
$n$& ${\hat p}_{01}(n)$ & ${\hat p}_{10}(n)$\\
\hline 
%
$0$&$0$&$1$\\
$1$&$1$&$0$\\
$2$&$\frac{4}{5}\,x$&$\frac{1}{5}$\\
$3$&${\frac {32}{35}}\,{x}^{2}-\frac{1}{7}$&${\frac {8}{35}}\,x$\\
$4$&${\frac {128}{105}}\,{x}^{3}-{\frac {16}{35}}\,x$&
${\frac {32}{105}}\,{x}^{2}-{\frac {1}{15}}$\\
$5$&${\frac {2048}{1155}}\,{x}^{4}-{\frac {416}{385}}\,{x}^{2}
+{\frac {5}{77}}$&
${\frac {512}{1155}}\,{x}^{3}-{\frac {232}{1155}}\,x$\\
\hline
\end{tabular}
\end{center}

We have seen in earlier sections that there are many connections
between the polynomials at $\alpha=1/2$ and $\alpha=3/2$.
For example, defining
$$ p_{1x}(n) = p_{10}(n) + x\, p_{01}(n) ,$$
Theorem~\ref{thm:thm1} suggests relations with the ${\hat p}$.
We have
$$ \frac{p_{1x}(n+1)-p_{1x}(n-1)}{x^2-1}
= \frac{3\, {\hat p}_{01}(n)}{4\, n} . $$
This also follows from equations~(\ref{eq:L1ppBarn})~(\ref{eq:L2ppBarn}).)

\medskip

Various quantities,
whose notation has been standardised, as in~\cite{HTF2}\S{10},
are associated with sequences of orthogonal polynomials.
These include the leading coefficient, often denoted $k_n$,
moments, often denoted $c_n$ and discussed in the later
subsection with the weight function,
and the norm squared of the polynomial, denoted $h_n$.
(The quantities $k_n$ and $h_n$ occur in~\cite{HTF2}10.3(11)
in connection with a Christoffel-Darboux formula,
a formula we will present in a later subsection.)
For the family $p_{01}$ we have
\begin{equation}
{\rm lcoeff}(p_{01}(n)) = k_n 
= \frac{2^{2n-2} (n-1)!}{(2n-1)!!} . 
\label{eq:lc01}
\end{equation}
The leading coefficients of the other sequences of polynomials
satisfy, for $n\ge{2}$,
\begin{eqnarray*}
{\rm lcoeff}(p_{01}(n))
&=& -4\, {\rm lcoeff}(p_{10}(n)) ,\\
{\rm lcoeff}({\hat p}_{01}(n))
&=& 4\, {\rm lcoeff}({\hat p}_{10}(n)) ,\\
{\rm lcoeff}({\hat p}_{01}(n))
&=& \frac{3}{2n+1}\, {\rm lcoeff}(p_{01}(n)) ,\\
{\rm lcoeff}({\hat p}_{10}(n))
&=& -\frac{3}{2n+1} {\rm lcoeff}(p_{10}(n)) .
\end{eqnarray*}
\cite{HTF2}10.3(7,8) writes the recurrence relation in the
form
$$ p_{n+1}= A_n x\, p_n - C_n p_{n-1} , $$
and notes that
$$ A_n = \frac{k_{n+1}}{k_n}, \qquad 
C_n= \frac{k_{n+1}\, k_{n-1}}{k_n^2}\,\frac{h_n}{h_{n-1}} . $$
For us at $\alpha=1/2$,
$$ A_n=\frac{4n}{2n+1},\ 
\frac{h_n}{h_m}=\frac{m}{n},\ C_n=\frac{2n-1}{2n+1} . $$

\medskip
Our interest is in $g_n$ and ${\hat g}_n$, so
$\alpha=1/2$ and $\alpha=3/2$.
There remains the possibility that investigating other values
of $\alpha$ might ultimately be useful, and, perhaps
lead to the discovery of further `contiguous relations'
such as those in $\mu$ for the Legendre functions.
In this connection we record some results on $\alpha=-1/2$.

\begin{center}
\begin{tabular}{|c|c|c|}
\hline
$n$& ${ p}_{01}(\alpha=-\frac{1}{2},n)$ & 
${ p}_{10}(\alpha=-\frac{1}{2},n)$\\
\hline 
$0$&$0$&$1$\\
$1$&$1$&$0$\\
$2$&$4\,x$&$-3$\\
$3$&${\frac {32}{3}}\,{x}^{2}-\frac{5}{3}$&$-8\,x$\\
$4$&${\frac {128}{5}}\,{x}^{3}-{\frac {48}{5}}\,x$&
$-{\frac {96}{5}}\,{x}^{2}+{\frac {21}{5}}$\\
$5$&${\frac {2048}{35}}\,{x}^{4}-{\frac {1248}{35}}\,{x}^{2}
+{\frac {15}{7}}$&
${-\frac{1536}{35}}\,{x}^{3}+{\frac {696}{35}}\,x$\\
\hline
\end{tabular}
\end{center}

We focus on connections between the solutions of
(RE(-1/2)) and those of (RE(1/2)) and (RE(3/2)).
As noted in Theorem~\ref{thm:thm3a},
$$ p_{01}(\alpha=-1/2,n)
= \frac{4 n^2 -1}{3}\, {\hat p}_{01}(n) . $$

\subsubsection{Derivatives of $p_{01}$ and of $p_{10}$}

\noindent{\bf Ordinary differential equations}
\smallskip

\noindent
The Legendre differential equation satisfied by $ Q_{n-1/2}(z)$ leads to coupled d.e.s
for  $p_{01}(n)$ and $p_{10}(n)$ (where we drop the reminder $\alpha=1/2$).
{\small
\begin{eqnarray}
 (1-x^2) {\frac {d^{2}p_{10}(n)}{d{x}^{2}}}-
 x{\frac {d p_{10}(n)}{dx}}+
 p_{10}(n){n}^{2}   +
{\frac {d p_{01}(n)}{dx}}
&=& 0 \label{eq:p10de} ,\\
(1-x^2)  {\frac {d^{2}p_{01}(n)}{d{x}^{2}}}-
3\, x\,  {\frac {d p_{01}(n)}{dx}} +
({n}^{2}-1)\,p_{01}(n)-
 {\frac {d p_{10}(n)}{dx}}
&=& 0 . \label{eq:p01de}
\end{eqnarray}
}

From these one can find 4-th order differential equations
satisfied by the polynomials.
For example $p_{01}(n)$, satisfies, with
\begin{eqnarray*}
P
&=& (1-x^2)^2
 , \\
Q
&=& 10 x (1-x^2)
 , \\
R
&=& 2n^2 (1-x^2)-9 +25 x^2
 , \\ 
S
&=& -3 x (2n^2 -5)
 , 
\end{eqnarray*}
the differential equation
$$ P\frac{d^4 p_{01}(n)}{d x^4} +
Q\frac{d^3 p_{01}(n)}{d x^3} +
R\frac{d^2 p_{01}(n)}{d x^2} +
S\frac{d p_{01}(n)}{d x} +
(n^2-1)^2 p_{01}(n) = 0. $$
The notation is as in~\cite{Me66}.
The conditions in the theorems in~\cite{Me66} are 
{\it not} satisfied, but some equations similar to his are.
We have, with $W_M=\sqrt{1-x^2}$ which is {\it not}
the weight function for us,
\begin{eqnarray*}
W_M Q &=& 2\frac{d}{dx}\left(W_M P\right) ,\\
W_M S 
&=&\frac{d}{dx}\left(W_M(R-1)\right)
- \frac{d^3}{dx^3}\left(W_M P\right) . 
\end{eqnarray*}

The remainder of this subsection is speculative
(and won't be included in any paper).
We also don't see any use for the main goals of the paper.
We have yet to write the d.e. in the form
$$ {\cal L}_4 u = f(x,n) {\cal L}_2 u , $$
with, in some appropriate sense, 
the 4th order ${\cal L}_4$ and the 2nd order ${\cal L}_2$
being somehow `self-adjoint'.
Associated with this effort, we would like to be able to
express derivatives of the $w(x)$ defined in
equation~(\ref{eq:weightK}) in terms of 
a rational function of $w(x)$, $x$ and $\sqrt{1-x^2}$
(or something similarly simple).

We believe there is much more to be discovered,
e.g. in connection with variational characterisations of
the (generalized eigen)solutions of the fourth order operator.
We expect that the qualitative behaviour would be akin
to the generalized matrix eigenvalue problem
$$ A u = \lambda B u ,$$
with matrices $A$ and $B$ self-adjoint (real-symmetric).
The function spaces involved might include the space of
all polynomials.
Having found (eigen)solutions polynomial of degree $n$,
Rayleigh Ritz might be appropriate to produce solutions of degree $n+1$.

Different inner products may be appropriate in different circumstances.
\medskip

\noindent{\bf The Wronskian of $p_{01}$ and $p_{10}$}
\smallskip

\noindent
We find, at any fixed $n$, using equations~(\ref{eq:dp10}),~(\ref{eq:dp01})
to eliminate the derivative terms and then equation~(\ref{eq:niceph})
to eliminate the polynomials indexed by $n-1$,
\begin{eqnarray}
{\rm Wronskian}(p_{01},p_{10})
&=&  (p_{01}\, p_{10}' - p_{10}\, p_{01}') ,
\nonumber \\
&=& \frac{(p_{01}+ p_{10})^2 +2\, (x-1)p_{01}\, p_{10} -1}
{2(x^2-1)} .
\label{eq:Wn}
\end{eqnarray}

\subsubsection{Contiguous relations}

We have already seen several instances of `contiguous relations'.
\begin{itemize}
\item
Some of these involve only the polynomials at a fixed value of
$\alpha$ as in (RE($\alpha$)) and, in this category, we
mention~(\ref{eq:niceph}) and~(\ref{eq:niceph3}) below.

\item
Some of these involve polynomials at a fixed $n$ and
the different sequences $p$ at $\alpha=1/2$ and
$\hat p$ at $\alpha=3/2$.
One such, namely~(\ref{eq:L1ppBarn}),~(\ref{eq:L2ppBarn}),
was a simple consequence of Theorem~\ref{thm:thm1}.

\item
Others involve not merely the polynomials, but also their
first derivatives.

\end{itemize}

We now see an  identity (\ref{eq:niceph})  for the polynomial solutions of (RE(1/2))
which is similar to (\ref{eq:niceWn}) applying to LegendreQ.
At $\alpha=1/2$ we have
\begin{equation}
 p_{10}(n)\, p_{01}(n+1)- p_{10}(n+1)\,  p_{01}(n)
= \frac{1}{2n+1} .
\label{eq:niceph}
\end{equation}
(See also the preceding subsection on  alternative formulations of the recurrence relation,
specifically (RE1x$a_n$).)
At $\alpha=3/2$ we have
\begin{equation}
 {\hat p}_{10}(n)\, {\hat p}_{01}(n+1)-
{\hat p}_{10}(n+1)\,  {\hat p}_{01}(n)
= \frac{-3}{(2n-1)(2n+1)(2n+3)} .
\label{eq:niceph3}
\end{equation}
(There is a similar result for $\alpha=0$:
$$ T_n(x)\, U_{n+1}(x) - T_{n+1}(x)\, U_n(x) = 1. $$
And, again, another similar result for $\alpha=1$.)

These are easily proved by induction, as given in the following
generalization which applies not just to polynomial solutions of
the difference equations.

\begin{theorem}\label{thm:induct}
Let $u_n$ and $U_n$ be two sequences satisfying {\rm (RE($\alpha$))}.
Then
$$ u_n\, U_{n+1} - u_{n+1}\, U_n
= c(x)\, \frac{\Gamma(n+1-\alpha)}{\Gamma(n+1+\alpha)} , $$
where
$$ c(x)
= \left(u_0\, U_{1} - u_{1}\, U_0\right)\,
\frac{\Gamma(1+\alpha)}{\Gamma(1-\alpha)} .
$$
\end{theorem}

\medskip
\noindent{\it Proof.}
The proof is by induction.
Define
$${\cal S}(n) 
= u_n\, U_{n+1} - u_{n+1}\, U_n
- c(x)\, \frac{\Gamma(n+1-\alpha)}{\Gamma(n+1+\alpha)} ,
$$
with $c(x)$ as in the statement of the theorem.
We have, from the definition of $c(x)$, that ${\cal S}(0)=0$.
We will next show that ${\cal S}(n-1)=0$ implies
that ${\cal S}(n)=0$.
Let ${\cal S}^{*}(n)$ be the expression obtained from
 ${\cal S}(n)$ on using (RE($\alpha$)) 
to write $u_{n+1}$ in terms of $u_n$ and $u_{n-1}$ and, similarly,
to write $U_{n+1}$ in terms of $U_n$ and $U_{n-1}$.
One finds
$${\cal S}^{*}(n)
= \frac{n-\alpha}{n+\alpha}\, {\cal S}(n-1) .$$
Thus ${\cal S}(n-1)=0$ implies ${\cal S}(n)=0$
so, as asserted, the theorem is established by induction.

\bigskip

Formulae involving the derivatives of the polynomials can also be found.

The starting point in deriving~(\ref{eq:dp10}) and~(\ref{eq:dp01})
below is
$$ (x^2-1)\, \frac{d Q_\nu}{d x} = \nu\, (x Q_\nu - Q_{\nu-1}) .$$
See \cite{HTF1}~3.8.10. Also
\begin{eqnarray*}
 (x^2-1)\, \frac{d Q_{1/2}}{d x} 
 &=& \frac{1}{2}\, (x Q_{1/2}-Q_{-1/2}) ,\\
(x^2-1)\, \frac{d Q_{-1/2}}{d x} 
 &=& -\frac{1}{2}\, (x Q_{-1/2}-Q_{1/2}) .
\end{eqnarray*}
This leads to the somewhat lengthy identities
{\small
\begin{eqnarray}
 2(x^2-1)\, \frac{d p_{10}(j)}{d x} 
 &=& 2 j x p_{10}(j)+p_{01}(j)-(2j-1)p_{10}(j-1) ,
 \label{eq:dp10}\\
2(x^2-1)\, \frac{d p_{01}(j)}{d x} 
 &=&  (2j-2)x p_{01}(j)-p_{10}(j)-(2j-1)p_{01}(j-1) .
 \label{eq:dp01}
\end{eqnarray}
}
 
There are also formulae involving the four sequences of polynomials,
both the $p$ sequences and the $\hat p$ sequences.
The starting point for equations~(\ref{eq:pB01n}) and~(\ref{eq:pB10n})
below is that derivatives with respect to $z$ of 
$Q_{n-1/2}^\mu(z)$ can be expressed in terms of the
undifferentiated function at a different value of $\mu$.
In particular, begin with equation~(\ref{eq:HTF1Qp3p8p9}).
After some calculation one finds
\begin{eqnarray}
(4n^2-1)\,{\hat p}_{01}(n)
&=&
6\frac{d}{dx}\left(x\,p_{01}(n)+p_{10}(n)\right) - 3p_{01}(n) ,
\label{eq:pB01n}
\\
(4n^2-1)\,{\hat p}_{10}(n) ,
&=&
2\frac{d}{dx}\left(x\,p_{10}(n)+p_{01}(n)\right) - 3p_{10}(n) .
\label{eq:pB10n}
\end{eqnarray}

\subsubsection{Christoffel-Darboux (CD)}

As in~\cite{HTF2}10.3(11) we have
$$ \sum_{j=1}^n j p_{01}(j)^2
= \frac{2n+1}{4}\, \left(
p_{01}(n)\, \frac{d p_{01}(n+1)}{d x} -
\frac{d p_{01}(n)}{d x} \, p_{01}(n+1) \right) .$$

The CD kernel, defined by
$$ K(x,y) = \sum_{j=1}^n j p_{01}(j,x)\, p_{01}(j,y) , $$
has the reproducing property
$$ K(x,y)= \langle K(x,\cdot),K(\cdot,y) \rangle .$$

\subsubsection{A convolution identity involving Legendre polynomials}

A special case of~\cite{BI82} equation (3.13) is
\begin{equation}
p_{10}(n)
= -\sum_{k=0}^n\frac{{\rm LegendreP}(k,x)\, {\rm LegendreP}(n-k,x)}{2k-1} .
\label{eq:p10Leg}
\end{equation}
Another identity (also in~\cite{Ch78}~(12.6) is
\begin{equation}
p_{01}(n)
= -\sum_{k=0}^{n-1}\frac{{\rm LegendreP}(k,x)\, {\rm LegendreP}(n-1-k,x)}{2k+1} .
\label{eq:p01Leg}
\end{equation}
(Equation~(\ref{eq:p01Leg}) can also be found using the generating function
given in (\ref{eq:GenDef}).)

We have yet to find similar concise identities for $\alpha=3/2$ and the 
$\hat p$ polynomials.
Of course, using equations~(\ref{eq:p01Leg}),~(\ref{eq:p10Leg}) in equation~(\ref{eq:L2ppBarn})
gives the $4x{\hat p}(n)-12{\hat p}(n)$
in terms of Legendre polynomials.

Presumably related to $\alpha=0$ or $\alpha=1$ we have the well-known
identity (\cite{Ch78}~(12.7)):
$$ {\rm ChebyshevU}(n,x)
= \sum_{k=0}^n {\rm LegendreP}(k,x)\, {\rm LegendreP}(n-k,x) .
$$
We have yet to find a similar identity for ChebyshevT.

There are, of course, many other ways the identities can be written
as, if $f$ is any real-valued function odd about 0
(and any functions of $(k,x)$ not just LegendreP),
$$ 0 
=  \sum_{k=0}^n f(k-\frac{n}{2}) {\rm LegendreP}(k,x)\, {\rm LegendreP}(n-k,x) .
$$
A typical example would be $f(k-\frac{n}{2})=k-\frac{n}{2}$.

Allowing other Gegenbauer polynomials than the Legendre $P$
polynomials above, we remark that one obtains convolution
identities on using 
$\gamma=-\beta=\alpha-1$
in~\cite{BI82} equation (3.13).

\subsection{$\alpha=1/2$: orthogonality and related results}\label{subapp:orthog}

Results from~\cite{BD67} are reported in~\cite{Ch78} and~\cite{Gr85}.
Equation~(12.3) of~\cite{Ch78} and equation(1) of~\cite{Gr85} defines
for the interval $(-1,1)$, a weight function $w$ 
(in the case $\nu=1/2$ appropriate to our problem) by
\begin{equation}
w(x) = \frac{4}{Q_{-1/2}(x)^2 +(\pi P_{-1/2}(x)/2)^2} .
\label{eq:wPQ}
\end{equation}
Using representations of the toroidal functions in terms of
complete elliptic integrals, EllipticK, this rewrites to
\begin{equation}
w(x) = \frac{4}{K\left(\sqrt{\frac{1-x}{2}}\right)^2+
K\left(\sqrt{\frac{1+x}{2}}\right)^2} . 
\label{eq:weightK}
\end{equation}
(See also~\cite{VZ07} where,
with $A_n$ satisfying their (1.5), 
$k^n A_n$ satisfies (RE(1/2)).
Also their $\psi_n$ satisfying their (4.10)
satisfy our (RE(1/2)), except for a shift.
Weight functions involving EllipticK occur in their (3.7),
(3.9), (5.19) and, close to the form we report in
(5.20) and (5.21).
See also~\cite{Rees45}.)
We have yet to investigate (RE(1/2)) in connection with study of examples of
`polynomials orthogonal on the (boundary of the) unit disk':
see equation~(\ref{eq:Kdisk}).

Define the inner product
$$\langle u,v\rangle
= \int_{-1}^1 u(x)\, v(x)\, w(x)\, dx . $$
Direct calculation (via wolframalpha) for several $(m,n)$ pairs confirms
\begin{equation}
\langle p_{01}(m),p_{01}(n)\rangle
= \frac{\delta_{mn}}{n} . 
\label{eq:p01ort}
\end{equation}
Shortly we will give an elementary derivation of this without reference to
the details of the inner product.

The moments are defined by
$$ c_n = \langle 1, x^n\rangle ,$$
and are obviously 0 when $n$ is odd. We have
$$ c_0=1, \ c_2=\frac{9}{2^5},\ c_4=\frac{39}{2^8},\ 
c_6=\frac{6633}{2^{16}} , $$
but have yet to find the formula for $c_{2n}$.
Once all the $c_n$ are calculated,
\cite{HTF2}~10.3(4) gives a determinantal formula for the
orthogonal polynomials.

\medskip
Our derivation of equation~(\ref{eq:p01ort}) involves considering
$$\phi(n)= \langle x^{n-1}, p_{01}(n) \rangle.$$
(See also~\cite{Ch78} Theorems 3.2 and 4.2.)
Multiplying the recurrence
$$ 2n x p_{01}(n) 
= (n+\frac{1}{2})\, p_{01}(n+1) + (n-\frac{1}{2}) \, p_{01}(n-1) , $$
by $x^{n-2}$ and taking the inner product gives
$$ 2 n \phi(n) = (n-\frac{1}{2})\, \phi(n-1) .$$
If one agrees that $\phi(1)=\langle 1,1\rangle =1$,
this first order recurrence for $\phi$ solves to
$$ \phi(n)=\frac{(2n-1)!!}{4^{n-1}\, n!}
=\frac{1}{n\, {\rm lcoeff}(p_{01}(n))} ,$$ 
where the leading coefficient lcoeff is given by equation~(\ref{eq:lc01}).
Next consider
$\langle p_{01}(n), p_{01}(n) \rangle$
beginning with considering
$$p_{01}(n)
={\rm lcoeff}(p_{01}(n))\, x^{n-1} +\ {\rm lower\ order} , $$
and noting that the inner product of the lower order monomials with $p_{01}(n)$ is 0.
Thus
$$\langle p_{01}(n), p_{01}(n) \rangle
= {\rm lcoeff}(p_{01}(n)) \, \langle x^{n-1}, p_{01}(n) \rangle
= {\rm lcoeff}(p_{01}(n)) \, \phi(n) . $$
From the preceding calculated value of $\phi(n)$ we find
the value stated in~(\ref{eq:p01ort}).

\subsubsection{$\alpha$ generally}

The orthogonality result is also available, in a
more general form in~\cite{HTF2},p220,\S10.21
where one of the Pollaczek sequences is treated.
The weight function, for our $\alpha=1/2$ case is given 
in~\cite{HTF2}10.21(14) when one sets
$$a=b=0=t, \qquad c=\lambda= \frac{1}{2} . $$
These values in (13) yield our (RE(1/2)) while, in (14),
the weight function involves, as before, EllipticK functions
as the hypergeometric expression in (14) is
$$ {}_2{F}_1(\frac{1}{2},\frac{1}{2};1;z)
=\frac{2}{\pi} K(\sqrt{z}) . $$
For $\alpha=3/2$ we have,in both (13) and (14),
$$a=b=0=t, \qquad c=\frac{3}{2},\ \lambda= -\frac{1}{2} . $$
This time the weight function given in (14) is
$ {}_2{F}_1(\frac{3}{2},\frac{3}{2};1;z)$ which,
again, can be expressed in terms of EllipticE and EllipticK.
We comment further on the $\alpha=3/2$ case below.

\subsubsection{More concerning $\alpha=3/2$}

We noted, in connection with Pollaczek polynomials above,
that the weight function for our $\alpha=3/2$ case is available.
However, the setting is considerably more general than
appropriate for our very special cases.

Gegenbauer polynomials (themselves special cases of
Jacobi polynomials) satisfy the recurrence
$$ (n+1)\, C_{n+1}^\lambda{(x)}
= 2(n+\lambda)x\, C_n^\lambda{(x)}-
(n+2\lambda-1) C_{n-1}^\lambda{(x)} . 
\eqno({\rm REG})$$
\begin{itemize}
\item
Our (RE(1)), satisfied by scaled Chebyshev polynomials, 
is the $\lambda=0$ case of this.
\item  (RE(0)) is the $\lambda=1$ case of (REG).
\end{itemize}
However (REG) does not include the (RE($\alpha$)) needed in
our paper.
The case $\lambda=1/2$ in (REG) gives Legendre polynomials.

\cite{BI82}\S3, equation (3.1) modify this to
$$ (n+\lambda+1)\, C_{n+1}^\lambda{(x,\beta)}
= 2x\,(n+\beta+\lambda)x\, C_n^\lambda{(x,\beta)}-
(2\beta + n+\lambda-1) C_{n-1}^\lambda{(x,\beta)} . 
\eqno({\rm BI3.2})$$
With $\lambda=\alpha-1=-\beta$ we have our (RE($\alpha$)).
In particular
\begin{itemize}
\item $\lambda=\beta=0$ is our (RE(1)).
\item $\lambda=-1=-\beta$ is our (RE(0)).
\item $\lambda=-1/2=-\beta$ is our (RE(1/2)).
\item $\lambda=1/2=-\beta$ is our (RE(3/2)).
\end{itemize}

\subsubsection{An aside on other cases in~\cite{Gr85}}

This is an aside as the only case which is amongst our (RE($\alpha$)) is the case $\alpha=1/2$.

The recurrence relation at the beginning of~\cite{Gr85}
is, with a shift so as to accord with our indexing,
$$ (n+\nu) p_{n+1}
= (2n+2\nu-1)x\, p_n - (n+\nu-1)\, p_{n-1} .$$
We have already used that when $\nu=1/2$ this is our (RE(1/2)).
The weight function for orthogonality is
$$ w(\nu,x)
=\frac{1}{(\nu Q_{\nu-1}(x))^2 +(2\nu P_{\nu-1}(x)/\pi)^2 } .
$$

Now, the limit of $w(\nu,x)$ as $\nu$ tends to 0 is 1,
so we will write
$$ w(0,x) =1 .$$
The recurrence is solved by the Legendre polynomials
and these are orthogonal with respect to the weight $w(0,x)$.

We also have
\begin{eqnarray*}
w(1,x)
&=& \frac{1}{Q_0(x)^2 + (\pi P_0(x)/2)^2} ,\\
&=& \frac{1}{\left(\frac{1}{2}\log|\frac{1+x}{1-x}|\right)^2 +
\frac{\pi^2}{4}} .
\end{eqnarray*}
Repeating work in~\cite{Gr85}, the usual notation
for the LegendreQ functions is
$$ Q_n(x) = P_n(x) Q_0(x) - W_{n-1}(x) .$$
In~\cite{Gr85} there is the demonstration that the
polynomials $W_n$ are orthogonal with respect to the weight
$w(1,x)$.

\subsection{Bessel functions}\label{subapp:Bessel}

In view of one of our goals being to sum series such as that for $u_\infty$
there is a case for investigating representations of the toroidal functions.

Connections with Bessel functions are available.
See~\cite{Sonn59} where the reference is to Watson's classic book on Bessel functions.
See also~\cite{CV12} equation (3) and references there.
The conditions on the parameters are not those we need here.
The formulae below are sometimes approached via Laplace transforms.
We have, for $q>1$,
\begin{equation}
\sqrt{q}\int_0^\infty \exp(-t\sqrt{q})\, {\rm BesselI}(0,t\cos(\psi))\, dt
= g(\psi) 
=\frac{1}{\sqrt{1- \frac{\cos(\psi)^2}{q}}} .
\label{eq:besselg}
\end{equation}
The Fourier coefficients of $g$, $Q_{n-1/2}(q_2)$ can also be written in terms of integrals of Bessel functions.
We have, for $s>1$,
$$\int_0^\infty \exp(-s\, t)\, {\rm BesselI}(n,\frac{t}{2})^2\, dt
= \frac{2}{\pi} Q_{n-1/2}(2 s^2- 1) ,$$
and its application here would have $s=\sqrt{q}$.
See~\cite{GR} 6.612.3.

\subsection{Boundary moments revisited}\label{subapp:BndryMom}

We remark that the perimeter can be expressed in terms of
a toroidal function: see equation~(\ref{eq:ab09P}).

Having introduced the $g_n$ and ${\hat g}_n$, we might note
that boundary moments can be expressed in terms of them.
For example, as
$$ r^2= x^2+ y^2
= a^2 (1-e^2)\left( 1-\frac{e^2}{2} + \frac{e^2}{2} \cos(2\psi)
\right) ,
$$
we have
\begin{eqnarray*}
i_2
&=& a \int_{-\pi}^\pi r^2\, {\hat g}(\psi)\, d\psi\\
&=& a^3 (1-e^2)\left(
( 1-\frac{e^2}{2})\frac{{\hat g}_0}{2} +
\frac{e^2}{2} {\hat g}_1
\right) .
\end{eqnarray*}
There may be useful recurrence relations involving the $i_n$.

It may be that there are neat expressions for any moment,
not merely the `polar' moments of the preceding paragraph,
involving toroidal functions, both of the LegendreP and
LegendreQ kinds.

\clearpage

\newpage

\begin{center}
{\Large{\textbf{\textsc{ Part III:
A general setting for variational approximation\\
$R\le{Q}$\\
 \qquad  \\
By G. Keady
}}}}
\end{center}
\section{A general lower bound} 

A published form of the important results of this part is~\cite{KW20}.
Further applications are given in~\cite{Ke20i,Ke20}

\subsection{A development from \cite{KM93}}\label{subsec:inS0Sinf}

Our notation is as in~\cite{KM93}.
For large $\beta$
$$ u \sim \beta\, \frac{|\Omega|}{|\partial\Omega|}  +  u_\infty + o(1) , $$
where $u_\infty$ solves
$$ - \frac{\partial^2 u_\infty}{\partial x^2}- \frac{\partial^2 u_\infty}{\partial y^2}
 = 1\, {\rm in}\,\Omega ,\,
\frac{\partial u_\infty}{\partial n}
= -\frac{|\Omega|}{|\partial\Omega|} \, {\rm on}\,\partial\Omega \,
{\rm and}\,  \int_{\partial\Omega} u_\infty = 0 .
\eqno({\rm P}(\infty))
$$
Define
\begin{equation}
\Sigma_\infty
= \int_\Omega u_\infty , \qquad  {\rm and\ \ } \Sigma_1= - \int_{\partial\Omega} u_\infty^2 ,
\label{eq:SigDef}
\end{equation}
with $u_\infty$ satisfying Problem (P($\infty$)).
The terms $\Sigma_\infty>0$ and $\Sigma_1\le{0}$ (both obviously independent of $\beta$) give coefficients in
the next terms in the asymptotic expansion
$$ \QQ(\beta) \sim \frac{\beta\, |\Omega|^2}{|\partial\Omega|} + \Sigma_\infty +\frac{\Sigma_1}{\beta}
\qquad {\rm for\ \ } \beta\rightarrow\infty . $$
Inequality~(4.10) of~\cite{KM93} (also (4.12)) states
\begin{equation}
 \frac{\beta\, |\Omega|^2}{|\partial\Omega|} + \Sigma_\infty +\frac{\Sigma_1}{\beta} \le \QQ(\beta)  .
\label{in:KM4p10}
\end{equation}

We now improve on the lower bounds on $\QQ(\beta)$ given in~(\ref{in:vBnds}) and in~(\ref{in:KM4p10}).
A reasonable choice for test functions $v$ to insert into ${\cal J}(v)$ is
\begin{equation}
 v
= c_0 + t_0 u_0 + t_\infty u_\infty .
\label{eq:vTest}
\end{equation}
In substituting the functions as in~(\ref{eq:vTest}) into ${\cal J}(v)$ it is useful to note 
that simple applications of the Divergence Theorem give
$$\QQ_0 =\int_\Omega |\nabla u_0|^2 \qquad {\rm and\ \ \ }
\Sigma_\infty =\int_\Omega |\nabla u_\infty|^2 .
$$
The Divergence Theorem also gives
$$\int_\Omega {\rm div} (u_0\nabla{u_\infty)}= 0 \qquad{\rm whence \ \ } 
\QQ_0 = \int_\Omega \nabla u_0 \cdot \nabla u_\infty .
$$
Now
\begin{eqnarray*}
{\cal A}(v) 
&=& 2c_0|\Omega| + t_0(2-t_0)\,\QQ_0 + t_\infty (2-t_\infty)\,\Sigma_\infty -2 t_0\, t_\infty\,  \int_\Omega \nabla u_0 \cdot \nabla u_\infty ,\\
&=& 2c_0|\Omega|+ t_0\,(2-t_0-2t_\infty)\, \QQ_0  +t_\infty (2-t_\infty)\,\Sigma_\infty  ,\\
{\cal B}(v)
&=& c_0^2 | \partial\Omega| - t_\infty^2\,\Sigma_1 \ \ .
\end{eqnarray*}
Next we notice that we have quadratics in $c_0$ and separately in $(t_0, \ t_\infty)$ to maximize.
Maximizing over $c_0$ gives
$$ c_{0,{\rm max}} = \beta\, \frac{|\Omega|}{ | \partial\Omega| } . $$

\cmdvmcode{VarlQuadr/genVarl2mpl.txt}{ }

Substituting for $c_0$ and considering the function so formed, quadratic in the $t$ variables,
we find that its hessian is negative semidefinite,
and that there is a unique maximum, at which $t_0+t_\infty=1$.
We calculate the maximizing $t$ and denoting them by $t_{\rm max}$, then use 
$ c_{0,{\rm max}}$ and $t_{0,\rm max}$, $t_{\infty,\rm max}$ to define, by equation~(\ref{eq:vTest}),
a function $v_{\rm max}$.
Then ${\cal J}(v_{\rm max})$ is a rational function of $\beta$.
Written in a form appropriate for $\beta$ not close to 0, 
\begin{eqnarray}
\QQ(\beta)\ge {\cal J}(v_{\rm max})
&=& R( \frac{|\Omega|^2}{|\partial\Omega|} ,\Sigma_\infty,\Sigma_1,\QQ_0;\beta), 
\label{in:lbbL}\\
R( \frac{|\Omega|^2}{|\partial\Omega|} ,\Sigma_\infty,\Sigma_1,\QQ_0;\beta)
&=& \left(\beta\, \frac{|\Omega|^2}{|\partial\Omega|} + \Sigma_\infty +\frac{1}{\beta} \Sigma_1\right) +
\frac{\Sigma_1^2}{\beta^2(\Sigma_\infty - \QQ_0) -\beta\Sigma_1} . 
\nonumber 
\end{eqnarray}
(On using $\Sigma_1\le{0}$ we see that this is an improvement on the left-hand side of~\cite{KM93}~(4.12).)
In a form appropriate for $\beta$ small, the rational function $R$ can be rewritten to give:
\begin{equation}
\QQ(\beta)\ge {\cal J}(v_{\rm max})
= \QQ_0 + \beta\, \frac{|\Omega|^2}{|\partial\Omega|} +
\frac{\beta(\Sigma_\infty-\QQ_0)^2}{\beta(\Sigma_\infty-\QQ_0)-\Sigma_1}
 . 
\label{in:lbbS}
\end{equation} 
This inequality improves on (4.4) of ~\cite{KM93}.

The function ${\cal J}(v_{\rm max})$
is increasing in $\beta$, a property which it shares with $\QQ(\beta)$.
 ${\cal J}(v_{\rm max})$ is concave in $\beta$:
 $$
\frac{\partial^2{\cal J}(v_{\rm max})}{\partial\beta^2}
= \frac{2\Sigma_1 (\Sigma_\infty - \QQ_0)^3}{(\beta (\Sigma_\infty -\QQ_0)-\Sigma_1)^3} .
$$

\section{Rectangular domains}\label{app:Rect}

\subsection{Rectangles: introduction}\label{subsec:Rect}

In this subsection we will consider the rectangle 
$$\Omega=(-a,a)\times(b,b)\qquad a=rh, \ b=\frac{h}{r} , $$
and usually treat the case $h=1$ when, of course we have the area of the rectangle as 4.
(Outside the context of flows with slip, the corresponding problem for box shapes 
in $R^n$ can be treated in exactly the same way. See~\cite{MK94}.)

We believe the value of our approximation, good lower bound~(\ref{in:lbbS},\ref{in:lbbS}),
is that it applies to all simply-connected domains.
We expect that it will be relatively uncommon for all five functionals involved in the inequality
to be explicitly available.
A more common situation would be that some would involve computation.
In this respect rectangles, 
besides being a geometry that is
relatively common in practice, is reasonably typical.

\subsubsection{The functionals $\QQ_0$, $\Sigma_\infty$ and $\Sigma_1$, etc.}

Before presenting the functionals which depend on solving pde problems we note
the geometric functionals which we use:
$$ |\Omega|=4a\, b = 4\, h^2,\ \  |\partial\Omega|=4 (a+b ) = 4\, h\, (r + \frac{1}{r} ) .
$$
 $B$ defined in~\cite{PoS51} is 
 $$B=4(a/b + b/a)
 = 4 (r^2 + \frac{1}{r^2} .)
 $$
 We use this expression for $B$ below in deriving a very simple lower bound for $\QQ(\beta)$.

The calculation of $\QQ_0$ is an undergraduate exercise in separation of variables:
see~\cite{LSU79} Problem 133, p74.
The quantity $\QQ_0$ is available only as an infinite series,
$$ \QQ_0	=	\frac{4 a^3b}{3}\left(
1-\frac{192}{\pi^5}\,\frac{a}{b}\sum_{n=1}^\infty \frac{\tanh(\pi(2n-1)b/(2a))}{(2n-1)^5})\right) .
$$
We don't know a closed form for its sum.
Also, indicating dependence on the sides as $\QQ_0(a,b)$, it is `physically'  evident that
$\QQ_0(a,b)=\QQ_0(b,a)$.
Also it is known (by different arguments in~\cite{PoS51}) that $\QQ_0(r,1/r)$ attains its
maximum for a square, $r=1$.
Neither of the last two sentences are immediate from the series.
However computation of $\QQ_0(a,b)$ isn't difficult, just one line of code.
The other quantities in the inequality~(\ref{in:lbbS},\ref{in:lbbS}) are extremely easy to calculate, and have very simple formulae, 
just rational functions of the rectangle's side lengths.
Once again it is clear that ${\cal J}(v_{\rm max})(a,b)$ is  unchanged if we swap $a$ and $b$
as this is true of all five functionals.

There are also series solutions for the case $\beta>0$ in $R^n$ as well as
the $R^2$ case relevant to flow in ducts: see~\cite{MK94,DM07,DG14,WWS14,KW16}.
The series are very substantially messier than the $\beta=0$ case, and require about 20 lines of code.
Once again $\QQ_\beta(a,b)-\QQ_\beta(b,a)$ should be zero
and we have found for rectangles with largish aspect ratios and truncated sums
there can be noticible differences.
If, however, an approximation is adequate, it is clear that the rational function of $\beta$
is far easier code.
Having said this, nevertheless, the value of our approximation is that it is simple, and general,
not merely for rectangular ducts.

Let's start the rectangle calculation.
The well-known series solution for $\QQ_0$ is easily coded:

\cmdvmcode{../../VarlLB/Maple/S0rectmpl.txt}{ }

Much easier are the calculations associated with Problem~(P($\infty$)).
It is easily verified that
$$ u_\infty
= \,{\frac {ab \left( {a}^{2}+6\,ab+{b}^{2} \right) }
{6 \left( a+b \right) ^{2}}}-{\frac {b{x}^{2}+a{y}^{2}}{2\,(a+b)}}
, $$
and hence
$$\Sigma_\infty
= {\frac {8 {a}^{3}{b}^{3}}{ 3\left( a+b \right) ^{2}}} .$$
From~\cite{KWK18} we note the easy result:\\
{\it Amongst rectangles with a given area, that which maximizes
$\Sigma_\infty$ is the square.}\\
These results (and further terms in the expansion for $\QQ_r(\beta)$ at large $\beta$) are presented in~\cite{KWK18}. 
Also
$$\Sigma_1
= - {\frac {4}{45}}\,{\frac {{a}^{2}{b}^{2}  \left( {a}^{4}+6\,b{a}^{3}-10\,{a}^{2}{b}^{2}+6\,a{b}^{3}+{b}^{4} \right) }
{ \left( a+b \right) ^{3}}}
$$
We remark that, both for ellipses and rectangles, when
$\beta$ is moderately large but the aspect ratio $a/b$ is
very large, the $\Sigma_1/\beta$ term can be relatively large
compared to $\Sigma_\infty$.
Asymptotics at large $\beta$ require care when the 
aspect ratio is large.

\cmdvmcode{../../VarlLB/Maple/rectLBmpl.txt}{ }

\subsubsection{Comparison of the lower bound with the series}

\begin{figure}[hb]
\centerline{\includegraphics[height=7cm,width=13cm]{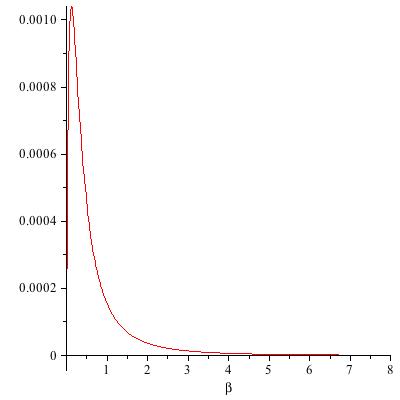}}
\caption{A plot at $a=5/4$, $b=4/5$ of the ratio 
$(S(\beta)-{\cal J}(v_{\rm max}))/{\cal J}(v_{\rm max})$ 
as a function of $\beta$.
This tends to 0 both as $\beta\rightarrow{0}$ and as  $\beta\rightarrow\infty$.
\label{fig:r5o4ratiojpg}}
\end{figure}

\medskip
We now compare the lower bound with the elaborate series solution: see Figure~\ref{fig:r5o4ratiojpg}.

\ifthenelse{\boolean{vmcode}}{
The code which produces this (and performs some checks) is}{}
\cmdvmcode{../../VarlLB/Maple/rectLBplotmpl.txt}{ }

\subsection{Using quadratic test functions}\label{subsec:VarlRectq}

As we did with the ellipse,
once again choose quadratic test functions.
This time we cannot expect useful results for $\beta$ very small.
Indeed, at $\beta=0$ we obtain just the trivial lower bound of 0.
However the results are good for $\beta$ large.
\begin{equation}
 Q_{r,LB}=
\,{\frac {4{a}^{2}{b}^{2}\beta\, 
 \left( 15\,a\beta+{b}^{2}+5\,ab\right)  
 \left( 5\,ab+{a}^{2}+15\,b\beta \right) }
 {3(30\,a{b}^{3}\beta+5\,{b}^{4}\beta+30\,{a}^{3}b\beta+5\,{a}^{4}\beta+
 75\,{a}^{2}b{\beta}^{2}+75\,a{b}^{2}{\beta}^{2}+
 2\,{a}^{2}{b}^{3}+2\,{b}^{2}{a}^{3})}}
\label{eq:QrLB}
\end{equation}
With $b=1/a$ when one plots both this lower bound and the series solution,
the curves are very close for all $\beta>1$.

\clearpage

\newpage

\begin{center}
{\Large{\textbf{\textsc{ Part IV:
$R$ approximated for ellipse\\
 \qquad \\
By G. Keady
}}}}
\end{center}

\section*{Abstract, Part IV}
Consider steady flows in a channel whose cross-section $\Omega$ is an ellipse, flows with the Navier slip boundary condition.
Denote the volume flow rate by $Q$.
We apply to elliptic cross-sections a recent simple approximation,  a rigorous lower bound $R$ on $Q$,
requiring, along with the channel's area and perimeter, the calculation of just the torsional rigidity and
two other domain functionals.
Here we make a further approximation, $R_a$, still, however finding a rigorous lower bound for $Q$.
$R_a$ is a good approximation to $Q$:
using it avoids the need for solving the partial differential equation repeatedly for
differing values of the slip parameter.

\section{Introduction}\label{sec:IntroRa}

Some of the papers concerning steady flows under a constant pressure gradient down microchannels are
~\cite{BYC,BTT09,TH11,Wa12,Wa12b}
and, for channels with elliptical 
cross-section~\cite{Duan07,SV12,Wa14,DT16}.

This paper is a sequel to~\cite{KW20} concerning approximating the steady slip flow down a channel with a given
cross-section, $\Omega$.
As noted in the abstract, just a few parameters are needed, and in cases tested to date, the approximation is remarkably accurate:
see~\cite{KW20} for rectangles and~\cite{Ke20i,Ke20} for regular polygons and isoceles triangles.
For the preceding shapes the `two other domain functionals', denoted below by $\Sigma_\infty$ and $\Sigma_1$, had simple exact formulae.
For the ellipse this is no longer the case
(as the exact formulae are infinite series with no obvious closed form), and this paper is a check that further approximation still leads to adequate results.

As in previous papers~\cite{KW20,KM93,Ke20i,Ke20} the slip-parameter is denoted by $\beta$.
The scaled velocity is denoted by $u(\beta)$; the boundary-value problem for $u(\beta)$
is given in 
Problem~(P($\beta$)) in~\cite{KW20} and in an earlier paper~\cite{KM93}.
The volume flow rate is defined by
\begin{equation}
 Q(\beta)= \int_\Omega u(\beta) .
\label{eq:QdefRa}
\end{equation}
For sufficiently smooth functions $v$, define, 
as in~\cite{KW20}, and~\cite{Ke20},~\cite{KM93}, 
with $\partial\Omega$ denoting the boundary of $\Omega$,
$\cal A$, $\cal B$  and $\cal J$ by
\begin{equation}
{\cal A}(v)= \int_\Omega(2 v -|\nabla v|^2),  \ \ {\cal B}(v) =\int_{\partial\Omega} v^2\ ,\ \ 
{\cal J}(v) ={\cal A}(v) - \frac{1}{\beta}{\cal B}(v) .
\label{eq:VarlRa}
\end{equation}
The maximiser of ${\cal J}(v)$ over all $v$ is $u(\beta)$ and
$$ \QQ(\beta)
= {\cal J}(u(\beta)) 
= \int_\Omega|\nabla u(\beta)|^2  - \frac{1}{\beta}{\cal B}(u(\beta)) .
$$
This variational characterisation of solutions $u(\beta)$  
led to the lower bound to and approximation of
$\QQ(\beta)$ by $R(\beta)$ defined, as in~\cite{KW20}, by
\begin{equation}
R(\beta)
= \QQ_0 + \beta \, \frac{ A ^2}{ P  } +
\frac{\beta (\Sigma_\infty-\QQ_0)^2}{\beta (\Sigma_\infty-\QQ_0)-\Sigma_1} .
\label{in:lbbSRa}
\end{equation}
The arguments of $R$ depending on the domain are:
$A$ the area of $\Omega$, $P$ its perimeter, 
$\QQ_0$ the torsional rigidity, and two other quantities $\Sigma_1$ and $\Sigma_\infty$.
The $\Sigma$ quantities were defined in~\cite{KM93}:
the definition is repeated in~\cite{KW20} and here in~\S\ref{subsec:Pinf}.
 (Another use of the approximation is that, from an experiment in which the flow rate $\QQ$ is measured,  
 the slip parameter can be estimated by simply solving,
 for $\beta$,  the polynomial equation $R(\beta)$ equals the measured flow $\QQ$.)
 
 For ellipses of given area $A$ and eccentricity $e$, the quantity $\QQ_0$ is elementary:
 see equation~(\ref{eq:Q0e}).
 The other quantities $P$ and approximations $\Sigma_{a,\infty}$ and $\Sigma_{a,1}$
 are given in terms of elliptic integrals at equation~(\ref{eq:perimERa}) and
 (on using (\ref{eq:i2}), (\ref{eq:C1})) equation~(\ref{eq:SigaInf}) and
(on using (\ref{eq:i4}), (\ref{eq:C2})) equation~(\ref{eq:Siga1}).
We emphasise that the elliptical microchannels have been studied extensively.
There is no claim that the numerical results from this paper are better than earlier ones.
What is new is the uniform treatment of potentially all cross-sections by the approximation given in $R$,
adding elliptical cross-sections to the list where this has been done.

\subsection{The torsion pde, Problem (P(0))}\label{subsec:P0}

The elastic torsion problem is to find a function $u_0$ satisfying
$$ - \frac{\partial^2 u_0}{\partial x^2}- \frac{\partial^2 u_0}{\partial y^2}
 = 1\, {\rm in}\,\Omega ,\,\qquad
u_0
= 0 \, {\rm on}\,\partial\Omega \, .
\eqno({\rm P}(0))
$$
Channel flows with no slip have $\beta=0$.
The {\it torsional rigidity} of $\Omega$ (volume flow in the fluids context) is,
\begin{equation}
\QQ_0
:= \int_\Omega u_0 .
\label{eq:Q0def}
\end{equation}

With $\cal A$ defined as in equation~(\ref{eq:VarlRa})
the solution $u_0$ maximizes $\cal A$ over functions vanishing on the boundary $\partial\Omega$.
Also ${\cal A}(u_0)=\QQ_0$.

Define also
\begin{equation}
\QQ_1
:=  \int_{\partial\Omega}  \left(\frac{\partial u_0}{\partial n}\right)^2 .
\label{eq:Q1def}
\end{equation}
For small $\beta$, as noted in equation~(4.5) of~\cite{KM93},
$$ \QQ(\beta) \sim \QQ_0 + \beta\, \QQ_1 \qquad {\rm as \ \ }\beta\rightarrow{0} . $$
This suggested the following inequality (proved in~\cite{Ke20}) which is used in~\S\ref{subsec:Further}
as a check on the quantities occurring in 
our approximation:
\begin{equation}
\frac{A^2}{P}-\frac{(\Sigma_\infty-\QQ_0)^2}{\Sigma_1}\le \QQ_1  .
\label{eq:QQ1R1a}
\end{equation}


\subsection{Problem (P($\infty$))}\label{subsec:Pinf}

Problem (P($\infty$) was defined in~\cite{KW20}
and, before that in
\cite{KM93}.
Our notation here is as before.
The function $u_\infty$ solves
$$ - \frac{\partial^2 u_\infty}{\partial x^2}- \frac{\partial^2 u_\infty}{\partial y^2}
 = 1\, {\rm in}\,\Omega ,\,
\frac{\partial u_\infty}{\partial n}
= -\frac{|\Omega|}{|\partial\Omega|} \, {\rm on}\,\partial\Omega \,
{\rm and}\,  \int_{\partial\Omega} u_\infty = 0 .
\eqno({\rm P}(\infty))
$$
Define, as in~\cite{KW20},~\cite{KM93} 
\begin{equation} 
\Sigma_\infty
= \int_\Omega u_\infty , \qquad  {\rm and\ \ } \Sigma_1= - \int_{\partial\Omega} u_\infty^2 ,
\label{eq:SigDefRa}
\end{equation}
with $u_\infty$ satisfying Problem (P($\infty$)).

Once again there is a variational characterisation of the solutions.
This time $u_\infty$ is the maximizer of ${\cal A}(v)$ as one varies over functions $v$
for which the integral of $v$ over the boundary $\partial\Omega$ is zero.
The inequality $\Sigma_\infty \ge \QQ_0$
is an immediate consequence of the variational approach.
(See also~\cite{KM93} equation~(4.9).) 
Also ${\cal A}(u_\infty)=\Sigma_\infty$.

As noted in~\cite{KM93}, for large $\beta$,
$$ \QQ(\beta ) \sim \frac{\beta \,  A ^2}{ P  } + \Sigma_\infty +\frac{\Sigma_1}{\beta }
\qquad {\rm for\ \ } \beta \rightarrow\infty . $$

\section{The ellipse}\label{ec:Ell}
Consider ellipses with area $A=\pi$, boundary specified, with $a\ge 1$, by
\begin{equation}
\frac{x^2}{a^2}+a^2 y^2  = 1 , 
\qquad A=\pi .
\label{eq:ellabRa}
\end{equation}
The eccentricity $e$ of the ellipse satisfies $e^2=1-1/a^4$.

The elliptic integral EllipticE is denoted by $\EllipticE$:
EllipticK is denoted by $\EllipticK$.
The perimeter of the ellipse, $P$, 
is
\begin{equation}
P = 4 a {\EllipticE}(e) \quad{\rm for\ \ } a\ge{1} .
\label{eq:perimERa}
\end{equation}

The approximations to the $\Sigma$ quantities (and $Q_1$) are readily
expressed in terms of boundary moments, these being defined by
 \begin{equation}
i_{2k} 
=  \int_{\partial\Omega} (x^2+y^2)^k . 
\label{eq:i2k}
\end{equation}
Evaluating these for the ellipse gives
\begin{eqnarray}
i_2
&=& \frac{4}{3 a}\left( (1+a^4)\, {\EllipticE}(e)+{\EllipticK}(e) \right) ,
\label{eq:i2} \\
i_4
&=& \frac{4}{5 a^3}\left( (1-a^4+a^8)\, {\EllipticE}(e) +2(1+a^4)\, {\EllipticK}(e) \right)
\label{eq:i4}
\end{eqnarray}

A well-known special case is the circular disk of radius 1 where
$u_{\odot}(\beta) =(1-x^2-y^2)/4+\beta/2$ and
$\QQ_\odot(\beta)=\pi(1+4\beta)/8$.
It is also known that for any ellipse with area $\pi$, $\QQ(\beta)\le\QQ_\odot(\beta)$.

\subsection{$\QQ_0$ and $\QQ_1$}\label{subsec:Qell}

We have
\begin{eqnarray}
u_0
&=& \frac{1-\frac{x^2}{a^2}- a^2 y^2}{2\,(\frac{1}{a^2}+a^2)}  ,
\label{eq:u0e}
\\  
\frac{\partial u_0}{\partial n}
&=& - \, \frac{\sqrt{(\frac{1}{a^2}+a^2)-(x^2+y^2)}}{\frac{1}{a^2}+a^2}
\quad {\rm on\ \ } \partial\Omega .
\label{eq:u0nBdry} 
\end{eqnarray}
and hence
\begin{eqnarray}
\QQ_0
&=& \frac{\pi a^4}{4\,(1+a^4)} ,
\label{eq:Q0e}
\\  
\QQ_1
&=& \frac{a^2 P}{1+a^4} -\frac{a^4 i_2}{(1+a^4)^2} ,\nonumber 
\\
&=& \frac{4}{3}\, \frac{a^3}{(1+a^4)^2}
\left( 2(1+a^4){\EllipticE}(e)- {\EllipticK}(e)\right)
\ {\rm for\ } a\ge{1} .
\label{eq:QQ1e} 
\end{eqnarray}

\subsection{Approximations $\Sigma_{a,\infty}$ and $\Sigma_{a,1}$}\label{subsec:SigApprox}

The quantities $\Sigma_\infty$ and $\Sigma_1$ can be approximated
by terms denoted below by  $\Sigma_{a,\infty}$ and $\Sigma_{a,1}$
which, as mentioned above, are calculated in terms of boundary moments.


We seek an approximation to $u_\infty$ by using test functions of the form
$$ v = u_0 + t\left( x^2-y^2 - \frac{C_1}{P}\right) , $$
where the constant $C_1$ is chosen so the boundary integral of $v$ is zero, i.e.
\begin{equation}
C_1= \int_{\partial\Omega}(x^2-y^2)
= \frac{ -2 a^2 P+ (1+a^4) i_2}{(-1+a^4)} .
\label{eq:C1}
\end{equation}
The quantity ${\cal A}(v)$ is quadratic in $t$
$$ {\cal A}(v)=  - \frac{\pi\left(4(1+a^4)^2 t^2 - 2(1+a^4)(-1+a^4 -4a^2 \frac{C_1}{P}) t -a^4\right))}{4 a^2 (1+a^4)} . $$
The maximum of this quadratic occurs at $t=t_*$ with
$$ t_* = \frac{ -1+a^4-4 a^2 \frac{C_1}{P}}{4(1+a^4)} .$$
Then, at $t_*$,
\begin{equation}
{\cal A}(v_*)= \frac{\pi\left(a^4 +\frac{1}{4}(-1+a^4 -4 a^2 \frac{C_1}{P})^2\right)}{4 a^2 (1+a^4) } . 
\label{eq:SigaInf}
\end{equation}
Define the approximation to $\Sigma_\infty$ by
$$\Sigma_{a,\infty} = \int_\Omega v_* , $$
and then, after performing the integrals, note the remarkable fact that 
$$\Sigma_{a,\infty} ={\cal A}(v_*) = \int_\Omega |\nabla v_*|^2 . $$
(The second equality follows immediately from the first.
The first equality is checked by performing the integrals.
We note it as `remarkable' because, although we can prove using applications of the
Divergence Theorem that  $\Sigma_\infty={\cal A}(u_\infty)$, it isn't obvious there would be
a parallel result for $v_*$.)
As ${\cal A}(u_\infty)>{\cal A}(v_*)$, our approximation underestimates $\Sigma_\infty$, i.e. 
$\Sigma_\infty>\Sigma_{a,\infty}$.

Higher boundary moments enter in the approximation to $\Sigma_1$.
We have
$$\int_{\partial\Omega} v^2
= t^2\left( C_2 - \frac{C_1^2}{P}   \right) .
$$
Substituting $t=t_*$ gives, with
\begin{equation}
C_2 
= \frac{ 4 a^4 P -4 a^2 (1+a^4) i_2+ (1+a^4)^2 i_4}{(-1+a^4)^2} , 
\label{eq:C2}
\end{equation}
\begin{equation}
\Sigma_{a,1}
= - \int_{\partial\Omega} v_*^2
= - \, \frac{\left(-4 a^2 \frac{C_1}{P}+a^4-1\right)^2  \left( C_2 
- \frac{C_1^2}{P} \right)}
  {16 \left(a^4+1\right)^2} .
\label{eq:Siga1}
\end{equation}

We have checked that the approximations $\Sigma_{a,\infty}(a)$, $\Sigma_{a,1}(a)$ when used in the obvious
places in inequality~(\ref{eq:QQ1R1a}) do satisfy the inequality.

\goodbreak

\subsection{Comparison of $R_a(\beta)$ with $\QQ(\beta)$}\label{subsec: Comparison}

\begin{figure}[h]
\centerline{\includegraphics[height=5cm,width=7.5cm]{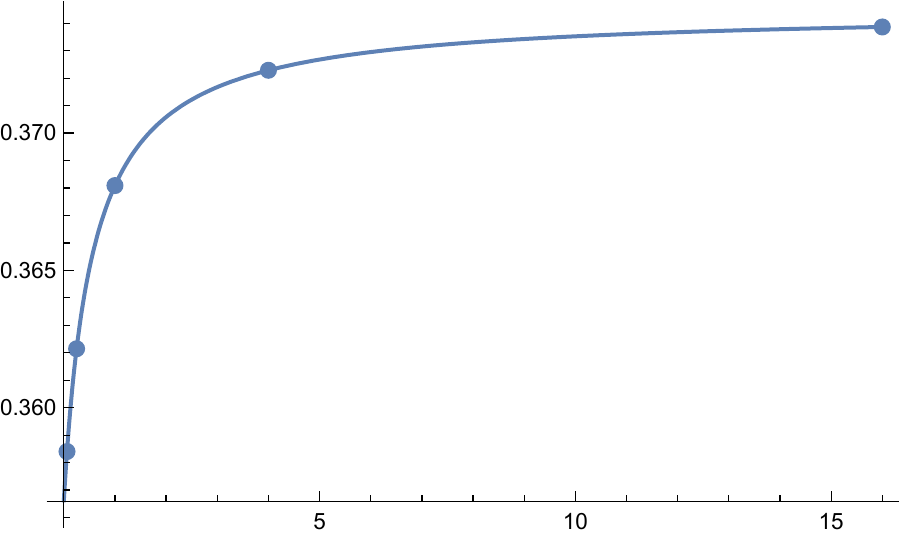}}
\vspace{-0.2cm}
\caption{For an ellipse with $a=5/4$, the quantity $R_a(\beta)-\beta A^2/P$
plotted against $\beta$.
The plot goes between $Q_0$ and $\Sigma_{a,\infty}$.
The solid dots are values when one replaces $R_a$ with $Q$ calculated numerically. }
\label{fig:ellRPlot}
\end{figure}

Write $R_a(\beta)$ as the expression obtained when one uses $\Sigma_{a,\infty}$ and  $\Sigma_{a,1}$ 
in the obvious places in the parameter list of $R$.
The numeric values for $\QQ(\beta)$ used in Figure~\ref{fig:ellRPlot} have been obtained by several methods.
The PDE Toolbox of Matlab is one source (and could be, and has been, used with other cross-sections).
The variational principle used for the theory in~\cite{Ke20} has been used for numerical values
for $\QQ(\beta)$ for ellipses in~\cite{Wa14}.
Fourier series for $u(\beta)$ is used in~\cite{Duan07} with, however, a further approximation:
that approximation is avoided in numerics performed by BW reported 
in~Part I.

\subsection{Further remarks}\label{subsec:Further}

\subsubsection{Approximation by quadratics for $\beta\ge{0}$}

The numeric values of the approximation  
coincide with those of the quadratic approximation of Part I.

There we used test functions of the form
$$ v_{\rm quad} = c + c_{xx} x^2 + c_{yy}y^2 , $$
in maximizing ${\cal J}(v)$, $\cal J$ defined at equations~(\ref{eq:VarlRa}).
This was completely ad hoc: it was suggested by the contours `looking elliptical'
e.g. in~\cite{SV12}.
Write $v_{\rm quad,max}$ for the maximiser.
A calculation shows
$${\cal J}{v_{\rm quad,max}} = R_a(\beta) . $$
Thus our approximate $R_a(\beta)$, as with the $R(\beta)$ involving
the true $\Sigma$ values,
is a lower bound for $\QQ(\beta)$.

\subsubsection{Higher order approximation for $u_\infty$}

We suspect that finding higher approximations to $u_\infty$ and the $\Sigma$ quantities
is likely to be unimportant for the ellipse,
but they could be found.
One method begins by defining an inner product
$$\langle \phi, \psi\rangle_{1\Omega}
= \int_\Omega \nabla \phi\cdot\nabla \psi , $$
and noting that $\langle u_0,h\rangle_{1\Omega}=0$ for any harmonic function.
Suppose one has a set of $N$ orthogonal harmonic functions $h_j$ each of which has
boundary integral zero.
Use test functions $v$ of the form
$$ v = u_0 + \sum t_j h_j , $$
in the variational characterisation.
It is easy to maximize ${\cal A}(v)$ over the $t$, and denoting the maximizers
with a subscript star
$$ t_{j*} = \frac{\int_\Omega h_j}{\langle h_j,h_j\rangle_{1\Omega}} . $$

Return now to the ellipse.
Take as the $h_j$
\begin{equation}
h_j = {\rm Re}\left(T_{2j}\left(\frac{x+iy}{a e}\right)\right) + C_{1,j} , 
\label{eq:hjDef}
\end{equation}
where $T_{2j}$ denote Chebyshev polynomials and the
$C_{1,j}$ are chosen so that the boundary integrals are zero.
The $h_j$ are orthogonal.
Also
$$ \int_\Omega h_1 = -\frac{\pi}{2} + A\, C_{1,1},\ \ 
 \int_\Omega h_j =  A\, C_{1,j}\ {\rm for\ \ } j>1,
$$
so that the integrals defining the $t_{j*}$ can be evaluated.

The series one gets approximating $u_\infty$ using this variational method is
the same as that found using separation of variables in elliptic coordinates
and the Fourier series which follows from that.

\section{Conclusion and future prospects}

The approximation $R(\beta)$ presented in~\cite{KW20}  for the flow $\QQ(\beta )$ has been determined to be good for elliptical cross-sections.

For any domain it is easy to use finite elements or other numerical methods to find
numerical approximations to the five parameters occuring in $R$.
However, having a few cross-sections where some or all of the parameters are available analytically is useful in checking numerics.

The ellipse now joins the list initiated in~\cite{KW20} for rectangles, and continued in~\cite{Ke20i,Ke20} for shapes like
regular polygons, triangles, rhombi, a list where all the parameters of $R$ are known, or estimated.
For many years there have been tables for $\QQ(0)$, the torsional rigidity, one of the parameters in $R$.
Looking to the future, a similar table -- online lookup -- of all five of the parameters would enable experimenters to calculate quickly an approximation to the flow for any value of $\beta $.

\clearpage

\newpage
\appendix

\section*{Outline of appendices}

The notes in the appendices were prepared for a student of BW 
further developing knowledge of $u_\infty$ for the ellipse.
This also influenced how GK wrote papers, in order to leave a clear plan 
for where the student's paper might be submitted.

\begin{enumerate}
\item[A] Quadratic test functions, general $\beta$

\item[B] (P($\infty$)), more terms: agreement with Fourier Series\\
B.1-B.3 General $\Omega$, including basics of the variational principle\\
B.4-B.5 More on polynomial test functions

\item[C] Conformal maps, ellipses

\item[D] More recurrences, moments

\ifthenelse{\boolean{vmcode}}{
\item[E] Code
}{ }
\end{enumerate}

There are very many open questions concerning $u_\infty$ for the ellipse.
\begin{itemize}
\item
We believe series expansions for small eccentricity should be found.

\item
For eccentricities tending to 1, matched asymptotic expansions might be
useful.

\item
It might be possible to make use of the known conformal map of the ellipse to 
unit disk.

\item Recurrences for $t_{j*}$, or $V_j$ of Part II might be found.\\
The geometric quantities $I_{2k}$ and $i_{2k}$ may have other uses and the 
recurrences explored in a later Appendix may have interest.
(GK has more results not presented in this report.)

\end{itemize}

\section{Agreement with quadratic approximation}\label{app:GenBetaQuad}

In Part I 
we used test functions of the form
$$ v_{\rm quad} = c + c_{xx} x^2 + c_{yy}y^2 , $$
in maximizing ${\cal J}(v)$, $\cal J$ defined at equations~(\ref{eq:VarlRa}).
The values of the $c$ such that ${\cal J}(v_{\rm quad})$ is maximized were found.
This was completely ad hoc: it was suggested by the contours `looking elliptical'
e.g. in~\cite{SV12}.
Write $v_{\rm quad,max}$ for the maximiser.
We did this for any geometry, but it probably only works well, at least for small $\beta$,
for the ellipse.
For the ellipse the relations between the moments
(recorded at ~\ref{subsec:momentsa}, ensure that
the $c$ come out so that the Laplacian of $v_{\rm quad,max}$ is minus one.
Thus, relating the two approaches, we have, for some $c$ and $\tau$ depending on $\beta$:
$$ v_{\rm quad,max} = c +\tau u_0 + (1-\tau) u_{a,\infty} .$$

Next, if we maximise ${\cal J}$ over functions as on the right of the preceding equation,
maximising over $(c,\tau)$ one finds a maximiser $v_*$ and
$${\cal J}{v_*} = R_a(\beta) . $$
Thus our approximate $R_a(\beta)$, like the $R(\beta)$ involving
the true $\Sigma$ values,
is a lower bound for $\QQ(\beta)$.


\section{More terms, agreement with Fourier Series}

\subsection{Introduction, and (P($\infty$)) for general $\Omega$}\label{subsec:GenOma}

One of the many questions yet to be treated well is
 how close are GK's approximations $\Sigma_{a,\infty}$ and $\Sigma_{a,1}$
to the true values of $\Sigma_\infty$ (which we know is bigger than  $\Sigma_{a,\infty}$)
and of $\Sigma_1$?

It isn't clear how best to represent higher approximations.
The methods include
\begin{itemize}
\item Fourier series 
which is equivalent to the variational characterization of $u_\infty$
using Chebyshev polynomials, and
\item the variational characterization of $u_\infty$ using harmonic polynomials, thus the next term after that used in $R_a$ might be approached as
\begin{eqnarray*}
v = u_0 &+& {\tilde t}_1\left( x^2-y^2 -\frac{\int_{\partial\Omega}(x^2-y^2)}{P}\right) +\\
 & &{\tilde t}_2\left( x^4-6 x^2 y^2 + y^4 -\frac{\int_{\partial\Omega}(x^4-6 x^2 y^2 + y^4 )}{P}\right) .
\end{eqnarray*}
An outline of how this might be done is given, in the context of general domains,
below in~\S\ref{subsec:BetterVa}. 
 \end{itemize}
 Of course, when there are several  methods, one must cross-check.
 
\subsection{The variational principle for (P($\infty$))}\label{subsec:GVarla}

The basics of the variational principle (somewhat formally,
as in the usual undergraduate treatment of Euler-Lagrange d.e.)
is as follows.

The standard approach to variational methods for applied mathematicians
(as opposed to more pure-maths approaches involving Frechet derivatives)
starts with a linearization, $\epsilon$ small
$$ {\cal A}(u_\infty +\epsilon \phi) - {\cal A}(u_\infty)
\sim 2\epsilon {\cal I}(u_\infty,\phi)\quad
{\rm where\ \ }  
{\cal I}(u_\infty,\phi)=\int_\Omega(\phi - \nabla u_\infty\cdot\nabla\phi) , $$
with
$$\int_{\partial\Omega}\phi =0 \qquad{\rm and\ \ } \int_{\partial\Omega} u_\infty = 0 . $$
(The boundary integrals have to be zero in order that the functions are allowed in the
variational competition.)

For $u_\infty$ to be the maximum, we must have the $O(\epsilon)$ term, $\cal I$, is 0.
There are two parts to showing that $u_\infty$ satisfies the pde Problem~(P($\infty$)),
namely (i) the pde (Laplacian is -1) and (ii) the boundary condition that the normal derivative is constant.

\noindent{(i)}
 $$ {\rm div}(\phi\nabla{u_\infty})
 = \phi\Delta u_\infty + \nabla u_\infty\cdot\nabla\phi , $$
 so
 \begin{equation}
0 =  {\cal I}(u_\infty,\phi)\ = \int_\Omega(\phi(1+\Delta{u_\infty}) - {\rm div}(\phi\nabla{u_\infty})) .
\label{eq:Idiv}
\end{equation}
 Now consider smooth functions $\phi$ which vanish on the boundary.
 We have $  \int_\Omega(\phi(1+\Delta{u_\infty})=0$ for all such $\phi$ so we must have
 $-\Delta u = 1$.
 (For suppose not and $(1+\Delta{u_\infty})>0$ in some domain.
 Then choose a nonzero function $\phi\ge{0}$ which has its support in that domain.
 The contradiction is clear.)

\noindent{(ii)} Now consider maximisers satisfying $-\Delta u = 1$, so
$$ \int_{\partial\Omega} \frac{\partial u_\infty}{\partial n} = - A,
\qquad{\rm and\ as\ before \ \ }  \int_{\partial\Omega}\phi = 0. $$
From equation~(\ref{eq:Idiv})
$$  0 = \int_\Omega {\rm div}(\phi\nabla{u_\infty} )
= \int_{\partial\Omega} \phi \frac{\partial u_\infty}{\partial n}  
=  \int_{\partial\Omega} \phi \left(\frac{\partial u_\infty}{\partial n}+\frac{A}{P}\right) ,$$
as the integral of $\phi$ around the boundary is 0.
Now
$$  \int_{\partial\Omega} \left(\frac{\partial u_\infty}{\partial n}+\frac{A}{P}\right) = 0 ,$$
so if we take
$$\phi = \left(-\frac{\nabla u_\infty\cdot \nabla u_0}{|\nabla u_0|}+\frac{A}{P}\right) $$
which is a function whose boundary integral is 0, we must then have
$$  \int_{\partial\Omega} \left(\frac{\partial u_\infty}{\partial n}+\frac{A}{P}\right)^2 = 0 ,$$
and hence the integrand 0, which establishes the boundary condition 
of Problem~(P($\infty$)).

\subsection{Better approximation, general $\Omega$}\label{subsec:BetterVa}

Define an inner product
$$\langle f_1, f_2\rangle_{1,\Omega}
= \int_\Omega \nabla f_1 \cdot \nabla f_2 . $$
For any harmonic function $h$, using the divergence theorem on
${\rm div}(u_0\nabla{h}$, one finds
$$ \langle u_0,h\rangle_{1,\Omega}
= 0 . $$

Given any sequence of harmonic functions ${\tilde h}_k$, 
for each of which the boundary integral is zero,one can apply Gram-Schmidt to
find an orthogonal sequence of harmonic functions $h_k$:
$$ h_k = {\tilde h}_k - \sum_{j=1}^{k-1}
\frac{\langle {\tilde h}_k, h_j\rangle_{1,\Omega}}{\langle h_j, h_j\rangle_{1,\Omega}}\, h_j .$$
If necessary change the sign of $h_k$ so that $\int_\Omega h_k >0$.
The boundary integrals of the $h_k$ are also zero.

The function
$$ v = u_0 +t_1 h_1 + t_2 h_2 + \ldots , $$
is suitable as a test function in the problem of maximizing
${\cal A}(v)$.
Now
\begin{eqnarray*}
{\cal A}(v) 
&=& 2\int_\Omega v - \langle v,v\rangle_{1,\Omega} ,\\
&=& 2\int_\Omega v - \langle u_0,u_0\rangle_{1,\Omega}
-\sum_{j=1}^N t_j^2 \langle h_j,h_j\rangle_{1,\Omega} ,
\end{eqnarray*}
on using the orthogonality properties we have noted above.
Thus
$${\cal A}(v)
= \QQ_0 -\sum_{j=1}^N \left( t_j^2  \langle h_j,h_j\rangle_{1,\Omega}
-2\,t_j \int_\Omega h_j \right) .
$$
There are no `cross terms' $t_j\,t_k$ with $j\ne{k}$.
It is easy to maximize over the $t$, and denoting the maximizers 
with a subscript star
$$ t_{j*} = \frac{\int_\Omega h_j}{\langle h_j,h_j\rangle_{1,\Omega}} . $$
We expect the $N$-th approximation to be
$$\Sigma_{\infty,N}
= {\cal A}(v_{*N}) 
=\QQ_0 + \sum_{j=1}^N t_{j*} \int_\Omega h_j . $$

\subsection{Harmonic polynomials for higher approximations via the variational principle?}

We now suggest how the general set-up described in~\S\ref{subsec:BetterVa}
might be applied to the case of the ellipse.

Without the orthogonality described in~\S\ref{subsec:BetterVa}
the approximations would be in the form
$$ u_N = u_0 + \sum_{j=1}^N c_j(N) {\hat h}_j(x,y) , $$
with the ${\hat h}_j$ appropriate harmonic functions with boundary integrals being $0$.
The $c_j(N)$ are then chosen so that $\cal A$ is maximized.
Amongst the possibilities for the ${\hat h}_j$ are the following
\begin{itemize}
\item ${\hat h}_j = {\rm Re}(x+i y)^{2j} - c_j$ with the $c_j$ chosen so the boundary integral of $h_j$ is $0$.
\item ${\hat h}_j = c \cosh(2j\eta) \cos(2j\psi) - {\hat c}_j$  as before with $c^2=a^2-1/a^2$, and again
with ${\hat c}_j$ such that the boundary integral of ${\hat h}_j$ is zero.\\
This is equivalent to the Chebyshev polynomials treated in the next subsection.
\end{itemize}

We discovered in~\S\ref{subsec:BetterVa} that it will be possible to have a sequence $h_j$ which has the $c_j(N)$ independent of $N$.
The Chebyshev polynomials of the next subsection provide such a sequence.


\subsection{Chebyshev polynomials with complex argument}

\subsubsection{The polynomials}

We begin with a few facts. As before $c=\sqrt{a^2 - 1/a^2}$. Then
$$ \int_\Omega T_{2}( \frac{x+iy}{c} )= -\frac{\pi}{2},\qquad
\int_\Omega T_{2j}( \frac{x+iy}{c} ) = 0\ {\rm for\ \ } j>1. $$
Also
$$ \frac{d T_n(\frac{z}{c})}{d z} = \frac{n}{c}\, U_{n-1}(\frac{z}{c}). $$
Refer now to equation~(\ref{eq:5}) of Part I concerning elliptic coordinates.
$$ T_n(\cosh(\eta)\cos(\psi) + i\, \sinh(\eta)\sin(\psi))
= T_n(\cosh(\eta+ i\, \psi)) 
= \cosh(n(\eta+ i\, \psi)) ,
$$
so
$$ {\rm Re}(T_n(\cosh(\eta)\cos(\psi) + i\, \sinh(\eta)\sin(\psi)))
= \cosh(n\eta)\cos(n\psi) .
$$

\subsubsection{Integrals over $\Omega$}

Let $z=x+iy$. Since
$$ z^m {\bar z}^n - {\bar z}^m z^n = |z|^n\, (z^{m-n}- {\bar z}^{m-n}) , $$
it is imaginary: thus
we have ${\Re}( (x+iy)^m (x-iy)^n)$ is even in $y$.
Similarly
$$ z^m {\bar z}^n + {\bar z}^m z^n = |z|^n\, (z^{m-n}+ {\bar z}^{m-n}) , $$
it is real: thus
we have ${\Im}( (x+iy)^m (x-iy)^n)$ is odd in $y$.
Thus
$$ \int_\Omega U_m(\frac{z}{c}) {\bar z}\ \ \
{\rm is\ real,} $$
and hence
$$ \int_\Omega U_m(\frac{z}{c}) U_n(\frac{\bar z}{c})\ \ \
{\rm is\ real.} $$
The value of the integral is, with 
$r=a +1/a$, is given by 
\begin{equation}
\int_\Omega U_m(\frac{x+i y}{c}) \, U_n(\frac{x-i y}{c}) 
= \frac{\pi c^2}{4(1+n)} \left(
(\frac{r}{c})^{2n+2} -(\frac{c}{r})^{2n+2} \right)\,
\delta_{n,m} . 
\label{eq:UmUnA}
\end{equation}
(See~\cite{ANPV} equation~(4.25).)
When $m=n$ we can write the right-hand side of (\ref{eq:UmUnA}) as follows.
First
$$
\frac{r^2}{c^2} = \frac{a+1/a}{a-1/a}
= \exp(2\eta_0) ,$$
on using $a=c\cosh(\eta_0)$, $1/a=c\sinh(\eta_0)$.
Hence
\begin{eqnarray}
\int_\Omega \left| U_n(\frac{x+iy}{c})\right|^2
&=& \frac{\pi c^2}{4(1+n)} \left(
(\frac{r}{c})^{2n+2} -(\frac{c}{r})^{2n+2} \right)\,
\nonumber \\
&=& \frac{\pi c^2}{2(1+n)}
\sinh(2(n+1)\,\eta_0) ,
\nonumber \\
&=& \frac{\pi }{\sinh(2\eta_0)(1+n)}
\sinh(2(n+1)\,\eta_0)
\label{eq:UmUnA1}
\end{eqnarray}

We now establish the orthogonality of the $h_j$ of equation~(\ref{eq:hjDef}).
Now
{\small
\begin{eqnarray}
<h_j,h_k>_{1,\Omega}
&=& \int_\Omega \nabla h_j \cdot \nabla h_k , 
\nonumber\\
&=& \int_\Omega \nabla\left({\rm Re}({T_{2j}(\frac{x+iy}{c}))}\right)
\cdot \nabla\left({{\rm Re}(T_{2k}(\frac{x+iy}{c}))}\right) .
\label{eq:hjhkReT}
\end{eqnarray}
} 
Next, for holomorphic functions $f_j$ and $f_k$,
\begin{eqnarray*}
{\rm Re}\int_\Omega f_j(z)\ {\overline{f_k(z)}}
&=& {\rm Re}\int_\Omega (f_{jr}+i f_{ji})\, (f_{kr} - i f_{ki}) ,\\
&=& \int_\Omega (f_{jr}f_{kr}+ f_{ji}f_{ki}) ,
\end{eqnarray*}
where the r and i in the subscripts indicate real and imaginary parts
respectively.
Applying this to derivatives and using the Cauchy-Riemann equations
gives
$${\rm Re}\int_\Omega \frac{d f_j}{ d z}\ {\overline{\frac{d f_k}{d z}}}
= \int_\Omega \nabla f_{jr}\cdot \nabla f_{kr} . $$
Using this with the Chebyshev $T$ functions in
equation~(\ref{eq:hjhkReT}) gives
\begin{eqnarray*}
<h_j,h_k>_{1,\Omega}
&=& \int_\Omega \nabla\left({\rm Re}({T_{2j}(\frac{x+iy}{c}))}\right)
\cdot \nabla\left({{\rm Re}(T_{2k}(\frac{x+iy}{c}))}\right) ,\\
&=& {\rm Re}\int_\Omega
\frac{d T_{2j}(\frac{x+iy}{c})}{ d z}\ 
{\overline{\left(\frac{d T_{2k}(\frac{x+iy}{c})}{d z}\right)}} ,\\
&=& \frac{4 j k}{c^2}
{\rm Re}\int_\Omega
\left({U_{2j-1}(\frac{x+iy}{c})}\right)
\, \left({U_{2k-1}(\frac{x-iy}{c})}\right) .
\end{eqnarray*}
This final expression is zero for $j\ne{k}$ and when $j=k$
its value is given as in equation~(\ref{eq:UmUnA}).
\begin{equation}
<h_j,h_j>_{1,\Omega}
= \pi j \, \sinh(4 j\,\eta_0)
\label{eq:ghjSqr}
\end{equation}

\subsection{Integrals over $\partial\Omega$}

We have
\begin{eqnarray}
L\, C_{1,j}
&=& -\int_{\partial\Omega} {\rm Re}(T_{2j}(\frac{x+iy}{c}))\, ds ,
\nonumber \\
&=& -2\int_0^\pi
{\rm Re}(T_{2j}(\cosh(\eta_0)\cos(\psi) + i\sinh(\eta_0)\sin(\psi)))\,
\sqrt{ a^2 \sin(\psi)^2+ \frac{\cos(\psi)^2}{a^2}}\, d\psi ,
\nonumber\\
&=& -2a\cosh(2j\eta_0) \int_0^\pi \cos(2j\psi)
\sqrt{1- e^2\cos(\psi)^2}\, d\psi ,
\nonumber \\
&=& -\pi a\cosh(2j\eta_0) {\hat g}_j .
\label{eq:C1jVarl}
\end{eqnarray}

From equations~(\ref{eq:C1jVarl}) and~(\ref{eq:ghjSqr}) one finds
$t_{j*}$ and finds that for $j>1$ it agrees with the Fourier series
calculation at~(\ref{eq:Vninf}).
 

\section{Conformal maps, ellipses}

We remark that the normal to $\partial\Omega$ is
\begin{itemize}

\item in the direction of $\nabla u_0$,

\item in the direction of $\log(|f(z)|)$ where $f$ is the
conformal map from $\Omega$ to the unit disk in which the origin is mapped to
the origin.\\
The function  $\log(|f(z)/z|)$ is harmonic in $\Omega$.

\end{itemize}

We remark that there is
a `Poisson Integral Formula' for solving the Dirichlet problem for
Laplace's equation in an ellipse.
This follows from the conformal map between ellipse and disk:
see~\cite{Kob}, p177 and~\cite{Sz50}.
See also~\cite{Ro64,Mi90,LK14,SD16}.


\section{More recurrences}\label{sec:Nicea}

\subsection{Geometric quantities: moments}\label{subsec:momentsa}

The quantities $\Sigma_\infty$ and $\Sigma_1$ can be expressed in terms of boundary moments
 -- moments about the incentre -- defined by
 \begin{eqnarray}
 i(m,n)
 &=& \int_{\partial\Omega} (x^{2m}\, y^{2n})
 \label{eq:imn}
\\  
i_{2k} 
&=& i(k,k) =  \int_{\partial\Omega} (x^2+y^2)^k . 
\nonumber 
\end{eqnarray}
We remark that the Cauchy-Schwarz inequality for the integrals implies that
$i_4\ge{i_2^2/P}$.

We now record results needed in the approximation leading to $v_{max}$
in~\S\ref{subsec:SigApprox}.
\begin{itemize}
\item For second moments we have the following:
$$ i(1,0)+i(0,1) = i_2, \qquad \frac{i(1,0)}{a^2} + a^2\,{i(0,1)} = P .$$
From this
$$ i(1,0) = \frac{ (a^2 (a^2 i2 - P))}{(-1 + a^4)}, \qquad
i(0,1) = \frac{ (-i2 + a^2 P)}{(-1 + a^4)} , $$
and
$$ C _1=  
i(1,0)-i(0,1) = \frac{ -2 a^2 P+ (1+a^4) i_2}{(-1+a^4)}.$$

\item For fourth moments we have the following:
\begin{eqnarray*}
i(2,0)+2\,i(1,1)+ i(0,2) 
&=& i_4, \\
\frac{i(2,0)}{a^2} +(a^2 +\frac{1}{a^2}) \,i(1,1)+ a^2\,{i(0,2)} 
&=& i_2 ,\\
\frac{i(2,0)}{a^4} +2\,i(1,1)+ a^4\,{i(0,2)} 
&=& P .
\end{eqnarray*}
A quantity occuring in the approximation $\Sigma_{a,1}$ is
$$C_2 = i(2,0)-2 i(1,1)+i(0,2) 
= \frac{ 2 a^4 P -4 a^2 (1+a^4) i_2+ (1+a^4)^2 i_4}{(-1+a^4)^2} . $$

\end{itemize}

The Cauchy-Schwarz inequality gives $C_2\ge{C_1^2/P}$.

\subsection{Further work to be done}

The further work will need more of the boundary moments $i_{2k}$.
We expect that there is a 3-term recurrence relation.

If (as we expect) all the $i_{2k}$ can be written in the form
$$ i_{2k} = r_K(k,a) {\rm EllipticK}(e) +r_E(k,a) {\rm EllipticE}(e) , $$
one can eliminate the EllipticK and EllipticE from pairs of such equations to get the
3-term equation
$$ i_{2k} = f_+(k,a) i_{2(k+1)} +  f_-(k,a) i_{2(k-1)} .$$
It remains to find the functions $f$, i.e. the recurrence.

It may also be possible to get to the recurrence from the integrals defining $i_{2k}$
as done for $g_k$ and ${\hat g}_k$.

The ellipse is given parametrically by
$$ x = a \cos(\psi), \qquad y =\frac{1}{a}\sin(\psi) . $$
Then
$$ r^2 = x^2+y^2 = \frac{1}{a^2} +(a^2 - \frac{1}{a^2}) \cos(\psi)^2 
= a^2 \left(  \frac{1}{a^4} + e^2  \cos(\psi)^2 \right) ,
$$
and
$$ (\frac{ds}{d\psi})^2
= a^2 \left(  1- e^2  \cos(\psi)^2 \right) .
$$
So, 
$$ \frac{ds}{d\psi} = a\, {\hat g}(\psi) , $$
so
$$ i_{2k}
= 4 a^{2k+1}\, \int_0^{\pi/2} \left(  \frac{1}{a^4} + e^2  \cos(\psi)^2 \right)^k\, {\hat g}(\psi)  .
$$
It should be possible to find a recurrence relation satisfied by the $i_{2k}$.

Of more immediate relevance
$$ C_1 = 4 a \, \int_0^{\pi/2} \left( (a^2+\frac{1}{a^2})\cos(\psi)^2 -a^2\right)\,  {\hat g}(\psi)  ,
$$
$$ C_2 = 4 a \, \int_0^{\pi/2} \left( (a^2+\frac{1}{a^2})^2\cos(\psi)^4 -2(1+a^4)\cos(\psi)^2 +a^4\right)\,  {\hat g}(\psi)  .
$$

One wonders also about `Pascal-triangle-like' organizations of the $i(m,n)$ quantities:
\begin{eqnarray*}
&i(0,0)&  \\
i(1,0)& \qquad& i(0,1) \\
i(2,0)\qquad&i(1,1)& \qquad i(0,2) \\
i(3,0) \qquad i(2,1)& \qquad&i(1,2)\qquad  i(0,3) 
\end{eqnarray*}

For the ellipse the first 3 rows of the above satisfy the determinantal equation:
\begin{equation}
\begin{vmatrix}
P& i_{xx}& i_{yy}\\
i_{xx}& i_{xxxx}& i_{xxyy}\\
i_{yy}& i_{xxyy}& i_{yyyy}
\end{vmatrix} 
= \begin{vmatrix}
i(0,0)& i(1,0)& i(0,1)\\
i(1,0)& i(2,0)& i(1,1)\\
i(0,1)& i(1,1)& i(0,2)
\end{vmatrix} = 0 .
\label{eq:iDet}
\end{equation}
The notation on the left is that in 
Part I.
The 3 by 3 matrix is nonnegative semidefinite as Cauchy-Schwarz gives
$$ i(1,0)^2 \le i(0,0) i(2,0), \  i(0,1)^2 \le i(0,0) i(0,2), \  i(1,1)^2 \le i(2,0) i(0,2) .$$

Returning to the `Pascal-triangle-like' figure, we note the following method of finding the $i(m,n)$.
As in the treatment of the Fourier series for $g$ and $\hat g$, set
$$ x=a\cos(\psi),\ \ y=\frac{\sin(\psi)}{a}, \ \ 
\frac{ds}{d\psi} = a\sqrt{1- e^2 \cos(\psi)^2} , \ \ e^2 = 1-\frac{1}{a^4} .$$
Using this
$$ i(m,n) = 2 a^{2(m-n)+1} j(m,n), $$
with
$$j(m,n) = 2\int_0^{\pi/2} (cos(\psi))^{2m} \, (1 -\cos(\psi)^2)^n\, \sqrt{1- e^2 \cos(\psi)^2} \, d\psi . $$
Setting $u=\cos(\psi)^2$
$$ j(m.n) = \int_0^1 u^m\, (1-u)^n \sqrt{\frac{1-e^2 u}{u(1-u)}} \, du .$$
It is trivial that
\begin{equation}
 j(m+1,n) + j(m,n+1) = j(m,n) .
 \label{eq:jtriv}
 \end{equation}
Each of the $j(m,n)$ has the form
$$ j(m,n) = f_K(m,n,e^2) {\rm EllipticK}(e) +
f_E(m,n,e^2) {\rm EllipticE}(e) , $$
for rational functions $f_E$ and $f_K$.
Defining
$$ u_m = j(m,0),\qquad v_n = j(0,n), $$
one can find the first couple of terms and recurrences to determine the rest.
This enables one to determine the entries at the ends of the rows in the triangle below
\begin{eqnarray*}
&j(0,0)&  \\
j(1,0)& \qquad& j(0,1) \\
j(2,0)\qquad&j(1,1)& \qquad j(0,2) \\
j(3,0) \qquad j(2,1)& \qquad&j(1,2)\qquad  j(0,3) 
\end{eqnarray*}
From the entries at the ends of the rows one can use equation~(\ref{eq:jtriv}) rearranged to
determine all the other entries.

The recurrences for the $u_n$ and $v_n$ are as follows.
$$ (2m-1) u_{m-1} - ((2m+2)e^2 +2m) u_m + (2m+3) e^2 u_{m+1} = 0 , $$
$$ -(2m - 1)\,(1 - e^2)\, v_{m - 1} -
2((2m + 1)e^2 - m)v_m + (2m + 3)e^2v_{ m + 1}=0 . $$
There are various identities that might be used to check, for example
(just from $\sin^2 +\cos^2=1$)
$$v_m = \sum_{j=0}^m (-1)^j{\rm Binomial}(m,j) u_j .$$





\end{document}